\documentclass[opre,nonblindrev]{informs3_modified1}

\DoubleSpacedXI 

\usepackage{endnotes}
\let\footnote=\endnote

\usepackage{natbib}
 \bibpunct[, ]{(}{)}{,}{a}{}{,}%
 \def\bibfont{\small}%
\TheoremsNumberedThrough     
\ECRepeatTheorems

\EquationsNumberedThrough    

\usepackage{bm}

\usepackage{caption,subcaption}
\usepackage{color}
\usepackage{algpseudocode}
\usepackage{algorithm}
\def\mbp{{\mathbf{p}}}

\begin{document}


\RUNAUTHOR{Ghosh and Lam}

\RUNTITLE{Robust Analysis in Stochastic Simulation}

\TITLE{Robust Analysis in Stochastic Simulation: Computation and Performance Guarantees}

\ARTICLEAUTHORS{%
\AUTHOR{Soumyadip Ghosh}
\AFF{IBM Research AI, IBM T.J. Watson Research Center, Yorktown Heights, NY 10598, \EMAIL{ghoshs@us.ibm.com}} 
\AUTHOR{Henry Lam}
\AFF{Department of Industrial Engineering and Operations Research, Columbia University, New York, NY 10027, \EMAIL{henry.lam@columbia.edu}}
} 

\ABSTRACT{%
Any performance analysis based on stochastic simulation is subject to the
errors inherent in misspecifying the modeling assumptions, particularly the input distributions. In situations with little support from data, we investigate the use of worst-case analysis to analyze these errors, by representing the partial, nonparametric knowledge of the input models via optimization constraints. We study the performance and robustness guarantees of this approach. We design and analyze a numerical scheme for solving a general class of simulation objectives and uncertainty specifications. The key steps involve a randomized discretization of the probability spaces, a simulable unbiased gradient estimator using a nonparametric analog of the likelihood ratio method, and a Frank-Wolfe (FW) variant of the stochastic approximation (SA) method (which we call FWSA) run on the space of input probability distributions. A convergence analysis for FWSA on non-convex problems is provided. We test the performance of our approach via several numerical examples.
}%



\maketitle

%


\section{Introduction}
Simulation-based performance analysis of stochastic models, or stochastic simulation, is
built on input model assumptions that to some extent deviate from the
truth. 
Consequently, a performance analysis subject to these input errors may lead to poor prediction and suboptimal decision-making.
To address this important problem, a typical framework in the stochastic simulation literature focuses on output variability measures or confidence bounds that account for the input uncertainty when input data are available. Established statistical techniques such as the bootstrap (e.g., \cite{barton1993uniform,barton2013quantifying}), goodness-of-fit tests (e.g., \cite{banks2000dm}), Bayesian inference and model selection (e.g., \cite{chick2001input,zouaoui2004accounting}) and the delta method (e.g., \cite{cheng1998two,cheng2004calculation}) have been proposed and have proven effective in many situations.

In this paper, we take a different approach for situations with insufficient data, or when the modeler wants to assess risk beyond what the data or the model indicates. Such situations can arise when the system, service target or operational policy in study is at a testing stage without much prior experience. 
To find reliable output estimates in these settings, we investigate a worst-case approach with respect to the input models. In this framework, the modeler represents the partial and nonparametric beliefs about the input models as constraints, and computes tight worst-case bounds among all models that satisfy them. More precisely, let $Z(P^1,\ldots,P^m)$ be a performance measure that depends on $m$ input models, each generated from a probability distribution $P^i$. The formulation for computing the worst-case bounds are
\begin{equation}
\min_{P^i\in\mathcal U^i, i=1,\ldots,m}Z(P^1,\ldots,P^m)\text{\ \ \ \ and\ \ \ \ }\max_{P^i\in\mathcal U^i, i=1,\ldots,m}Z(P^1,\ldots,P^m) \label{generic}
\end{equation}
The set $\mathcal U^i$ encodes the collection of all possible $P^i$ from the knowledge of the modeler. The decision variables in the optimizations in \eqref{generic} are the unknown models $P^i,i=1,\ldots,m$.

The primary motivation for using \eqref{generic} is the robustness against model misspecification, where a proper construction of the set $\mathcal U^i$ avoids making specific assumptions beyond the modeler's knowledge. The following three examples motivate and explain further.

\begin{example}[Robust bounds under expert opinion]
When little information is available for an input model, a common practice in stochastic simulation is to summarize its range (say $[a,b]$) and mean (or mode) as a triangular distribution, where the base of the triangle denotes the range and the position of the peak is calibrated from the mean. This specific distribution only crudely describes the knowledge of the modeler and may deviate from the true distribution, even if $a,b,\mu$ are correctly specified. Instead, using
\begin{equation}
\mathcal U^i=\{P^i:E_{P^i}[X^i]=\mu,\ \text{supp\ }P^i=[a,b]\}\label{triangle}
\end{equation}
in formulation \eqref{generic}, where $X^i$ is the random variate, $E_{P^i}[\cdot]$ is the expectation under $P^i$, and $\text{supp\ }P^i$ is the support of $P^i$, will give a valid interval that covers the true performance measure whenever $a,b,\mu$ are correctly specified. Moreover, when these parameters are not fully known but instead specified within a range, \eqref{triangle} can be relaxed to
$$\mathcal U^i=\{P^i:\underline\mu\leq E_{P^i}[X^i]\leq\overline\mu,\ \text{supp\ }X^i=[\underline a,\overline b]\}$$
where $[\underline\mu,\overline\mu]$ denotes the range of the mean and $\underline a,\overline b$ denote the lower estimate of the lower support end and upper estimate of the upper support end respectively. The resulting bound will cover the truth as long as these ranges are supplied correctly.\label{triangle ex}
\Halmos
\end{example}

\begin{example}[Dependency modeling]
In constructing dependent input models, common approaches in the simulation literature fit the marginal description and the correlation of a multivariate model to a specified family. Examples include Gaussian copula (e.g., \cite{lurie1998approximate,channouf2009fitting}; also known as normal-to-anything (NORTA), e.g. \cite{cario1997modeling}) and chessboard distribution (\cite{ghosh2002chessboard}) that uses a domain discretization. These distributions are correctly constructed up to their marginal description and correlation, provided that these information are correctly specified. However, dependency structure beyond correlation can imply errors on these approaches (e.g., \cite{lam2016serial}), and formulation \eqref{generic} can be used to get bounds that address such dependency. For example, suppose $P^i$ is a bivariate input model with marginal distributions $P^{i,1},P^{i,2}$, marginal means $\mu^{i,1},\mu^{i,2}$ and covariance $\rho^i$. We can set
$$\mathcal U^i=\{P^i:P_{P^{i,1}}(X^{i,1}\leq q_j^{i,1})=\nu_j^1,j=1,\ldots,l_1,\ P_{P^{i,2}}(X^{i,2}\leq q_j^{i,2})=\nu_j^2,j=1,\ldots,l_2,\ E[X^{i,1}X^{i,2}]=\rho^i+\mu^{i,1}\mu^{i,2}\}$$
where $(X^{i,1},X^{i,2})$ denote the random vector under $P^i$, and $q_j^{i,1},q_j^{i,2},\nu_j^{i,1},\nu_j^{i,2}$ are pre-specified quantiles and probabilities of the respective marginal distributions. Unlike previous approaches, \eqref{generic} outputs correct bounds on the truth given correctly specified marginal quantiles and correlation, regardless of the dependency structure. \label{dependency ex}
\Halmos
\end{example}

\begin{example}[Model risk]
Model risk refers broadly to the uncertainty in analysis arising from the adopted model not being fully accurate. This inaccuracy occurs as the adopted model (often known as the baseline model), typically obtained from the best statistical fit or expert opinion, deviates from the truth due to the real-world non-stationarity and the lack of full modeling knowledge or capability. To assess model risk, a recently surging literature studies the use of statistical distance as a measurement of model discrepancy (e.g., \cite{gx12b,lam2013robust}). Given the baseline model $P_b^i$, the idea is to represent the uncertainty in terms of the distance away from the baseline via a neighborhood ball
\begin{equation}
\mathcal U^i=\{P^i:d(P^i,P_b^i)\leq\eta^i\}\label{ball}
\end{equation}
where $d$ is a distance defined on the nonparametric space of distributions (i.e., without restricting to any parametric families). The bounds drawn from formulation \eqref{generic} assess the effect of model risk due to the input models, tuned by the ball size parameter $\eta^i$ that denotes the uncertainty level. Besides risk assessment, this approach can also be used to obtain consistent confidence bounds for the true performance measure, when $P_b^i$ is taken as the empirical distribution and $\eta$ and $d$ are chosen suitably (discussed further in Section \ref{sec:discretization}).
\label{model ex}
\Halmos
\end{example}


Our worst-case approach is inspired from the literature of robust optimization (\cite{ben2009robust,bertsimas2011theory}), which considers decision-making under uncertainty and advocates optimizing decisions over worst-case scenarios. In particular, when the uncertainty lies in the probability distributions that govern a stochastic problem, the decision is made to optimize under the worst-case distributions, a class of problems known as distributionally robust optimization (e.g. \cite{delage2010distributionally,lim2006model}). Such an approach has also appeared in so-called robust simulation or robust Monte Carlo in the simulation literature (\cite{hu2012robust,gx12b}). However, the methodologies presented in the above literature focus on structured problems where the objective function is tractable, such as linear or linearly decomposable. In contrast, $Z(\cdot)$ for most problems in stochastic simulation is nonlinear and unstructured, obstructing the direct adaptation of the existing methods. In view of this, our main objective is to design an efficient simulation-based method to compute the worst-case bounds for formulation \eqref{generic} that can be applied to broad classes of simulation models and input uncertainty representations.

\subsection{Our Contributions}
We study a simulation-based iterative procedure for the worst-case optimizations \eqref{generic}, based on a modified version of the
celebrated stochastic approximation (SA) method (e.g. \cite{kushner}). Because of the iterative nature, it is difficult to directly operate on the space of continuous distributions except in very special cases. Thus, our first contribution (Section \ref{sec:discretization}) is to provide a randomized discretization scheme that can provably approximate the continuous counterpart. This allows one to focus on discrete distributions on fixed support points as the decision variable to feed into our SA algorithm. 


We develop the SA method in several aspects. In Section \ref{sec:gradient}, we construct an unbiased gradient estimator for $Z$ based on the idea of the Gateaux derivative for functionals of probability distributions (\cite{serfling2009approximation}), which is used to obtain the direction in each subsequent SA iterate. The need for such a construction is motivated by the difficulty in na\"{i}ve implementation of standard gradient estimators: An arbitrary perturbation of a probability
distribution, which is the decision variable in the optimization, may shoot outside the probability simplex and results in a gradient that does not bear any probabilistic meaning and subsequently does not support simulation-based estimation. Our approach effectively restricts the direction of perturbation to points within the probability simplex, leading to a simulable gradient estimator. We justify our approach as a nonparametric version of the classical
likelihood ratio method (or the score function method)
(\cite{glynn1990likelihood,reiman1989sensitivity,rubinstein1986score}).


Next, in Sections \ref{sec:procedure} and \ref{sec:guarantees}, we design and analyze our SA scheme under the uncertainty constraints. We choose to use a stochastic counterpart of the so-called Frank-Wolfe (FW) method (\cite{frank1956algorithm}), known synonymously as the conditional gradient method in deterministic nonlinear programming. For convenience we call our scheme FWSA. Note that a standard SA iteration follows the estimated gradient up to a
pre-specified step size to find the next candidate iterate. When the
formulation includes constraints, the common approach in the SA literature projects the candidate solution onto the feasible region in order to define the next
iterate (e.g. \cite{kushner}). Instead, our method looks in advance for a feasible direction along which the next iterate is guaranteed to lie in the (convex) feasible region. In order to find
this feasible direction, an optimization subproblem with a linear
objective function is solved in
each iteration. We base our choice of using FWSA on its computational benefit in solving these subproblems, as their linear objectives allow efficient solution scheme for high-dimensional decision variables for many choices of the set $\mathcal U^i$.


We characterize
the convergence rate of FWSA in terms of the step size and the
number of simulation replications used to estimate the gradient at each iteration. The form of our convergence bounds suggests
prescriptions for the step-size and sample-size sequences that are efficient with respect to the cumulative
number of sample paths simulated to generate all
the gradients until the current iterate. The literature on the stochastic FW methods for non-convex problems is small. \cite{kushner1974stochastic} proves almost sure convergence under assumptions that can prescribe algorithmic specifications only for one-dimensional settings. During the review process of this paper, two other convergence rate studies \cite{reddi2016stochastic} and \cite{lafond2016} have appeared. Both of them assume the so-called $G$-Lipschitz condition on the gradient estimator that does not apply to our setting. Consequently, our obtained convergence rates are generally inferior to their results. Nonetheless, we will point out how our rates almost match theirs under stronger assumptions on the behavior of the iterates that we will discuss. 


Finally, in Section \ref{sec:numerics} we provide numerical validation of our approach
using two sets of experiments, one testing the performance of our proposed randomized discretization strategy, and one on the convergence of FWSA.

\subsection{Literature Review}
We briefly survey three lines of related work. First, our paper is related to the literature on input model uncertainty. In the parametric regime, studies have focused on the construction of confidence intervals or variance decompositions to account for both parameter and stochastic uncertainty using data, via for instance the delta method (\cite{cheng1998two,cheng2004calculation}), the bootstrap (\cite{barton2013quantifying,cheng1997sensitivity}), Bayesian approaches (\cite{zouaoui2003accounting,xie2014bayesian,saltelli2010variance,saltelli2008global}), and metamodel-assisted analysis (\cite{xie2014bayesian,xie2013statistical}). Model selection beyond a single parametric model can be handled through goodness-of-fit or Bayesian model selection and averaging (\cite{chick2001input,zouaoui2004accounting}). Fully nonparametric approaches using the bootstrap have also been investigated (\cite{barton1993uniform,barton2001resampling,song2015quickly}).

Second, formulation \eqref{generic} relates to the literature on robust stochastic control (\cite{pjd00,iyengar2005robust,nilim2005robust,doi:10.1287/moor.1120.0540}) and distributionally robust optimization (\cite{delage2010distributionally,goh2010distributionally,ben2013robust,wiesemann2014distributionally}), where the focus is to make decision rules under stochastic environments that are robust against the ambiguity of the underlying probability distributions. This is usually cast in the form of a minimax problem where the inner maximization is over the space of distributions. This idea has spanned across multiple areas like economics (\cite{hsAER01,hansen2008robustness}), finance (\cite{gx12a,lsw11}), queueing (\cite{bertsimas2007semidefinite,jls10}), dynamic pricing (\cite{ls07}), inventory management (\cite{xin2015distributionally}), physical sciences (\cite{dupuis2016path}), and more recently machine learning (\cite{shafieezadeh2015distributionally,blanchet2016robust}). In the simulation context, \cite{hu2012robust} compared different global warming policies using Gaussian models with uncertain mean and covariance information. \cite{gx12b,glasserman2016bounding} studied approaches based on sample average approximation for solving distance-based constrained optimizations to quantify model risk in finance. \cite{lam2013robust,lam2016serial} investigated infinitesimal approximations for related optimizations to quantify model errors arising from sequences of uncertain input variates. \cite{bandi2012tractable} studied the view of deterministic robust optimization to compute various stochastic quantities. Simulation optimization under input uncertainty has also been studied via the robust optimization framework (\cite{fan2013robust,ryzhov2012ranking}), and the closely related approach using risk measures (\cite{qian2015composite,zhou2015simulation}). Lastly, optimizations over probability distributions have also arisen as generalized moment problems, applied to decision analysis (\cite{smith95,smith1993moment,bertsimas2005optimal}) and stochastic programming (\cite{birge1987computing}).  



Our algorithm relates to the literature on the FW method (\cite{frank1956algorithm}) and constrained SA. The former is a nonlinear programming technique initially proposed for convex optimization, based on sequential linearization of the objective function using the gradient at the solution iterate. The classical work of \cite{canon1968tight}, \cite{dunn1979rates} and \cite{dunn1980convergence} analyzed convergence properties of FW for deterministic convex programs. More recently, \cite{jaggi2013revisiting}, \cite{freund2014new} and \cite{hazan2016variance} carried out finite-time analysis for the FW method motivated by machine learning applications. For stochastic FW on non-convex problems (viewed as a class of constrained SA), \cite{kushner1974stochastic} focused on almost sure convergence based on a set of assumptions about the probabilistic behavior of the iterations, which were then used to tune the algorithm for one-dimensional problems. While this paper was under review, \cite{reddi2016stochastic} provided a complexity analysis in terms of the sample size in estimating gradients and the number of calls of the linear optimization routine. \cite{lafond2016} studied the performance in terms of regret in an online setting. Both \cite{reddi2016stochastic} and \cite{lafond2016} relied on the $G$-Lipschitz condition that our gradient estimator violated. Other types of constrained SA schemes include the Lagrangian method (\cite{buche2002rate}) and mirror descent SA (\cite{nemirovski2009robust}). Lastly, general convergence results for SA can be found in \cite{fu1994optimization}, \cite{kushner} and \cite{pasupathy2011stochastic}. 

\section{Formulation and Assumptions}\label{sec:formulation}
We focus on $Z(P^1,\ldots,P^m)$ that is a finite horizon performance measure generated from i.i.d. replications from the independent input models $P^1,\ldots,P^m$. Let $\mathbf X^i=(X_t^i)_{t=1,\ldots,T^i}$ be $T^i$ i.i.d. random variables on the space $\mathcal X^i\subset\mathbb R^{v^i}$, each generated under $P^i$. The performance measure can be written as
\begin{equation}
Z(P^1,\ldots,P^m)=E_{P^1,\ldots,P^m}[h(\mathbf X^1,\ldots,\mathbf X^m)]=\int\cdots\int h(\mathbf x^1,\ldots,\mathbf x^m)\prod_{t=1}^{T^1}dP(x_t^1)\cdots\prod_{t=1}^{T^m}dP(x_t^m)\label{perf measure}
\end{equation}
where $h(\cdot):\prod_{i=1}^m(\mathcal X^i)^{T^i}\to\mathbb R$ is a cost function, and $E_{P^1,\ldots,P^m}[\cdot]$ denotes the expectation associated with the generation of the i.i.d. replications. We assume that $h(\cdot)$ can be evaluated by the computer given the inputs. In other words, the performance measure \eqref{perf measure} can be approximated by running simulation.

\eqref{perf measure} is the stylized representation for transient performance measures in discrete-event simulation. For example, $\mathbf X^1$ and $\mathbf X^2$ can be the sequences of interarrival and service times in a queue, and $P^1$ and $P^2$ are the interarrival time and service time distributions. When $h(\mathbf X^1,\mathbf X^2)$ is the indicator function of the waiting time exceeding a threshold, \eqref{perf measure} will denote the corresponding threshold exceedance probability.

Next we discuss the constraints in \eqref{generic}. Following the terminology in robust optimization, we call $\mathcal U^i$ the \emph{uncertainty set} for the $i$-th input model. Motivated by the examples in the Introduction, we focus on two types of convex uncertainty sets:
\begin{enumerate}
\item\emph{Moment and support constraints: }We consider
\begin{equation}
\mathcal U^i=\{P^i:E_{P^i}[f_l^i(X^i)]\leq\mu_l^i,l=1,\ldots,s^i,\ \text{supp\ }P^i=A^i\}\label{moment set}
\end{equation}
where $X^i$ is a generic random variable under distribution $P^i$, $f_l^i:\mathcal X^i\to\mathbb R$, and $A^i\subset\mathcal X^i$. For instance, when $\mathcal X^i=\mathbb R$, $f_l^i(x)$ being $x$ or $x^2$ denotes the first two moments. When $\mathcal X^i=\mathbb R^2$, $f_l^i(x_1,x_2)=x_1x_2$ denotes the cross-moment. Equalities can also be represented via \eqref{moment set} by including $E_{P^i}[-f_l^i(X^i)]\leq-\mu_l^i$. Thus the uncertainty set \eqref{moment set} covers Examples \ref{triangle ex} and \ref{dependency ex} in the Introduction.

    Furthermore, the neighborhood measured by certain types of statistical distance (Example \ref{model ex}) can also be cast as \eqref{moment set}. For instance, suppose $d$ is induced by the sup-norm on the distribution function on $\mathbb R$. Suppose $P^i$ is a continuous distribution and the baseline distribution $P_b^i$ is discrete with support points $y_j,j=1,\ldots,n^i$. The constraint
    \begin{equation}
    \sup_{x\in\mathbb R}|F^i(x)-F_b^i(x)|\leq\eta^i\label{KS}
    \end{equation}
    where $F^i$ and $F_b^i$ denote the distribution functions for $P^i$ and $P_b^i$ respectively, can be reformulated as
    $$F_b^i(y_j+)-\eta^i\leq F^i(y_j)\leq F_b^i(y_j-)+\eta^i,\ j=1,\ldots,n^i$$
    where $F_b^i(y_j-)$ and $F_b^i(y_j+)$ denote the left and right limits of $F_b^i$ at $y_j$, by using the monotonicity of distribution functions. Thus
    $$\mathcal U^i=\{P^i:F_b^i(y_j+)-\eta^i\leq E^i[I(X^i\leq y_j)]\leq F_b^i(y_j-)+\eta^i,\ j=1,\ldots,n^i,\ \text{supp\ }P^i=\mathbb R\}$$
    where $I(\cdot)$ denotes the indicator function, falls into the form of \eqref{moment set}. \cite{bertsimas2014robust} considers this reformulation for constructing uncertainty sets for stochastic optimization problems, and suggests to select $\eta^i$ as the quantile of the Kolmogorov-Smirnov statistic if $F_b^i$ is the empirical distribution function constructed from continuous i.i.d. data.



 \item\emph{Neighborhood of a baseline model measured by $\phi$-divergence: }Consider
 \begin{equation}
 \mathcal U^i=\{P^i:d_\phi(P^i,P_b^i)\leq\eta^i\}\label{phi set}
 \end{equation}
 where $d_\phi(P^i,P_b^i)$ denotes the $\phi$-divergence from a baseline distribution $P_b^i$ given by
 $$d_\phi(P^i,P_b^i)=\int\phi\left(\frac{dP^i}{dP_b^i}\right)dP_b^i$$
 which is finite only when $P^i$ is absolutely continuous with respect to $P_b^i$. The function $\phi$ is a convex function satisfying $\phi(1)=0$. This family covers many widely used distances. Common examples are $\phi(x)=x\log x-x+1$ giving the KL divergence, $\phi(x)=(x-1)^2$ giving the (modified) $\chi^2$-distance, and $\phi(x)=(1-\theta+\theta x-x^\theta)/(\theta(1-\theta)),\ \theta\neq0,1$ giving the Cressie-Read divergence. Details of $\phi$-divergence can be found in, e.g.,  \cite{pardo2005statistical,ben2013robust,bayraksan2015data}. 


\end{enumerate}

As precursed in the Introduction, in the context of simulation analysis where $(P^1,\ldots,P^m)$ are the input models, $Z(\cdot)$ in \eqref{perf measure} is in general a complex nonlinear function. This raises challenges in solving \eqref{generic} beyond the literature of robust control and optimization that considers typically more tractable objectives. Indeed, if $Z(\cdot)$ is a linear function in $P^i$'s, then optimizing over the two types of uncertainty sets above can both be cast as specialized classes of convex programs that can be efficiently solved. But linear $Z(\cdot)$ is too restrictive to describe the input-output relation in simulation. To handle a broader class of $Z(\cdot)$ and to address its simulation-based nature, we propose to use a stochastic iterative method. 
The next sections will discuss our methodology in relation to the performance guarantees provided by \eqref{generic}.

\section{Performance Guarantees and Discretization Strategy}\label{sec:discretization}
This section describes the guarantees provided by our framework. Section \ref{sec:randomized} first presents the motivation and justification of a discretization scheme for continuous input distributions. Section \ref{sec:stat} then discusses the statistical implications in more details.

\subsection{Randomized Discretization}\label{sec:randomized}
Suppose there is a ``ground true" distribution $P_0^i$ for each input model. Let $Z_*$ and $Z^*$ be the minimum and maximum values of the worst-case optimizations \eqref{generic}. Let $Z_0$ be the true performance measure, i.e. $Z_0=Z(P_0^1,\ldots,P_0^m)$. The following highlights an immediate implication of using \eqref{generic}:
\begin{proposition}
If $P_0^i\in\mathcal U^i$ for all $i$, then $Z_*\leq Z_0\leq Z^*$.\label{basic guarantee}
\end{proposition}

In other words, the bounds from the worst-case optimizations form an interval that covers the true performance measure if the uncertainty sets contain the true distributions.

We discuss a discretization strategy for the worst-case optimizations for continuous input distributions. We will show that, by replacing the continuous distribution with a discrete distribution on support points that are initially sampled from some suitably chosen distribution, we can recover the guarantee in Proposition \ref{basic guarantee} up to a small error. The motivation for using discretization comes from the challenges in handling decision variables in the form of continuous distributions when running our iterative optimization scheme proposed later.

We focus on the two uncertainty sets \eqref{moment set} and \eqref{phi set}. The following states our guarantee:
\begin{theorem}
Consider $Z(P^1,\ldots,P^m)$ in \eqref{perf measure}. Assume $h$ is bounded a.s.. 
Let $n^i,i=1,\ldots,m$ and $n$ be positive integers such that $n^i=nw^i$ for some fixed $w^i>0$, for all $i$. For each input model $i$, we sample $n^i$ i.i.d. observations $\{y_1^i,\ldots,y_{n^i}^i\}$ from a distribution $Q^i$ such that the true distribution $P_0^i$ is absolutely continuous with respect to $Q^i$, with $L^i=dP_0^i/dQ^i$ satisfying $\|L^i\|_\infty<\infty$, where $\|L^i\|_\infty$ denotes the essential supremum of $L^i$ under $Q^i$. Consider the optimizations
\begin{equation}
\hat Z_*=\min_{P^i\in\hat{\mathcal U}^i, i=1,\ldots,m}Z(P^1,\ldots,P^m)\text{\ \ \ \ and\ \ \ \ }\hat Z^*=\max_{P^i\in\hat{\mathcal U}^i, i=1,\ldots,m}Z(P^1,\ldots,P^m)\label{sample counterpart}
\end{equation}
where each $\hat{\mathcal U}^i$ contains discrete distributions supported on $\{y_1^i,\ldots,y_{n^i}^i\}$, defined in one of the two cases below. For each case, we also make additional assumptions as follows:
\begin{enumerate}
\item
Set
\begin{equation}
\hat{\mathcal U}^i=\{P^i:E_{P^i}[f_l^i(X^i)]\leq\mu_l^i,l=1,\ldots,s^i,\ \text{supp\ }P^i\subset\{y_1^i,\ldots,y_{n^i}^i\}\}\label{moment discrete}
\end{equation}
Moreover, assume that $P_0^i$ satisfies $E_{P_0^i}|f_l^i(X^i)|<\infty$ and $E_{P_0^i}[f_l^i(X^i)]<\mu_l^i$ for all $l=1,\ldots,s^i$.


\item 
The distribution $Q^i$ is chosen such that $P_b^i$ is absolutely continuous with respect to $Q^i$, and we denote $L_b^i=dP_b^i/dQ^i$. Set
\begin{equation}
\hat{\mathcal U}^i=\{P^i:d_\phi(P^i,\hat P_b^i)\leq\eta^i\}\label{phi discrete}
\end{equation}
where $\hat P_b^i$ is defined as
$$\hat P_b^i=\sum_{j=1}^{n^i}\frac{L_b^i(y_j^i)}{\sum_{r=1}^{n^i}L_b^i(y_r^i)}\delta(y_j^i)$$
with $\delta(y)$ denoting the delta measure at $y$.
Moreover, assume $P_0^i$ satisfies $E_{P_b^i}|\phi(dP_0^i/dP_b^i)|<\infty$ and $d_\phi(P_0^i,P_b^i)<\eta^i$. Additionally, assume $\phi(\cdot)$ satisfies the continuity condition $|\phi(t(1+\lambda))-\phi(t)|\leq|\phi(t)|\kappa_1(\lambda)+\kappa_2(\lambda)$ for any $t\geq0$ and $\lambda$ in a fixed neighborhood of 0, where $\kappa_1(\cdot)$ and $\kappa_2(\cdot)$ are two functions such that $\kappa_1(\lambda)=O(\lambda)$ and $\kappa_2(\lambda)=O(\lambda)$ as $\lambda\to0$.
\end{enumerate}

Then we have
\begin{equation}
\hat Z_*\leq Z_0+O_p\left(\frac{1}{\sqrt n}\right)\leq\hat Z^*\label{approx bound}
\end{equation}\label{sample thm}
\end{theorem}

Here $O_p(1/\sqrt n)$ is an error term $e_n$ that is of stochastic order $1/\sqrt n$, i.e., for any $0<\epsilon<1$, there exist $M,N>0$ such that $P(|\sqrt ne_n|<M)>1-\epsilon$ for any $n>N$. Theorem \ref{sample thm} is proved in Appendix \ref{sec:proofs}. We have a few immediate remarks:
\begin{enumerate}
\item
Optimizations \eqref{sample counterpart} are the sample counterparts of the original worst-case optimizations \eqref{generic} with uncertainty sets given by \eqref{moment set} or \eqref{phi set}, which optimize discrete distributions over support points that are sampled from generating distributions $Q^i$'s. Theorem \ref{sample thm} guarantees that, if the original worst-case optimizations give valid covering bounds for the true performance measure (in the spirit of Proposition \ref{basic guarantee}), then so are the sample counterparts, up to an error $O_p(1/\sqrt n)$ where $n$ denotes the order of the sample size used to construct the sets of support points. The constant implicit in this $O_p(1/\sqrt n)$ error depends on the sensitivity of $Z$ with respect to the input distributions, as well as the discrepancies between the true input distributions and the support-generating distributions.

\item The condition $\|L^i\|_\infty<\infty$ implies that $Q^i$ has a tail at least as heavy as $P_0^i$. In practice, the tail of the true distribution $P_0^i$ is not exactly known a priori. This means that it is safer to sample the support points from a heavy-tailed distribution. Additionally, in the case of $\phi$-divergence, the generating distribution should also support the baseline. One easy choice is to merely use the baseline as the generating distribution.

\item The conditions $E_{P_0^i}[f_l^i(X^i)]<\mu_l^i$ and $d_\phi(P_0^i,P_b^i)<\eta^i$ state that $E_{P_0^i}[f_l^i(X^i)]$ and $d_\phi(P_0^i,P_b^i)$ are in the interior of $\{(z_1,\ldots,z_{s^i}):z_l\leq\mu_l^i,\ l=1,\ldots,s^i\}$ and $\{z:z\leq\eta^i\}$ respectively. These conditions guarantee that $P_0^i$ projected on a sample approximation of the support is asymptotically feasible for \eqref{sample counterpart}, which helps lead to the guarantee \eqref{approx bound}. In general, the closer $P_0^i$ is to the boundary of the uncertainty set, i.e., the smaller the values of $\mu_l^i-E_{P_0^i}[f_l^i(X^i)]$ and $\eta^i-d_\phi(P_0^i,P_b^i)$, the larger the sample size is needed for the asymptotic behavior in \eqref{approx bound} to kick in, a fact that is not revealed explicitly in Theorem \ref{sample thm}. One way to control this required sample size is to expand the uncertainty set by a small margin, say $\epsilon>0$, i.e., use $E_{P^i}[f_l^i(X^i)]\leq\mu_l^i+\epsilon$ and $d_\phi(P^i,P_b^i)\leq\eta^i+\epsilon$, in \eqref{moment discrete} and \eqref{phi discrete}. Note that, in the case of moment equality constraint, say $E_{P^i}[f_l^i(X^i)]=\mu_l^i$, one does have to deliberately relax the constraint to $\mu_l^i-\epsilon\leq E_{P^i}[f_l^i(X^i)]\leq\mu_l^i+\epsilon$ for the interior-point conditions to hold.

\item The continuity assumption imposed on $\phi(\cdot)$ in Case 2 is satisfied by many common choices, including KL, (modified) $\chi^2$-distance, and Burg entropy (see the definitions in \cite{ben2013robust}).

\item As $n^i$ increases, the sampled uncertainty set $\hat{\mathcal U}^i$ enlarges as it contains distributions supported on more values. As a result, $\hat Z_*$ becomes smaller and $\hat Z^*$ larger as $n^i$ increases. Moreover, since $\hat{\mathcal U}^i\subset\mathcal U^i$, we have $\hat Z_*\geq Z_*$ and $\hat Z^*\leq Z^*$. This means that as the generated support size increases, the interval $[\hat Z_*,\hat Z^*]$ progressively widens and is always contained by the interval $[Z_*,Z^*]$.

\end{enumerate}

\subsection{Statistical Implications}\label{sec:stat}
We further discuss the statistical guarantees implied from Section \ref{sec:randomized}. First, a probabilistic analog of Proposition \ref{basic guarantee} is:
\begin{proposition}
Suppose $\mathcal U^i$ contains the true distribution $P_0^i$ for all $i$ with confidence $1-\alpha$, i.e. $\mathbb P(\mathcal U^i\ni P_0^i\text{\ for all\ }i=1,\ldots,m)\geq1-\alpha$, then $\mathbb P(Z_*\leq Z_0\leq Z^*)\geq1-\alpha$, where $\mathbb P$ denotes the probability generated from a combination of data and prior belief.\label{simple conf}
\end{proposition}

Proposition \ref{simple conf} follows immediately from Proposition \ref{basic guarantee}. In the frequentist framework, $\mathbb P$ refers to the probability generated from data. However, Proposition \ref{simple conf} can also be cast in a Bayesian framework, in which $\mathbb P$ can represent the prior (e.g., from expert opinion) or the posterior belief.

Proposition \ref{simple conf} reconciles with the established framework in distributionally robust optimization that the uncertainty set $\mathcal U^i$ should be chosen as a confidence set for the true distribution, in order to provide a guarantee for the coverage probability on the true objective, in the case that $\mathbb P$ represents the generation of data under a true model. Some strategies for constructing confidence sets are:
\begin{enumerate}
\item For moment constraint $E_{P^i}[f_l^i(X^i)]\leq\mu_l^i$, one can choose $\mu_l^i$ as the upper confidence bound of the moment.
\item For the sup-norm constraint in \eqref{KS}, supposing that $P^i$ is continuous, $\eta^i$ chosen as the $(1-\alpha)$-quantile of $\sup_{x\in[0,1]}B(t)/\sqrt{n^i}$, where $B(t)$ is a standard Brownian bridge, gives an approximate $(1-\alpha)$ confidence region. This follows from the limiting distribution of the Kolmogorov-Smirnov statistic (see, e.g., \cite{bertsimas2014robust}). This calibration becomes conservative (but still correct) when $P^i$ is discrete, and one could use the bootstrap as a remedy. Note that the Kolmogorov-Smirnov-based confidence region is crude for the tail in that it can include a wide range of tail behaviors, and thus is not recommended if the performance measure of interest is sensitive to the tail.
\item For the $\phi$-divergence-based constraint in \eqref{phi set}, under the assumption that $P^i$ has finite support of size $r^i$, \cite{ben2013robust} proposes using $\eta^i=(\phi''(1)/(2n^i))\chi^2_{r^i-1,1-\alpha}$ in the case $P_b^i$ is taken as the empirical distribution, where $\chi^2_{r^i-1,1-\alpha}$ is the $(1-\alpha)$-quantile of a $\chi^2$-distribution with degree of freedom $r^i-1$. This leads to an approximate $(1-\alpha)$ confidence region by using the asymptotics of goodness-of-fit statistics (\cite{pardo2005statistical}). The resulting region from this approach, however, can be conservative as the involved degree of freedom can be large. Recent works such as \cite{lam2015quantifying,duchistatistics,lam2016recovering} investigate the tightening of divergence-based regions and extend their use to continuous data using the empirical likelihood theory. This theory can also potentially shed insights on the (second-order) accuracies achieved using different divergences (\cite{owen2001empirical}). Other alternatives include using the Wasserstein distance; see, e.g., \cite{esfahani2015data,blanchet2016quantifying,gao2016distributionally} for these developments and the involved ball-size calibration methods.
\end{enumerate}

When discretization is applied, the probabilistic analog of Theorem \ref{sample thm} is:
\begin{theorem}
Suppose all assumptions in Theorem \ref{sample thm} are in place except that $E_{P_0^i}[f_l^i(X^i)]<\mu_l^i$ or $d_\phi(P_0^i,P_b^i)<\eta^i$ now holds true jointly for all $i$ with confidence $1-\alpha$ under $\mathbb P$. Then $\mathbb P(\hat Z_*\leq Z_0+O_p(1/\sqrt n)\leq\hat Z^*)\geq1-\alpha$.\label{sample conf}
\end{theorem}

Theorem \ref{sample conf} follows immediately from Theorem \ref{sample thm}. Like before, Theorem \ref{sample conf} translates \eqref{generic}, whose input models can be continuously represented, to \eqref{sample counterpart} that is imposed over discrete distributions, by paying a small price of error. In the next section we discuss our algorithm run over discrete distributions and point out clearly why the discretization is necessary when the input distributions are continuous.

We close this section with two cautionary remarks. First, while our discretization strategy works for problems involving independent low-dimensional input distributions (which occur often in stochastic simulation), high-dimensional joint dependent models may greatly inflate the constant implicit in the error term, and we do not advise using our strategy in such settings. Second, in general, the finer the discretization scale (i.e., the more generated support points), the higher is the decision space dimension for the resulting optimization problem, and there is a tradeoff on the discretization scale between the approximation error and the optimization effort. Obviously, when the input model is finite discrete, the sampling step depicted in Theorems \ref{sample thm} and \ref{sample conf} is unnecessary, and our subsequent results regarding the algorithm applies readily to this case.

\section{Gradient Estimation on Probability Simplices via a Nonparametric Likelihood Ratio Method}\label{sec:gradient}
Since we work in the discrete space, for simplicity we denote $\mathbf p^i=(p_j^i)_{j=1,\ldots,n^i}\in\mathbb R^{n^i}$ as the vector of probability weights for the discretized input model $i$. This probability vector is understood to apply on the support points $\{y_1^i,\ldots,y_{n^i}^i\}$. Moreover, let $\mathbf p=\text{vec}(\mathbf p^i:i=1,\ldots,m)\in\mathbb R^N$ where $\text{vec}$ denotes a concatenation of the vectors $\mathbf p^i$'s as a single vector, and $N=\sum_{i=1}^mn^i$. We denote $\mathcal P_l=\{(p_1,\ldots,p_l)\in\mathbb R^l:\sum_{j=1}^lp_j=1, p_j\geq0, j=1,\ldots,l\}$ as the $l$-dimensional probability simplex. Hence $\mathbf p^i\in\mathcal P_{n^i}$. For convenience, let $\mathcal P=\prod_{i=1}^m\mathcal P_{n^i}$, so that $\mathbf p\in\mathcal P$. The performance measure in \eqref{sample counterpart} can be written as $Z(\mathbf p)$. Furthermore, denote $T=\max_{i=1,\ldots,m}T^i$ as the maximum length of replications among all input models. We also write $\mathbf X=(\mathbf X^1,\ldots,\mathbf X^m)$ and $h(\mathbf X)=h(\mathbf X^1,\ldots,\mathbf X^m)$ for simplicity.
Recall that $I(E)$ denotes the indicator function for the event $E$. In the rest of this paper, $'$ denotes transpose, and $\|\mathbf x\|$ denotes the Euclidean norm of a vector $\mathbf x$. We also write $Var_{\mathbf p}(\cdot)$ as the variance under the input distribution $\mathbf p$. Inequalities for vectors are defined component-wise.

We shall present an iterative simulation-based scheme for optimizing \eqref{sample counterpart}. The first step is to design a method to extract the gradient information of $Z(\mathbf p)$. Note that the standard gradient of $Z(\mathbf p)$, which we denote as $\nabla Z(\mathbf p)$, obtained through differentiation of $Z(\mathbf p)$, may not lead to any simulable object. This is because an arbitrary perturbation of $\mathbf p$ may shoot out from the set of probability simplices, and the resulting gradient will be a high-degree polynomial in $\mathbf p$ that may have no probabilistic interpretation and thus is not amenable to simulation-based estimation.

We address this issue by considering the set of perturbations within the simplices. Our approach resembles the Gateaux derivative on a functional of probability distribution (\cite{serfling2009approximation}) as follows. Given any $\mathbf p^i$, define a mixture distribution $(1-\epsilon)\mathbf p^i+\epsilon\mathbf 1_j^i$, where $\mathbf 1_j^i$ represents a point mass on $y_j^i$, i.e. $\mathbf 1_j^i=(0,0,\ldots,1,\ldots,0)\in\mathcal P_{n^i}$ and 1 is at the $j$-th coordinate. The number $0\leq\epsilon\leq1$ is the mixture parameter. When $\epsilon=0$, this reduces to the given distribution $\mathbf p^i$. We treat $\epsilon$ as a parameter and differentiate $Z(\mathbf p^1,\ldots,\mathbf p^{i-1},(1-\epsilon)\mathbf p^i+\epsilon\mathbf 1_j^i,\mathbf p^{i+1},\ldots,\mathbf p^m)$ with respect to $\epsilon$ for each $i,j$.

More precisely, let
$$\psi_j^i(\mathbf p)=\frac{d}{d\epsilon}Z(\mathbf p^1,\ldots,\mathbf p^{i-1},(1-\epsilon)\mathbf p^i+\epsilon\mathbf 1_j^i,\mathbf p^{i+1},\ldots,\mathbf p^m)\big|_{\epsilon=0}$$
Denote $\bm\psi^i(\mathbf p)=(\psi_j^i(\mathbf p))_{j=1,\ldots,n^i}\in\mathbb R^{n^i}$, and $\bm\psi(\mathbf p)=\text{vec}(\bm\psi^i(\mathbf p):i=1,\ldots,m)\in\mathbb R^N$. We show that $\bm\psi$ possesses the following two properties:

\begin{theorem}
Given $\mathbf p\in\mathcal P$ such that $\mathbf p>\mathbf 0$, we have: 
\begin{enumerate}
\item \begin{equation}
\nabla Z(\mathbf p)'(\mathbf q-\mathbf p)=\sum_{i=1}^m\nabla^iZ(\mathbf p)'(\mathbf q^i-\mathbf p^i)=\sum_{i=1}^m\bm\psi^i(\mathbf p)'(\mathbf q^i-\mathbf p^i)=\bm\psi(\mathbf p)'(\mathbf q-\mathbf p)\label{gradient equivalence}
\end{equation}
for any $\mathbf q^i\in\mathcal P_{n^i}$ and $\mathbf q=\text{vec}(\mathbf q^i:i=1,\ldots,m)$, where $\nabla^iZ(\mathbf p)\in\mathbb R^{n^i}$ is the gradient of $Z$ taken with respect to $\mathbf p^i$.

\item
\begin{equation}
\psi_j^i(\mathbf p)=E_{\mathbf p}[h(\mathbf X)s_j^i(\mathbf X^i)]\label{score function}
\end{equation}
where $s_j^i(\cdot)$ is defined as
\begin{equation}
s_j^i(\mathbf x^i)=\sum_{t=1}^{T^i}\frac{I(x_t^i=y_j^i)}{p_j^i}-T^i\label{score function2}
\end{equation}
for $\mathbf x^i=(x_1^i,\ldots,x_{T^i}^i)\in\mathbb R^{T^i}$. \end{enumerate}
\label{prop:gradient}
\end{theorem}

%

The proof of Theorem \ref{prop:gradient} is in Appendix \ref{sec:proofs}. The first property above states that $\bm\psi(\mathbf p)$ and $\nabla Z(\mathbf p)$ are identical when viewed as directional derivatives, as long as the direction lies within $\mathcal P$. Since the feasible region of optimizations \eqref{sample counterpart} lies in $\mathcal P$, it suffices to focus on $\bm\psi(\mathbf p)$. The second property above states that $\bm\psi(\mathbf p)$ can be estimated unbiasedly in a way similar to the classical likelihood ratio method (\cite{glynn1990likelihood,reiman1989sensitivity}), with $s_j^i(\cdot)$ playing the role of the score function. Since this representation holds without assuming any specific parametric form for $\mathbf p$, we view it as a nonparametric version of the likelihood ratio method.


From \eqref{score function}, an unbiased estimator for $\psi_j^i(\mathbf p)$ using a single simulation run is $(h(\mathbf X)s_j^i(\mathbf X^i))_{i=1,\ldots,m}$, where $\mathbf X=(\mathbf X^1,\ldots,\mathbf X^m)$ is the sample path. The following provides a bound on the variance of this estimator (See Appendix \ref{sec:proofs} for proof):
\begin{lemma}
Assume $h(\mathbf X)$ is bounded a.s., i.e. $|h(\mathbf X)|\leq M$ for some $M>0$, and that $\mathbf p>\mathbf 0$. Each sample for estimating $\psi_j^i(\mathbf p)$, given by $h(\mathbf X)s_j^i(\mathbf X^i)$ using one sample path of $\mathbf X$, possesses a variance bounded from above by $M^2T^i(1-p_j^i)/p_j^i$.
\label{prop:var}
\end{lemma}

The function $\bm\psi(\mathbf p)$ derived via the above Gateaux derivative framework can be interpreted as a discrete version of the so-called influence function in robust statistics (\cite{hampel1974influence,hampel2011robust}), which is commonly used to approximate the first order effect on a given statistics due to contamination of data. In general, the gradient represented by the influence function is defined as an operator on the domain of the random object distributed under $\mathbf p$. Thus, in the continuous case, this object has an infinite-dimensional domain and can be difficult to compute and encode. This is the main reason why we seek for a discretization in the first place.

\section{Frank-Wolfe Stochastic Approximation (FWSA)}\label{sec:procedure}
With the implementable form of the gradient $\bm\psi(\mathbf p)$ described in Section \ref{sec:gradient}, we design a stochastic nonlinear programming technique to solve \eqref{sample counterpart}. We choose to use the Frank-Wolfe method because, for the types of $\hat{\mathcal U}^i$ we consider in Section \ref{sec:discretization}, effective routines exist for solving the induced linearized subproblems.





\subsection{Description of the Algorithm}\label{sec:description}
For convenience denote $\hat{\mathcal U}=\prod_{i=1}^{n^i}\hat{\mathcal U}^i$. We focus on the choices of $\mathcal U^i$ depicted in Section \ref{sec:formulation}, which are all convex and consequently $\hat{\mathcal U}^i$ and also $\hat{\mathcal U}$ are convex.

FWSA works as follows. To avoid repetition we focus only on the minimization formulation in \eqref{generic}. First, pretending that $\nabla Z(\mathbf p)$ can be computed exactly, it iteratively updates a solution sequence $\mathbf p_1,\mathbf p_2,\ldots$ by, given a current solution $\mathbf p_k$, solving
\begin{equation}
\min_{\mathbf p\in\hat{\mathcal U}}\nabla Z(\mathbf p_k)'(\mathbf p-\mathbf p_k)
\label{step optimization}
\end{equation}
Let the optimal solution to \eqref{step optimization} be $\mathbf q_k$. The quantity $\mathbf q_k-\mathbf p_k$ gives a feasible minimization direction starting from $\mathbf p_k$ (recall that $\hat{\mathcal U}$ is convex). This is then used to update $\mathbf p_k$ to $\mathbf p_{k+1}$ via $\mathbf p_{k+1}=\mathbf p_k+\epsilon_k(\mathbf q_k-\mathbf p_k)$ for some step size $\epsilon_k$. This expression can be rewritten as $\mathbf p_{k+1}=(1-\epsilon_k)\mathbf p_k+\epsilon_k\mathbf q_k$, which can be interpreted as a mixture between the distributions $\mathbf p_k$ and $\mathbf q_k$.

When $\nabla Z(\mathbf p_k)$ is not exactly known, one can replace it by an empirical counterpart. Theorem \ref{prop:gradient} suggests that we can replace $\nabla Z(\mathbf p_k)$ by $\bm\psi(\mathbf p_k)$, and so the empirical counterpart of \eqref{step optimization} is
\begin{equation}
\min_{\mathbf p\in\hat{\mathcal U}}\hat{\bm\psi}(\mathbf p_k)'(\mathbf p-\mathbf p_k)\\
\label{step optimization1}
\end{equation}
where $\hat{\bm\psi}(\mathbf p_k)$ is an estimator of $\bm\psi(\mathbf p_k)$ using a sample size $R_k$. Note that all components of $\hat{\bm\psi}(\mathbf p_k)$ can be obtained from these $R_k$ sample paths simultaneously. Letting $\hat{\mathbf q}_k$ be the optimal solution to \eqref{step optimization1}, the update rule will be $\mathbf p_{k+1}=(1-\epsilon_k)\mathbf p_k+\epsilon_k\hat{\mathbf q}_k$ for some step size $\epsilon_k$. The sample size $R_k$ at each step needs to grow suitably to compensate for the bias introduced in solving \eqref{step optimization1}. All these are summarized in Procedure \ref{SCG}.
\begin{algorithm}
\label{proc}
  \caption{FWSA for solving \eqref{generic}}
  \textbf{Initialization: }$\mathbf p_1\in\mathcal P$ where $\mathbf p_1>\mathbf 0$.

  \textbf{Input: }Step size sequence $\epsilon_k$, sample size sequence $R_k$, $k=1,2,\ldots$.

\textbf{Procedure: }For each iteration $k=1,2,\ldots$, given $\mathbf p_k$:

  \begin{algorithmic}
  \State \textbf{1. }Repeat $R_k$ times: Compute
  $$h(\mathbf X)s_j^i(\mathbf X^i)\text{\ \ \ \ for all $i=1,\ldots,m$}$$
  using one sample path $\mathbf X=(\mathbf X^1,\ldots,\mathbf X^m)$, where $s_j^i(\mathbf X^i)=\sum_{t=1}^{T^i}I(X_t^i=y_j^i)/p_j^i-T^i$ for $j=1,\ldots,n^i$ and $i=1,\ldots,m$. Call these $R_k$ i.i.d. replications $\zeta_j^i(r)$, for $j=1,\ldots,n^i,\ i=1,\ldots,m,\ r=1,\ldots,R_k$.

  \State \textbf{2. } Estimate $\bm\psi(\mathbf p_k)$ by
  $$\hat{\bm\psi}(\mathbf p_k)=(\hat\psi_j^i(\mathbf p_k))_{i=1,\ldots,m,\ j=1,\ldots,n^i}=\left(\frac{1}{R_k}\sum_{r=1}^{R_k}\zeta_j^i(r)\right)_{i=1,\ldots,m,\ j=1,\ldots,n^i}.$$


\State \textbf{3. } Solve $\hat{\mathbf q}_k\in\text{argmin}_{\mathbf p\in\hat{\mathcal U}}\hat{\bm\psi}(\mathbf p_k)'(\mathbf p-\mathbf p_k)$.

\State \textbf{4. } Update $\mathbf p_{k+1}=(1-\epsilon_k)\mathbf p_k+\epsilon_k\hat{\mathbf q}_k$.

  \end{algorithmic}\label{SCG}
\end{algorithm}

\subsection{Solving the Subproblem}\label{sec:subprogram}
By \eqref{gradient equivalence} and the separability of uncertainty set $\hat{\mathcal U}=\prod_{i=1}^m\hat{\mathcal U}^i$, the subproblem at each iteration can be written as
\begin{equation}
\min_{\mathbf q\in\hat{\mathcal U}}\sum_{i=1}^m\hat{\bm\psi}^i(\mathbf p)'(\mathbf q^i-\mathbf p^i)=\sum_{i=1}^m\min_{\mathbf q^i\in\hat{\mathcal U}^i}\hat{\bm\psi}^i(\mathbf p)'(\mathbf q^i-\mathbf p^i)\label{step multiple}
\end{equation}
where $\hat{\bm\psi}^i(\mathbf p)=(\hat\psi_j^i(\mathbf p))_{j=1,\ldots,n^i}$ is the empirical counterpart of $\bm\psi^i(\mathbf p)$ obtained in Algorithm \ref{SCG}. Hence \eqref{step multiple} can be solved by $m$ separate convex programs. The update step follows by taking $\mathbf p_{k+1}=\text{vec}(\mathbf p_{k+1}^i:i=1,\ldots,m)$, where $\mathbf p_{k+1}^i=(1-\epsilon_k)\mathbf p_k^i+\epsilon_k\hat{\mathbf q}_k^i$ and $\hat{\mathbf q}_k^i$ is the solution to the $i$-th separate program. 

The separate programs in \eqref{step multiple} can be efficiently solved for the uncertainty sets considered in Section \ref{sec:discretization}. To facilitate discussion, we denote a generic form of each separate program in \eqref{step multiple} as
\begin{equation}
\min_{\mathbf p^i\in\hat{\mathcal U}^i}\bm\xi'\mathbf p^i
\label{step optimization2}
\end{equation}
for an arbitrary vector $\bm\xi=(\xi_j)_{j=1,\ldots,n^i}\in\mathbb R^{n^i}$.
\\

\noindent\emph{Case 1 in Theorem \ref{sample thm}: Moment and support constraints. }Consider $\hat{\mathcal U}^i=\{\mathbf p^i\in\mathcal P^{n^i}:{\mathbf f_l^i}'\mathbf p^i\leq\mu_l^i,l=1,\ldots,s^i\}$ where $\mathbf f_l^i=(f_l(y_j^i))_{j=1,\ldots,n^i}\in\mathbb R^{n^i}$. Then \eqref{step optimization2} is a linear program.
\\


\noindent\emph{Case 2 in Theorem \ref{sample thm}: $\phi$-divergence neighborhood. }Consider
\begin{equation}
\hat{\mathcal U}^i=\{\mathbf p^i\in\mathcal P^{n^i}:d_\phi(\mathbf p^i,\mathbf p_b^i)\leq\eta^i\}\label{phi uncertainty reformulated}
\end{equation}
where $\mathbf p_b^i=(p_{b,j}^i)_{j=1,\ldots,n^i}\in\mathcal P^{n^i}$ and $d_\phi(\mathbf p^i,\mathbf p_b^i)=\sum_{j=1}^{n^i}p_{b,j}^i\phi(p_j^i/p_{b,j}^i)$. We have:

\begin{proposition}
Consider \eqref{step optimization2} with $\hat{\mathcal U}^i$ presented in \eqref{phi uncertainty reformulated}, where $\mathbf p_b^i>\mathbf 0$. Let $\phi^*(t)=\sup_{x\geq0}\{tx-\phi(x)\}$ be the conjugate function of $\phi$, and define $0\phi^*(s/0)=0$ if $s\leq0$ and $0\phi^*(s/0)=+\infty$ if $s>0$. Solve the program
\begin{equation}
(\alpha^*,\lambda^*)\in\text{argmax}_{\alpha\geq0,\lambda\in\mathbb R}\left\{-\alpha\sum_{j=1}^{n^i}p_{b,j}^i\phi^*\left(-\frac{\xi_j+\lambda}{\alpha}\right)-\alpha\eta^i-\lambda\right\}\label{opt phi}
\end{equation}
An optimal solution $\mathbf q^i=(q_j^i)_{j=1,\ldots,n^i}$ for \eqref{step optimization2} is
\begin{enumerate}
\item If $\alpha^*>0$, then
\begin{equation}
q_j^i=p_{b,j}^i\cdot \text{argmax}_{r\geq0}\left\{-\frac{\xi_j+\lambda^*}{\alpha^*}r-\phi(r)\right\}\label{opt phi1}
\end{equation}

\item If $\alpha^*=0$, then
\begin{equation}
q_j^i=\left\{\begin{array}{ll}
\frac{p_{b,j}^i}{\sum_{j\in\mathcal M^i}p_{b,j}^i}&\text{\ for\ }j\in\mathcal M^i\\
0&\text{\ otherwise}
\end{array}\right.\label{opt phi2}
\end{equation}
\end{enumerate}
where $\mathcal M^i=\text{argmin}_j\xi_j$, the set of indices $j\in\{1,\ldots,n^i\}$ that have the minimum $\xi_j$.
\label{phisolution}
\end{proposition}

Operation \eqref{opt phi} involves a two-dimensional convex optimization. Note that both the function $\phi^*$ and the solution to the $n^i$ one-dimensional maximization \eqref{opt phi1} have closed-form expressions for all common $\phi$-divergence (\cite{pardo2005statistical}). The proof of Proposition \ref{phisolution} follows closely from \cite{ben2013robust} and is left to Appendix \ref{sec:proofs}.


In the special case where $\phi=x\log x-x+1$, i.e. KL divergence, the solution scheme can be simplified to a one-dimensional root-finding problem. More precisely, we have
\begin{proposition}
Consider \eqref{step optimization2} with $\hat{\mathcal U}^i$ presented in \eqref{phi uncertainty reformulated}, where $\phi(x)=x\log x-x+1$ and $\mathbf p_b^i>\mathbf 0$. Denote $\mathcal M^i=\text{argmin}_j\xi_j$ as in Proposition \ref{phisolution}. An optimal solution $\mathbf q^i=(q_j^i)_{j=1,\ldots,n^i}$ for \eqref{step optimization2} is:
\begin{enumerate}
\item If $-\log\sum_{j\in\mathcal M^i}p_{b,j}^i\leq\eta^i$, then
\begin{equation}
q_j^i=\left\{\begin{array}{ll}
\frac{p_{b,j}^i}{\sum_{j\in\mathcal M^i}p_{b,j}^i}&\text{\ for\ }j\in\mathcal M^i\\
0&\text{\ otherwise}
\end{array}\right.\label{opt2}
\end{equation}

\item If $-\log\sum_{j\in\mathcal M^i}p_{b,j}^i>\eta^i$, then
\begin{equation}
q_j^i=\frac{p_{b,j}^ie^{\beta\xi_j}}{\sum_{j=1}^{n^i}p_{b,j}^ie^{\beta\xi_j}}\label{opt1}
\end{equation}
for all $j$, where $\beta<0$ satisfies
\begin{equation}
\beta{\varphi_{\bm\xi}^i}'(\beta)-\varphi_{\bm\xi}^i(\beta)=\eta^i\label{root}
\end{equation}
Here $\varphi_{\bm\xi}^i(\beta)=\log\sum_{j=1}^{n^i}p_{b,j}^ie^{\beta\xi_j}$ is the logarithmic moment generating function of $\bm\xi$ under $\mathbf p_b^i$.
\end{enumerate}\label{KLsolution}
\end{proposition}
The proof of Proposition \ref{KLsolution} follows from techniques in, e.g., \cite{hansen2008robustness}, and is left to Appendix \ref{sec:proofs}.

\section{Theoretical Guarantees of FWSA} \label{sec:guarantees}

This section shows the convergence properties of our proposed FWSA. We first present results on almost sure convergence, followed by a local convergence rate analysis. Throughout our analysis we assume that the subproblem at any iteration can be solved using deterministic optimization routine to a negligible error.

\subsection{Almost Sure Convergence}
An important object that we will use in our analysis is the so-called Frank-Wolfe (FW) gap (\cite{frank1956algorithm}): For any $\tilde{\mathbf p}\in\hat{\mathcal U}$, let $g(\tilde{\mathbf p})=-\min_{\mathbf p\in\hat{\mathcal U}}\bm\psi(\tilde{\mathbf p})'(\mathbf p-\tilde{\mathbf p})$, which is the negation of the optimal value of the next subproblem when the current solution is $\tilde{\mathbf p}$. Note that $g(\tilde{\mathbf p})$ is non-negative for any $\tilde{\mathbf p}\in\hat{\mathcal U}$, since one can always take $\mathbf p=\tilde{\mathbf p}$ in the definition of $g(\tilde{\mathbf p})$ to get a lower bound 0. In the case of convex objective function, it is well-known that $g(\tilde{\mathbf p})$ provides an upper bound of the actual optimality gap (\cite{frank1956algorithm}). However, we shall make no convexity assumption in our subsequent analysis, and will see that $g(\tilde{\mathbf p})$ still plays an important role in bounding the local convergence rate of our procedure under the conditions we impose.

Our choices on the step size $\epsilon_k$ and sample size per iteration $R_k$ of the procedure are as follows:
\begin{assumption}
We choose $\epsilon_k,k=1,2,\ldots$ that satisfy
$$\sum_{k=1}^\infty\epsilon_k=\infty\text{\ \ \ \ and\ \ \ \ }\sum_{k=1}^\infty\epsilon_k^2<\infty$$\label{tuning}
\end{assumption}
\begin{assumption}
The sample sizes $R_k,k=1,2,\ldots$ are chosen such that
$$\sum_{k=1}^\infty\frac{\epsilon_k}{\sqrt{R_k}}\prod_{j=1}^{k-1}(1-\epsilon_j)^{-1/2}<\infty$$
where for convenience we denote $\prod_{j=1}^0(1-\epsilon_j)^{-1/2}=1$.\label{sample size tuning}
\end{assumption}

Note that among all $\epsilon_k$ in the form $c/k^\alpha$ for $c>0$ and $\alpha>0$, only $\alpha=1$ satisfies both Assumptions \ref{tuning} and \ref{sample size tuning} and avoids a super-polynomial growth in $R_k$ simultaneously (recall that $R_k$ represents the simulation effort expended in iteration $k$, which can be expensive). To see this, observe that Assumption \ref{tuning} asserts $\alpha\in(1/2,1]$. Now, if $\alpha<1$, then it is easy to see that $\prod_{j=1}^{k-1}(1-\epsilon_j)^{-1/2}$ grows faster than any polynomials, so that $R_k$ cannot be polynomial if Assumption \ref{sample size tuning} needs to hold. On the other hand, when $\alpha=1$, then $\prod_{j=1}^{k-1}(1-\epsilon_j)^{-1/2}$ grows at rate $\sqrt k$ and it is legitimate to choose $R_k$ growing at rate $k^\beta$ with $\beta>1$.

Assumption \ref{tuning} is standard in the SA literature. The growing per-iteration sample size in Assumption \ref{sample size tuning} is needed to compensate for the bias caused by the subproblem in FWSA. Note that in standard SA, a solution update is obtained by moving in the gradient descent direction, and Assumption \ref{tuning} suffices if this direction is estimated unbiasedly. In FWSA, the subprogram introduces bias on the feasible direction despite the unbiasedness of the gradient estimate. The increasing simulation effort at each iteration is introduced to shrink this bias as the iteration proceeds. We also note that the expression $\prod_{j=1}^{k-1}(1-\epsilon_j)^{-1/2}$ in Assumption \ref{sample size tuning} is imposed to compensate for a potentially increasing estimation variance, due to the form of the gradient estimator depicted in \eqref{score function} and \eqref{score function2} that possesses $p_j^i$ in the denominator and thus the possibility of having a larger variance as the iteration progresses.


We state our result on almost sure convergence in two parts. The first part only assumes the continuity of $g(\cdot)$. The second part assumes a stronger uniqueness condition on the optimal solution, stated as:
\begin{assumption}
There exists a unique minimizer $\mathbf p^*$ for $\min_{\mathbf p\in\hat{\mathcal U}}Z(\mathbf p)$. Moreover, $g(\cdot)$ is continuous over $\hat{\mathcal U}$ and $\mathbf p^*$ is the only feasible solution such that $g(\mathbf p^*)=0$.\label{main assumption}
\end{assumption}

In light of Assumption \ref{main assumption}, $g$ plays a similar role as the gradient in unconstrained problems. The condition $g(\mathbf p^*)=0$ in Assumption \ref{main assumption} is a simple implication of the optimality of $\mathbf p^*$ (since $g(\mathbf p^*)>0$ would imply the existence of a better solution).

Our convergence result is:
\begin{theorem}
Suppose that $h(\mathbf X)$ is bounded a.s. and that Assumptions \ref{tuning}-\ref{sample size tuning} hold. We have the following properties on $\mathbf p_k$ generated in Algorithm \ref{SCG} :
\begin{enumerate}
\item Assume that $g(\cdot)$ is continuous and an optimal solution exists. Then $D(Z(\mathbf p_k),\mathcal Z^*)\to0$ a.s., where $\mathcal Z^*=\{Z(\mathbf p):\mathbf p\text{\ satisfies\ }g(\mathbf p)=0\}$ and $D(x,A)=\inf_{y\in A}\|x-y\|$ for any point $x$ and set $A$ in the Euclidean space.\label{as part1}
\item Under Assumption \ref{main assumption}, $\mathbf p_k$ converge to $\mathbf p^*$ a.s..\label{as part2}
\end{enumerate}
\label{as}
\end{theorem}

Part \ref{as part1} of Theorem \ref{as} states that the objective value generated by Algorithm \ref{SCG} will eventually get close to an objective value evaluated at a point where the FW gap is zero. Part \ref{as part2} strengthens the convergence to the unique optimal solution $\mathbf p^*$ under Assumption \ref{main assumption}. In practice, this uniqueness condition may not hold, and we propose combining Algorithm \ref{SCG} with multi-start of the initial solution $\mathbf p_1$ as a remedy. Section \ref{expt:mltcls} and Appendix \ref{sec:multi-start} show some numerical results on this strategy.

\subsection{Local Convergence Rate}\label{sec:local rate}
We impose several additional assumptions. The first is a Lipchitz continuity condition on an optimal solution for the generic subproblem \eqref{step optimization2}, with respect to the coefficients in the objective in a neighborhood of the gradient evaluated at $\mathbf p^*$. Denote $\mathbf v(\bm\xi)$ as an optimal solution of \eqref{step optimization2}.
\begin{assumption}
We have
$$\|\mathbf v(\bm\xi_1)-\mathbf v(\bm\xi_2)\|\leq L\|\bm\xi_1-\bm\xi_2\|$$
for some $L>0$, for any $\bm\xi_1,\bm\xi_2\in\mathcal N_\Delta(\bm\psi(\mathbf p^*))$, where $\mathcal N_\Delta(\bm\psi(\mathbf p^*))$ denotes a Euclidean neighborhood of $\bm\psi(\mathbf p^*)$ with radius $\Delta$, and $\mathbf p^*$ is assumed to be the unique optimal solution for $\min_{\mathbf p\in\hat{\mathcal U}}Z(\mathbf p)$.\label{bias}
\end{assumption}

%

Next, we denote $\mathbf q(\tilde{\mathbf p})$ as an optimizer in the definition of the FW gap at $\tilde{\mathbf p}$, i.e. $\mathbf q(\tilde{\mathbf p})\in\text{argmin}_{\mathbf p}\bm\psi(\tilde{\mathbf p})'(\mathbf p-\tilde{\mathbf p})$.
\begin{assumption}
$$g(\mathbf p)\geq c\|\bm\psi(\mathbf p)\|\|\mathbf q(\mathbf p)-\mathbf p\|$$
for any $\mathbf p\in\hat{\mathcal U}$, where $c>0$ is a small constant.\label{angle}
\end{assumption}

\begin{assumption}
$$\|\bm\psi(\mathbf p)\|>\tau>0$$
for any $\mathbf p\in\hat{\mathcal U}$, for some constant $\tau$.\label{nonzero gradient}
\end{assumption}


Assumption \ref{angle} guarantees that the angle between the descent direction and the gradient must be bounded away from $90^\circ$ uniformly at any point $\mathbf p$. This assumption has been used in the design and analysis of gradient descent methods for nonlinear programs that are singular (i.e. without assuming the existence of the Hessian matrix; \cite{bertsekas1999nonlinear}, Proposition 1.3.3).

The non-zero gradient condition in Assumption \ref{nonzero gradient} effectively suggests that a local optimum must occur at the relative boundary of $\hat{\mathcal U}$ (i.e. the boundary with respect to the lower-dimensional subspace induced by the probability simplex constraint), which warrants further explanation. Note that the other alternate scenario for local optimality will be that it occurs in the interior region of the feasible set $\hat{\mathcal U}$. In the latter scenario, the gradient at the optimal solution is zero. While the convergence analysis can be simplified (and plausibly give a better rate) under this scenario, the statistical implication brought by this scenario is rather pathological. Note that our optimizations are imposed on decision variables that are input probability distributions. As discussed at the end of Section \ref{sec:gradient}, the gradient vector $\bm\psi(\mathbf p)$ is the influence function for the performance measure $Z(\cdot)$. If the influence function is zero, it is known that a Gaussian limit does not hold in the central limit theorem as the input sample size gets large (where the central limit theorem is on the difference between a simulation driven by empirical distributions and the truth). Instead, a $\chi^2$-limit occurs (\cite{serfling2009approximation}, Section 6.4.1, Theorem B). Such type of limit is unusual and has never been reported in simulation analysis. Indeed, in all our experiments, the obtained local optimal solution is always at the boundary. For this reason we impose Assumption \ref{nonzero gradient} rather than a more straightforward zero-gradient type condition. 



The following are our main results on convergence rate, first on the FW gap $g(\mathbf p_k)$, and then the optimality gap $Z(\mathbf p_k)-Z(\mathbf p^*)$, in terms of the number of iterations $k$. Similar to almost sure convergence, we assume here that the deterministic routine for solving the subproblems can be carried out with high precision.

\begin{theorem}
Suppose $|h(\mathbf X)|\leq M$ for some $M>0$ and that Assumptions \ref{tuning}-\ref{nonzero gradient} hold. Additionally, set
$$\epsilon_k=\frac{a}{k}\text{\ \ \ \ and\ \ \ \ }R_k=bk^\beta$$
when $k>a$, and arbitrary $\epsilon_k<1$ when $k\leq a$. Given any $0<\varepsilon<1$, it holds that, with probability $1-\varepsilon$, there exists a large enough positive integer $k_0$ and small enough positive constants $\nu,\vartheta,\varrho$ such that $0<g(\mathbf p_{k_0})\leq\nu$, and for $k\geq k_0$,
\begin{equation}
g(\mathbf p_k)\leq\frac{A}{k^C}+B\times\left\{\begin{array}{ll}\frac{1}{(C-\gamma)k^\gamma}&\text{if $0<\gamma<C$}\\
\frac{1}{(\gamma-C)(k_0-1)^{\gamma-C}k^C}&\text{if $\gamma>C$}\\
\frac{\log((k-1)/(k_0-1))}{k^C}&\text{if $\gamma=C$}
\end{array}\right.\label{main}
\end{equation}
where
$$A=g(\mathbf p_{k_0})k_0^C,$$
$$B=\left(1+\frac{1}{k_0}\right)^C\left(a\varrho+\frac{2a^2\varrho K}{c\tau k_0}\left(\frac{\nu}{c\tau}+L\vartheta\right)\right)$$
and
\begin{equation}
C=a\left(1-\frac{2KL\vartheta}{c\tau}-\frac{2K\nu}{c^2\tau^2}\right)\label{C}
\end{equation}
Here the constants $L,c,\tau$ appear in Assumptions \ref{bias}, \ref{angle} and \ref{nonzero gradient} respectively. The sample size power $\beta$ needs to be chosen such that $\beta>2\gamma+a+1$. More precisely, the constants $a,b,\beta$ that appear in the specification of the algorithm, the other constants $k_0,\vartheta,\varrho,\gamma,K$, and two new constants $\rho>1$ and $\delta>0$ are chosen to satisfy Conditions \ref{c1}-\ref{c8} listed in Appendix \ref{sec:proofs}.
\label{rate thm}
\end{theorem}

\begin{corollary}
Suppose that all the assumptions are satisfied and all the constants are chosen as indicated in Theorem \ref{rate thm}. Then with probability $1-\varepsilon$, there exists a large enough positive integer $k_0$ and small enough positive constants $\nu,\vartheta,\varrho$ such that $0\leq g(\mathbf p_{k_0})\leq\nu$, and for $k\geq k_0$,
\begin{equation}
Z(\mathbf p_k)-Z(\mathbf p^*)\leq\frac{D}{k-1}+\frac{E}{(k-1)^C}+F\times\left\{\begin{array}{ll}\frac{1}{(C-\gamma)\gamma(k-1)^\gamma}&\text{\ if\ }0<\gamma<C\\
\frac{1}{(\gamma-C)(k_0-1)^{\gamma-C}C(k-1)^C}&\text{\ if\ }\gamma>C\\
\frac{\log((k-1)/(k_0-1))}{C(k-1)^C}&\text{\ if\ }\gamma=C
\end{array}\right.\label{main1}
\end{equation}
where
$$D=a^2K,\ \ E=\frac{aA}{C},\ \ F=aB$$
and $a,A,B,C,K$ are the same constants as in Theorem \ref{rate thm}.\label{rate cor}
\end{corollary}


A quick summary extracted from Theorem \ref{rate thm} and Corollary \ref{rate cor} is the following: Consider the local convergence rate denominated by workload, i.e. the number of simulation replications. To achieve the most efficient rate, approximately speaking, $a$ should be chosen to be $1+\omega$ and $\beta$ chosen to be $5+\zeta+\omega$ for some small $\omega,\zeta>0$. The local convergence rate is then $O(W^{-1/(6+\zeta+\omega)})$ where $W$ is the total number of simulation replications. 

Note that the bounds in Theorem \ref{rate thm} and Corollary \ref{rate cor} are local asymptotic statements since they only hold starting from $k\geq k_0$ and $g(\mathbf p_k)\leq\nu$ for some large $k_0$ and small $\nu$. It should be cautioned that they do not say anything about the behavior of the algorithm before reaching the small neighborhood of $\mathbf p^*$ as characterized by $0\leq g(\mathbf p_{k_0})\leq\nu$.
The above summary therefore should be interpreted in the way that, given the algorithm has already run $k_0$ number of replications and $g(\mathbf p_k)\leq\nu$ for a suitably small $\nu$ (which occurs with probability 1 by Theorem \ref{as}), the convergence rate of $O(W^{-1/(6+\zeta+\omega)})$ for the optimality gap is guaranteed with probability $1-\varepsilon$ starting from that point. 



The summary above is derived based on the following observations:
\begin{enumerate}
\item The local convergence rate of the optimality gap, in terms of the number of iterations $k$, is at best $O(1/k^{C\wedge\gamma\wedge1})$. This is seen by \eqref{main1}.
\item We now consider the convergence rate in terms of simulation replications. Note that at iteration $k$, the cumulative number of replications is of order $\sum_{j=1}^kj^\beta\approx k^{\beta+1}$. Thus from Point 1 above, the convergence rate of the optimality gap in terms of replications is of order $1/W^{(C\wedge\gamma\wedge1)/(\beta+1)}$.

\item The constants $C$ and $\gamma$ respectively depend on $a$, the constant factor in the step size, and $\beta$, the geometric growth rate of the sample size, as follows:
\begin{enumerate}
\item \eqref{C} defines $C=a(1-2KL\vartheta/(c\tau)-2K\nu/(c^2\tau^2))$. For convenience, we let $\omega=2KL\vartheta/(c\tau)+2K\nu/(c^2\tau^2)$, and so $C=a(1-\omega)$.
\item From Condition \ref{c5} in Theorem \ref{rate thm} (shown in Appendix \ref{sec:proofs}), we have $\beta=2\gamma+\rho a+2+\zeta$ for some $\zeta>0$. In other words $\gamma=(\beta-\rho a-\zeta-2)/2$.
\end{enumerate}
\item Therefore, the convergence rate in terms of replications is $1/W^{((a(1-\omega))\wedge((\beta-\rho a-\zeta-2)/2)\wedge1)/(\beta+1)}$. Let us focus on maximizing
\begin{equation}
\frac{(a(1-\omega))\wedge((\beta-\rho a-\zeta-2)/2)\wedge1}{\beta+1}\label{index}
\end{equation}
over $a$ and $\beta$, whose solution is given by the following lemma:
\begin{lemma}
The maximizer of \eqref{index} is given by
$$a=\frac{1}{1-\omega},\ \ \beta=\frac{\rho}{1-\omega}+\zeta+4$$
and the optimal value is
$$\frac{1}{\rho/(1-\omega)+\zeta+5}$$\label{index max}
\end{lemma}
The proof is in Appendix \ref{sec:proofs}. With Lemma \ref{index max}, let us choose $\vartheta$ and $\nu$, and hence $\omega$, to be small. We also choose $\rho$ to be close to 1. (Unfortunately, these choices can lead to a small size of neighborhood around $\mathbf p^*$ in which the convergence rate holds.) This gives rise to the approximate choice that $a\approx1+\omega$ and $\beta\approx5+\zeta+\omega$. The convergence rate is then $O(W^{-1/(6+\zeta+\omega)})$.
\end{enumerate}

We compare our results to some recent work in stochastic FW. \cite{hazan2016variance} showed that to achieve $\epsilon$ error in terms of the optimality gap one needs $O(1/\epsilon^{1.5})$ number of calls to the gradient estimation oracle, when the objective function is strongly convex. \cite{reddi2016stochastic} showed that the number needed increases to $O(1/\epsilon^4)$ for non-convex objectives, and suggested several more sophisticated algorithms to improve the rate. Corollary \ref{rate cor} and our discussion above suggests that we need $O(1/\epsilon^{6+\zeta+\omega})$ sample size, for some small $\zeta,\omega>0$, a rate that is inferior to the one achieved in \cite{reddi2016stochastic}. However, \cite{reddi2016stochastic} has assumed that the gradient estimator is uniformly bounded over the feasible space, a condition known as $G$-Lipschitz (Theorem 2 in \cite{reddi2016stochastic}), which does not hold in our case due to the presence of $p_j^i$ in the denominator in \eqref{score function2} that gives a potentially increasing estimation variance as the iteration progresses. This complication motivates our sample size and step size sequences depicted in Assumption \ref{sample size tuning} and the subsequent analysis. On the other hand, if Assumption \ref{bias} is relaxed to hold for any $\bm\xi_1,\bm\xi_2\in\mathbb R^N$, it can be seen that by choosing $\beta\approx3+\zeta+\omega$ our complexity improves to $O(1/\epsilon^{4+\zeta+\omega})$, which almost matches the one in \cite{reddi2016stochastic} (see Remark \ref{remark:rate} in Appendix \ref{sec:proofs}). However, such a relaxed condition would not hold if the constraints are linear, because the optimal solutions of the subproblems are located at the corner points and will jump from one to the other under perturbation of the objective function.

\section{Numerical Experiments} \label{sec:numerics}

\newcommand{\tndi}{\rightarrow \infty}
\newcommand{\msx}{{\mathbf x}}
\newcommand{\mbX}{{\mathbf X}}
\newcommand{\mse}{{\mathsf{e}}}
\newcommand{\tndo}{\rightarrow 0}
\newcommand{\E}{\mathbb{E}}
\newcommand{\ind}{\mathbb{I}}
\newcommand{\DefAs}{\mbox{$\,\stackrel{\bigtriangleup}{=}\,$}}
\newcommand{\ProofEnd}{\mbox{$\Box$}}
\newcommand{\ubar}[1]{\underline{#1}}
\newcommand{\obar}[1]{{\overline{#1}}}

This section describes two sets of numerical experiments. The first set (Section \ref{expt:mltcls}) studies the performance guarantees from Section \ref{sec:discretization} regarding our randomized discretization strategy and the tightness of the bounds coming from moment constraints. The second set of experiments (Section \ref{sec:gg1example}) studies the numerical convergence of FWSA. The appendix provides additional details and results. Unless specified, in all experiments we terminate the FWSA algorithm at iteration $k$ if at least one of the
following criteria is met (as an indication that the convergence studied in Section \ref{sec:guarantees} is attained):
\begin{itemize}
\item The cumulative simulation replications $W_k$ reaches $5\times
  10^8$, or
\item The relative difference between objective value $Z(\mbp_k)$ and
  the average of the observed values in $30$ previous iterations,
  $(\sum_{v=1}^{30} Z(\mbp_{k-v}) )/30$, is below $5\times 10^{-5}$,
  or
\item The gradient estimate $\hat{\bm\psi}(\mbp_k)$ has an $l_2$-norm
  smaller than $1\times 10^{-3}$.
\end{itemize}


\subsection{Performance Bounds for Multiple Continuous and Unbounded Input Models}\label{expt:mltcls}
We use the example of a generic multi-class $M/G/1$ queue where jobs from
three distinct classes arrive and are attended to by one server. Such structures are common in service systems such as call-centers. Let
$\mathbf P = \{P^1, P^2, P^3\}$ represent all the constituent probability measures, where each $P^i = \{P^{i,j}\},\,\,i=1,2,3$ with
$j=1$ for interarrival and $j=2$ for service, denotes the joint
measure of the interarrival and service distributions of jobs of class $i$. The performance measure
of interest is the weighed average waiting time:
\begin{equation}
Z(\mathbf P) = E_{\mathbf P} \left[ \sum_{i=1}^3 \left(c_i\,\,\frac 1 {T^i} \sum_{t=1}^{T^i}  W_t^{i}\right) \right],
\end{equation}
where the average is observed up to a (fixed) $T^i=500$ customers of
class $i$ and $c_i$ is the cost assigned to its waiting times. Jobs within each class are served on a
first-come-first-served basis. The server uses a fixed
priority ordering based on the popular $c\mu$ rule (\cite{kleinrockv2}), which prioritizes the class on the next serving in decreasing order of the product of $c_i$ and the mean service rate $\mu^i$ of class
$i$ (as discussed momentarily, the $\mu^i$'s are unknown, so we fix a specific guess throughout this example).


To handle the uncertainty in specifying the interarrival and service time
distributions of each class (due to, e.g., the novelty of the service
operation with little pre-existing data), we use the uncertainty set based on moment constraints on the $P^i$ as:
\begin{equation}
\mathcal U=\prod_i \mathcal U^i,\,\,\,\text{where}\,\,\mathcal U^i = \{P^i:\ubar{\mu}^{i,j}_l \leq E_{P^i}[(X^{i,j})^l]\leq \obar{\mu}^{i,j}_l,\,\,l=1,2,\,\,j=1,2\}\label{indep moment set}
\end{equation}
where the index $l=1,2$ represents the first two moments of marginals
$P^{i,j}$. This set is motivated from queueing theory that mean system responses could depend on the mean and variance of the input distributions.
The moment bounds $\ubar{\mu}^{i,j}_l$ and $\obar{\mu}^{i,j}_l$ can be
specified from prior or expert opinion. Here, to test the information
value with respect to the accuracy of the moments, we specify the
bounds from a confidence interval on the corresponding moments
calculated from $N_s$ synthetically generated observations for each
$i,j$. For example,
\[\obar{\mu}^{i,j}_l = \hat{\mu}^{i,j}_l + t_{\alpha/2,N_s-1}
\hat{\sigma}^{i,j}_l / \sqrt{N_s},
\]
where $t_{\alpha/2,N_s-1}$ is the $(1-\alpha/2)$-quantile of the Student-t
distribution with degree of freedom $N_s-1$, $\hat{\mu}^{i,j}_l$ is the empirical $l-$th moment
and $\hat{\sigma}^{i,j}_l$ is the associated sample standard deviation as observed from the
$N_s$ data points. Suppose that the true marginal distribution
of interarrival times for each class is exponential with rate $0.5$
and the true service distribution of the three classes are
exponentials with rates $2.25, 2.0 $ and $1.75$ respectively, to yield
an overall traffic intensity of $0.75$.

The FWSA algorithm is run by first sampling a discrete approximate
support from bivariate independent-marginal lognormal distributions as
representative of each $P^i$ with support size $n=50, 100, 250$ (we
assume the support size corresponding to each distribution $P^i$ is
all equal to $n$). Theorem \ref{sample thm} suggests that selecting lognormal distributions is reasonable if the modeler conjectures that the true distributions are
light-tailed. Here we set the means and standard deviations of the
lognormals to $1$.  The parameter $n$ should ideally be large to
minimize discretization error, but this
pays a penalty in the slowness of the FWSA algorithm.

\begin{figure}[!htbp]
  \centering
  \begin{subfigure}[b]{0.45\textwidth}
    \includegraphics[width=\textwidth]{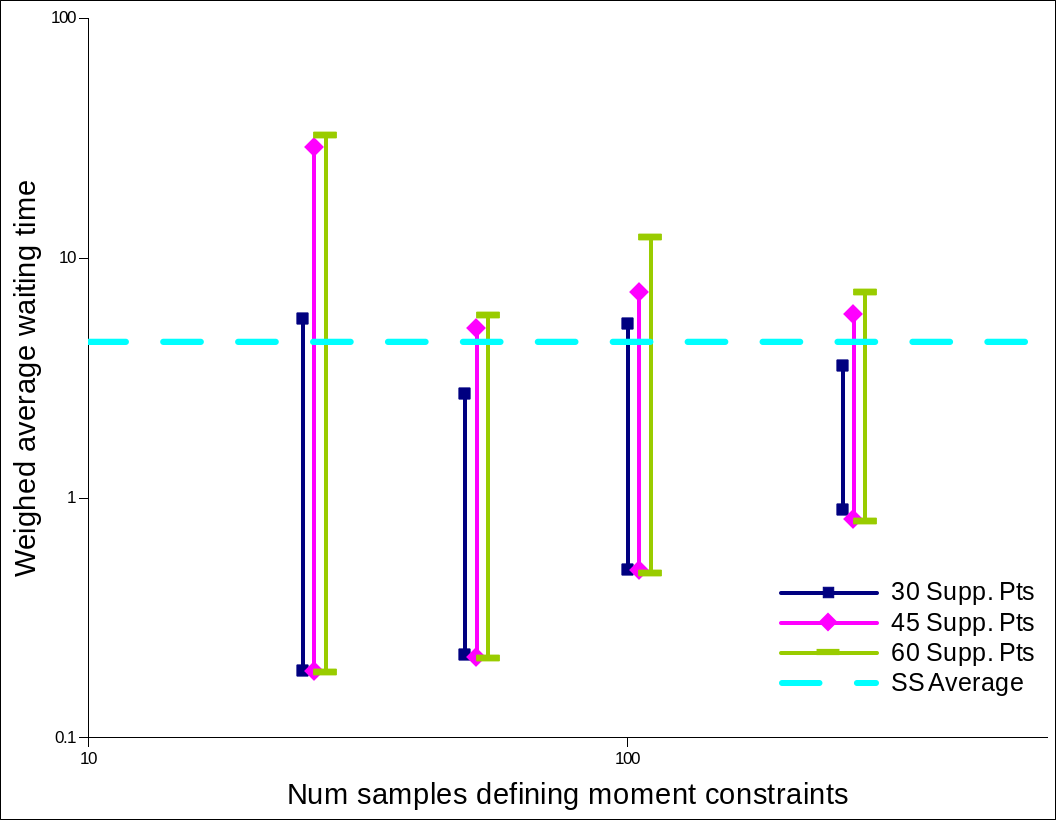}
    \caption{lognormal for discretization\label{fig:mm13lgnml}}
\end{subfigure}
\quad
  \begin{subfigure}[b]{0.45\textwidth}
    \includegraphics[width=\textwidth]{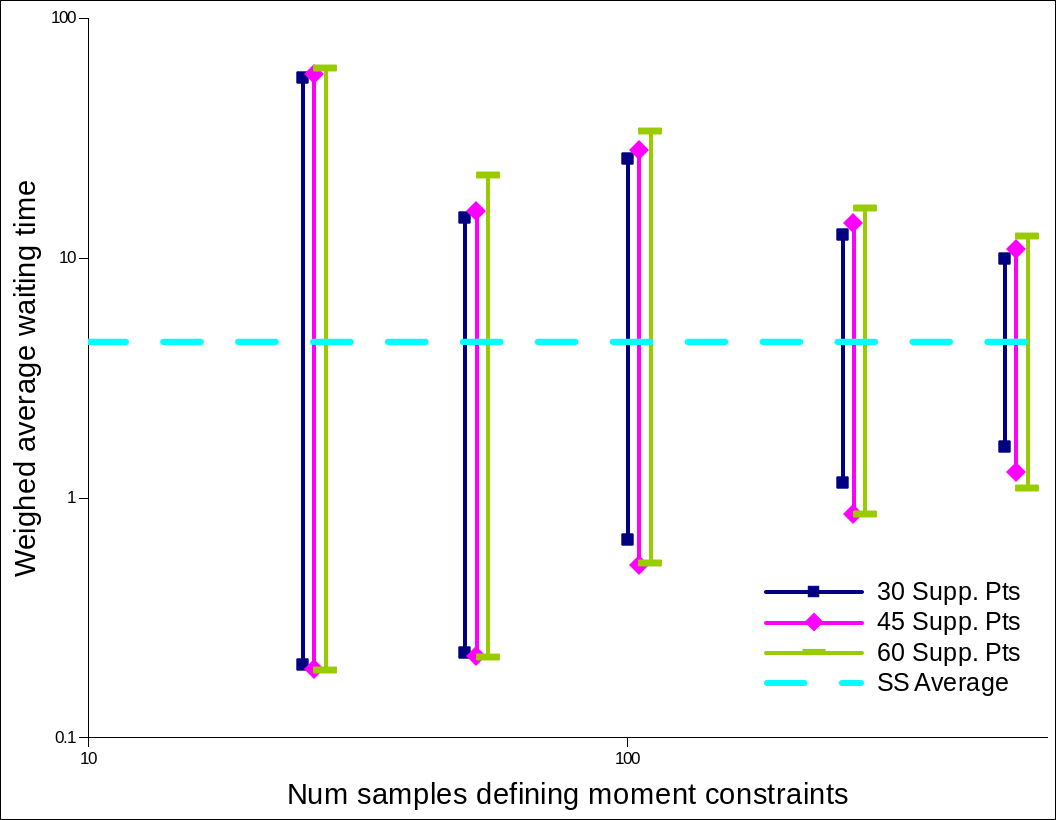}
    \caption{exponential for discretization\label{fig:mm13exp}}
  \end{subfigure}
  \caption{The range from max to min worst-case objectives when $N_s$
    and $n$ vary as indicated.\label{fig:mm13perf} The dotted-line
    indicates the expected steady-state performance under the true distribution.
    }
\end{figure}

Figure~\ref{fig:mm13lgnml} shows the output of our approach over
various $n$ and $N_s$ to illustrate the effect of discretization and the accuracy of moment
information on the tightness of our bounds. The true steady-state performance measure of the
multiclass $M/M/1$ system, available in explicit form (\cite{kleinrockv2}), is indicated as the dotted-line in each plot. The bounds provided by
our method are all seen to cover the
true performance value when $n\geq45$. This is predicted by Theorem \ref{sample thm}
as the moment constraints are all correctly calibrated (i.e. contain
the true moments) in this example. Moreover, as predicted by discussion point 5 in Section \ref{sec:randomized}, the obtained intervals widen as $n$ increases, since the expansion of support size enlarges the feasible region. On the other hand, the intervals shrink as $N_s$ increases, since this tightens the moment constraints and consequently reduces the feasible region.
The effect of the support size does not appear too sensitive in this example. Thus, taking into account the optimization efficiency, a use of support size of about 45 points appears sufficient.

Figure~\ref{fig:mm13exp} plots the performance when the
supports of the distributions are sampled from the true distributions. The performance trends are similar to Figure~\ref{fig:mm13lgnml}. However, the obtained bounds are slightly looser. Note that Theorem \ref{sample thm} guarantees that the obtained bounds under the generated support points cover the truth with high confidence, when the generating distributions satisfy the heavier-tail condition. In this example, both lognormal and exponential distributions (the latter being the truth) satisfy these conditions and lead to correct bounds. On the other hand, the tightness of the bounds, which is not captured in Theorem \ref{sample thm}, depends on the size and geometry of the feasible region that is determined by a complex interplay between the choice of the uncertainty set and the support-generating distributions. The feasible region using the true exponential distributions include probability weights that are close to uniform weights (since the moment constraints are calibrated using the same distribution). The region using the lognormal, however, does not contain such weights; in fact, when $N_s=500$, the resulting optimizations can be infeasible for $n\leq60$, signaling the need to use more support-generating samples, whereas they are always feasible using the exponential, whose values are shown in the rightmost set of intervals in Figure~\ref{fig:mm13exp}.



The results above are implemented with an initialization that assigns
equal probabilities to the support points. Appendix \ref{sec:multi-start} shows the results applied on different initializations to provide evidence that the formulation in this example has a unique global optimal solution or similar local optimal solutions.

\subsection{Convergence of FWSA and Worst-case Input Distributions} \label{sec:gg1example}
We test the numerical convergence of FWSA. The key parameters in the algorithm are the sample-size growth
rate $\beta$ and the step-size constant $a$. Varying these two parameters, we empirically test the rate of convergence of the FW gap to zero analyzed in Theorem~\ref{rate thm}, and the objective function
$Z(\mbp_k)$ to the true optimal value $Z(\mbp^*)$ analyzed in
Corollary~\ref{rate cor}.
We also investigate the magnitude of the optimal objective value and the form of the identified optimal solution.

Here we consider an $M/G/1$ queue where the arrival process is Poisson known with high accuracy to have rate $1$. On the
other hand, the service time $X_t$ for the $t$-th customer is uncertain but assumed i.i.d.. 
A simulation model is being used to estimate the expected long-run average of the waiting
times $Z(\mbp)=E_{\mbp}[h(\mbX)]$, where 
\[
h(\mathbf X)=\frac 1 T \sum_1^T W_t
\]
and $W_t$ is the waiting time obtained from Lindley's recursion.

We test our FWSA with a KL-divergence-based uncertainty set for $X_t$ as
\begin{equation}
\hat{\mathcal U}=\left\{\mathbf p:\sum_{j=1}^n p_j \log \left(\frac {p_j}{p_{b,j}} \right) \le \eta\right\}\label{mg1prob}
\end{equation}
where $\mathbf p_b=(p_{b,j})_{j=1,\ldots,n}$ is a baseline model chosen to be
a discretized mixture of beta distribution given by $0.3\times\text{Beta}(2,6)+0.7\times\text{Beta}(6,2)$. 
The discrete supports are obtained by uniformly discretizing the interval
$[0,1]$ into $n$ points, i.e. $y_j=(j+1)/n$. 
The set \eqref{mg1prob} provides a good testing ground because steady-state analysis allows obtaining an approximate optimal solution directly which serves as a benchmark for verifying the convergence of our FWSA algorithm (see Appendix \ref{sec:numerics appendix} for further details of this approximate optimal solution).



\begin{figure}[htbp]
  \centering
  \begin{subfigure}[t]{0.45\textwidth}
    \vspace{0pt}
    \includegraphics[width=\textwidth]{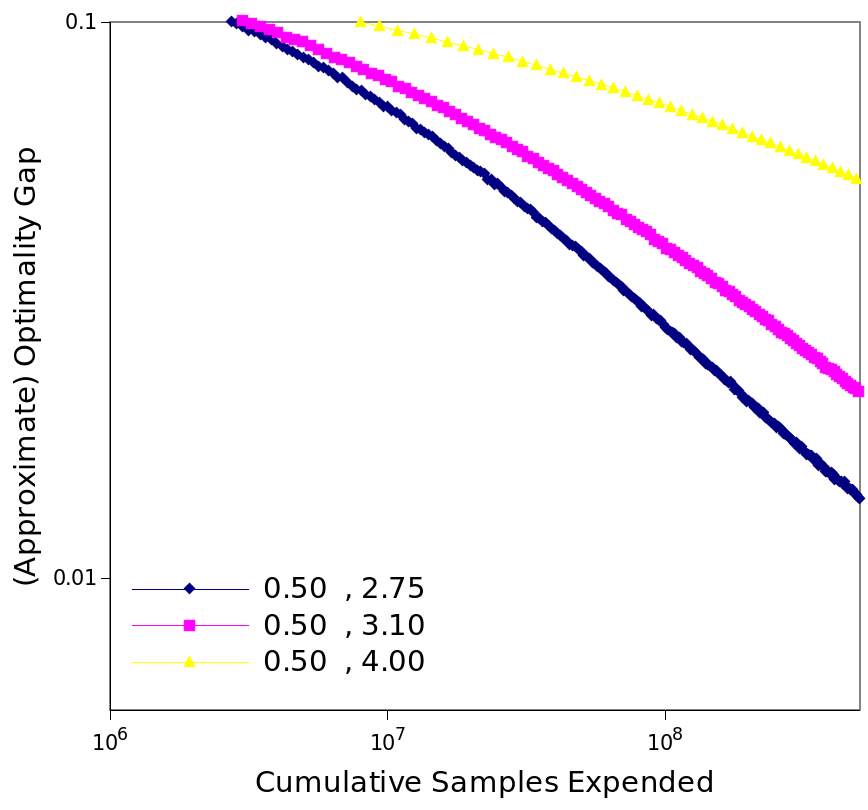}
    \caption{small $a$, $\beta$ varied as shown}
    \label{fig:mg1_small_a}
  \end{subfigure}
\qquad
  \begin{subfigure}[t]{0.45\textwidth}
    \vspace{0pt}
    \includegraphics[width=\textwidth]{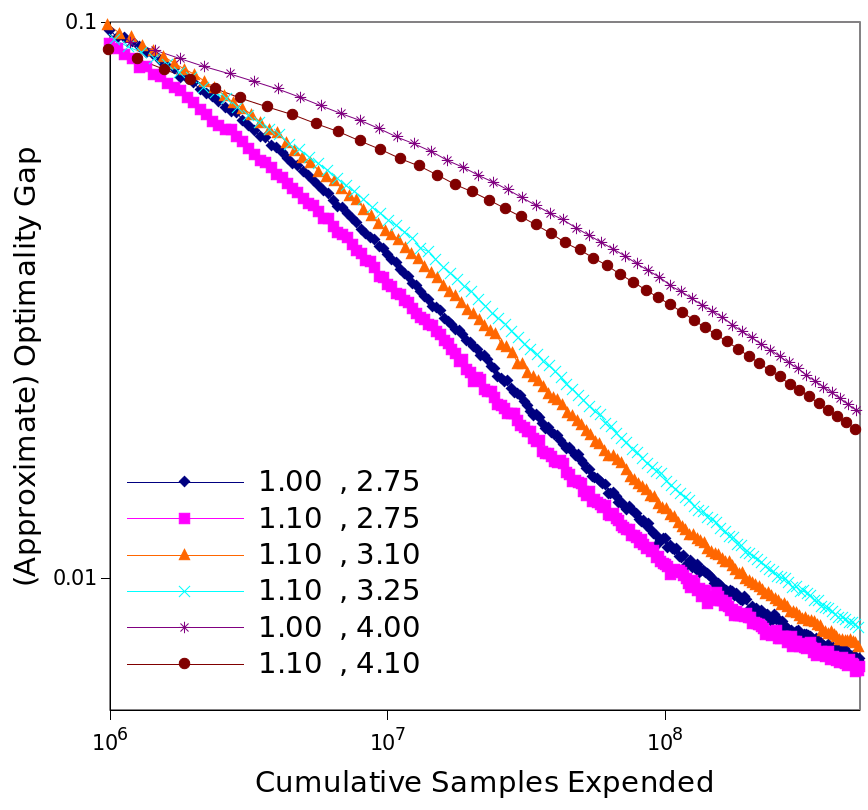}
    \caption{$a=1$, $\beta$ varied as shown}
    \label{fig:mg1_a_1}
  \end{subfigure}

  \begin{subfigure}[b]{0.45\textwidth}
    \includegraphics[width=\textwidth]{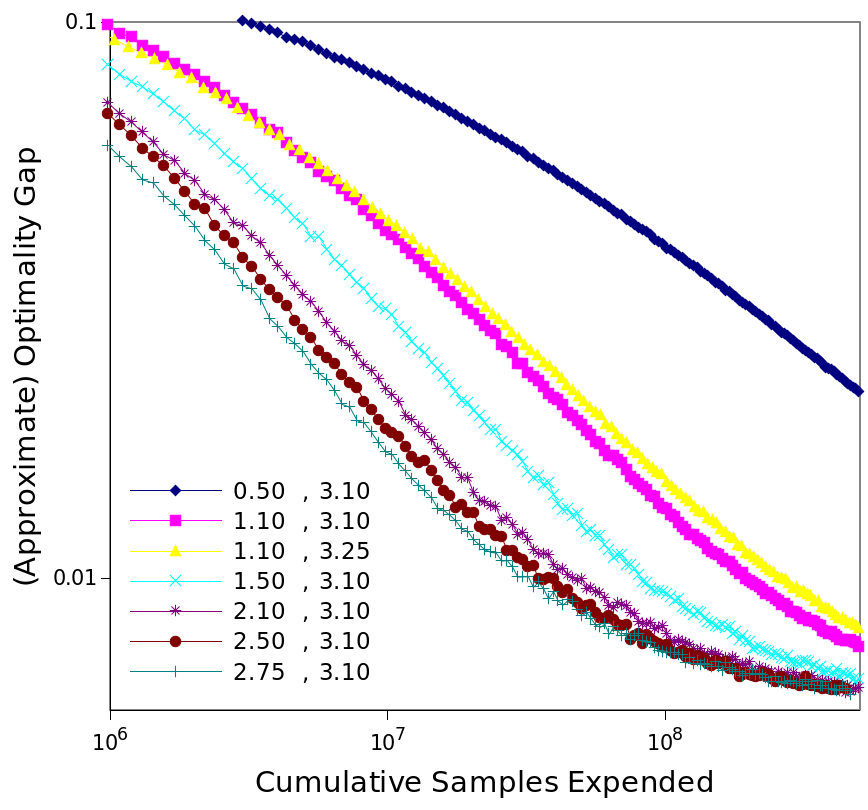}
    \caption{$\beta=3.1$, $a$ varied}
    \label{fig:mg1_beta_3}
  \end{subfigure}
\qquad
  \begin{subfigure}[b]{0.45\textwidth}
    \includegraphics[width=\textwidth]{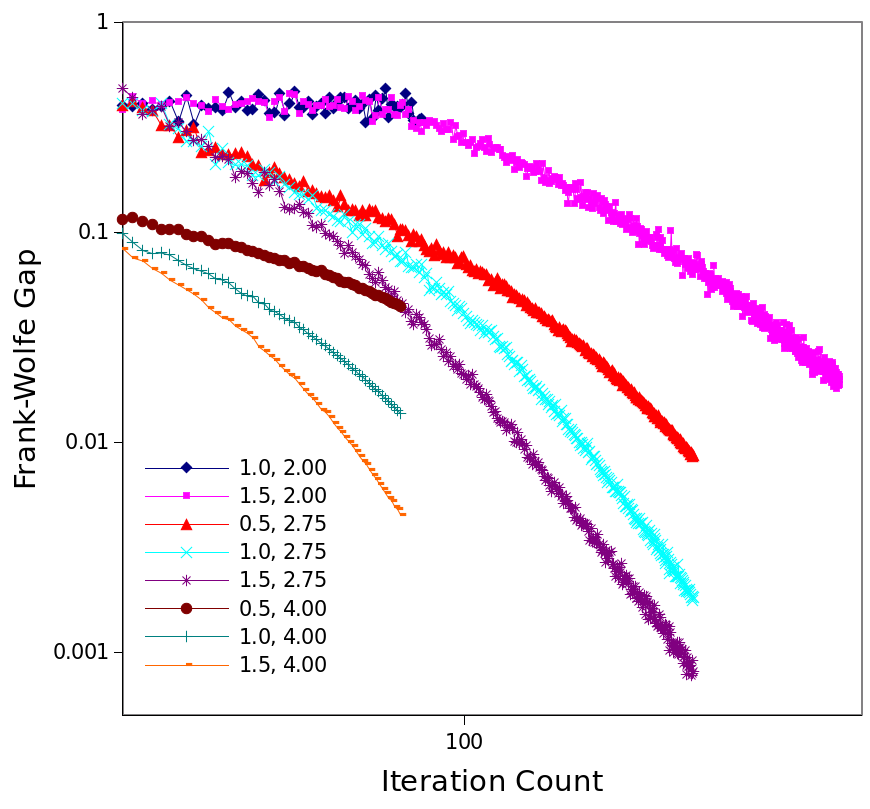}
    \caption{Frank-Wolfe gap vs iteration count}
    \label{fig:mg1fwgap}
  \end{subfigure}
  \caption{Figs~\ref{fig:mg1_small_a},~\ref{fig:mg1_a_1}
    and~\ref{fig:mg1_beta_3} plot the optimality gap of the FWSA algorithm for the $M/G/1$ example as
      function of cumulative simulation samples, under
      various combinations of step-size parameter $a$ and sample-size
      growth parameter $\beta$. The three figures have the same range
      of values in both axes (note the log scale). Fig~\ref{fig:mg1fwgap} shows the
      FW gap as a function of iteration count. All figures provide the legend as $a , \beta$.}
  \label{fig:mg1perf}
\end{figure}

Figure~\ref{fig:mg1perf} captures the
performance of our FWSA algorithm as a function of the $a$ and $\beta$
parameters. Figures~\ref{fig:mg1_small_a}--\ref{fig:mg1_beta_3} plot the
(approximate) optimality gap as a function of the cumulative
simulation replications $W_k$ 
for the maximization problem under~\eqref{mg1prob}. We set the parameters $\eta=0.025$, $n=100$ and
$T=500$. 
Figures~\ref{fig:mg1_small_a},~\ref{fig:mg1_a_1}
and~\ref{fig:mg1_beta_3} provide further insights into the actual
observed finite-sample performance (When interpreting these graphs, note that they are plotted in log-log scale and thus, roughly speaking, the slope of the curve represents the power of the cumulative samples whereas the intercept represents the multiplicative constant in the rate):
\begin{itemize}
\item {\em Fig.~\ref{fig:mg1_small_a}
  v.s.~\ref{fig:mg1_a_1}--\ref{fig:mg1_beta_3}:} Convergence is much slower when $a<1$ no matter the value of $\beta$.
\item {\em Fig.~\ref{fig:mg1_a_1}:} For $a>1$, convergence is again slow if $\beta >4$.
\item {\em Fig.~\ref{fig:mg1_a_1}:} For $a$ slightly greater than $1$, the convergence rates are similar for $\beta\in [2.75, 3.25]$ with better
  performance for the lower end.
\item {\em Fig.~\ref{fig:mg1_beta_3}:} For $\beta=3.1$, the rate of
  convergence generally improves as $a$ increases in the range $[1.10,2.75]$.
\item {\em Figs.~\ref{fig:mg1_small_a},~\ref{fig:mg1_a_1}
  and~\ref{fig:mg1_beta_3}:} The approximation $Z^*_{\infty}$ of the true $Z(\mbp^*)$ (from
  {\bf (SS)} in Appendix \ref{sec:numerics appendix}) has an error of about $0.006$ for
  the chosen $T$, as observed by the leveling off of all plots
  around this value as the sampling effort grows.
\end{itemize}
Figure~\ref{fig:mg1fwgap} shows the FW gap as a function of the
iteration count. In general, the sample paths
with similar  $\beta$ are clustered together, indicating that more effort expended in
estimating the gradient at each iterate leads to a faster drop in the
FW gap per iteration. Within each cluster, performance is
inferior when $a<1$, consistent with Theorem~\ref{rate thm}. Since most runs
terminate when the criterion on the maximum allowed budget of simulation
replications is expended, the end points of the curves indicate that a combination
of $a\geq1$ and a $\beta$ of around $3$ gains the best finite-sample
performance in terms of the FW gap. These choices seem to reconcile with the discussion at the end of Section \ref{sec:local rate} when Assumption \ref{bias} is relaxed to hold for any $\bm\xi_1,\bm\xi_2\in\mathbb R^N$.


We provide further discussion on the shape of the obtained optimal distributions in Appendix \ref{sec:numerics appendix1}.

\section{Conclusion}\label{conclusion}
In this paper we investigated a methodology based on worst-case analysis to quantify input errors in stochastic simulation, by using optimization constraints to represent the partial nonparametric information on the model. The procedure involved a randomized discretization of the support and running FWSA using a gradient estimation technique akin to a nonparametric version of the likelihood ratio or the score function method. We studied the statistical guarantees of the discretization and convergence properties of the proposed FWSA. We also tested our method and verified the theoretical implications on queueing examples. 

We suggest several lines of future research. First is the extension of the methodology to dependent models, such as Markovian inputs or more general time series inputs, which would involve new sets of constraints in the optimizations. Second is the design and analysis of other potential alternate numerical procedures and comparisons with the proposed method. Third is the utilization of the proposed worst-case optimizations in various classes of decision-making problems.

\section*{Acknowledgments}
We thank the Area Editor, the Associate Editor and the three referees for many helpful suggestions that have greatly improved the paper. We gratefully acknowledge support from the National Science Foundation under grants CMMI-1542020, CMMI-1523453 and CAREER CMMI-1653339.

\bibliographystyle{ormsv080}
\bibliography{bibliography}

\ECSwitch


\ECHead{Appendix}

\section{Technical Proofs}\label{sec:proofs}
\proof{Proof of Theorem \ref{sample thm}.}
Let $\delta(y)$ be the delta measure at $y$. For each $i=1,\ldots,m$, define
$$\tilde P^i=\sum_{j=1}^{n^i}\frac{L^i(y_j^i)}{\sum_{r=1}^{n^i}L^i(y_r^i)}\delta(y_j^i)$$
i.e., the distribution with point mass $L^i(y_j^i)/\sum_{r=1}^{n^i}L^i(y_r^i)$ on each $y_j^i$, where $L^i=dP_0^i/dQ^i$. We first show that as $n\to\infty$, the solution $(\tilde P^i)_{i=1,\ldots,m}$ is feasible for the optimization problems in \eqref{sample counterpart} in an appropriate sense.

Consider Case 1. For each $l=1,\ldots,s^i$, by a change measure we have $E_{Q^i}|f_l^i(X^i)L(X^i)|=E_{P_0^i}|f_l^i(X^i)|<\infty$ by our assumption. Also note that $E_{Q^i}L^i=1$. Therefore, by the law of large numbers,
$$E_{\tilde P^i}[f_l^i(X^i)]=\frac{\sum_{j=1}^{n^i}L^i(y_j^i)f_l^i(y_j^i)}{\sum_{j=1}^{n^i}L^i(y_j^i)}=\frac{(1/n^i)\sum_{j=1}^{n^i}L^i(y_j^i)f_l^i(y_j^i)}{(1/n^i)\sum_{j=1}^{n^i}L^i(y_j^i)}\to E_{Q^i}[f_l^i(X^i)L(X^i)]\text{\ \ a.s.}$$
Since $E_{Q^i}[f_l^i(X^i)L(X^i)]=E_{P_0^i}[f_l^i(X^i)]<\mu_j^i$ by our assumption, we have $E_{\tilde P^i}[f_l^i(X^i)]\leq\mu_l^i$ eventually as $n^i\to\infty$.

Consider Case 2. We have
\begin{eqnarray*}
d_\phi(\tilde P^i,\hat P_b^i)
&=&\sum_{j=1}^{n^i}\phi\left(\frac{L^i(y_j^i)/\sum_{r=1}^{n^i}L^i(y_r^i)}{L_b^i(y_j^i)/\sum_{r=1}^{n^i}L_b^i(y_r^i)}\right)\frac{L_b^i(y_j^i)}{\sum_{r=1}^{n^i}L_b^i(y_r^i)}\notag\\
&=&\frac{1}{n^i}\sum_{j=1}^{n^i}\phi\left(\tilde L^i(y_j^i)\frac{(1/n^i)\sum_{r=1}^{n^i}L_b^i(y_r^i)}{(1/n^i)\sum_{r=1}^{n^i}L^i(y_r^i)}\right)\frac{L_b^i(y_j^i)}{(1/n^i)\sum_{r=1}^{n^i}L_b^i(y_r^i)}
\end{eqnarray*}
where $\tilde L^i=dP_0^i/dP_b^i$. Consider, for a given $\epsilon>0$,
\begin{eqnarray}
&&P(|d_\phi(\tilde P^i,\hat P_b^i)-d_\phi(P_0^i,P_b^i)|>\epsilon)\notag\\
&=&P\left(\left|\frac{1}{n^i}\sum_{j=1}^{n^i}\phi\left(\tilde L^i(y_j^i)\frac{(1/n^i)\sum_{r=1}^{n^i}L_b^i(y_r^i)}{(1/n^i)\sum_{r=1}^{n^i}L^i(y_r^i)}\right)\frac{L_b^i(y_j^i)}{(1/n^i)\sum_{r=1}^{n^i}L_b^i(y_r^i)}-d_\phi(P_0^i,P_b^i)\right|>\epsilon\right)\notag\\
&\leq&P\left(\left|\frac{1}{n^i}\sum_{j=1}^{n^i}\phi\left(\tilde L^i(y_j^i)\frac{(1/n^i)\sum_{r=1}^{n^i}L_b^i(y_r^i)}{(1/n^i)\sum_{r=1}^{n^i}L^i(y_r^i)}\right)\frac{L_b^i(y_j^i)}{(1/n^i)\sum_{r=1}^{n^i}L_b^i(y_r^i)}-\frac{1}{n^i}\sum_{j=1}^{n^i}\phi(\tilde L^i(y_j^i))\frac{L_b^i(y_j^i)}{(1/n^i)\sum_{r=1}^{n^i}L_b^i(y_r^i)}\right|>\frac{\epsilon}{2}\right){}\notag\\
&&+P\left(\left|\frac{1}{n^i}\sum_{j=1}^{n^i}\phi(\tilde L^i(y_j^i))\frac{L_b^i(y_j^i)}{(1/n^i)\sum_{r=1}^{n^i}L_b^i(y_r^i)}-d_\phi(P_0^i,P_b^i)\right|>\frac{\epsilon}{2}\right)\label{interim revised}
\end{eqnarray}
We analyze the two terms in \eqref{interim revised}. For any sufficiently small $\lambda>0$, the first term is bounded from above by
\begin{eqnarray}
&&P\Bigg(\left|\frac{1}{n^i}\sum_{j=1}^{n^i}\phi\left(\tilde L^i(y_j^i)\frac{(1/n^i)\sum_{r=1}^{n^i}L_b^i(y_r^i)}{(1/n^i)\sum_{r=1}^{n^i}L^i(y_r^i)}\right)\frac{L_b^i(y_j^i)}{(1/n^i)\sum_{r=1}^{n^i}L_b^i(y_r^i)}-\frac{1}{n^i}\sum_{j=1}^{n^i}\phi(\tilde L^i(y_j^i))\frac{L_b^i(y_j^i)}{(1/n^i)\sum_{r=1}^{n^i}L_b^i(y_r^i)}\right|>\frac{\epsilon}{2}{}\notag\\
&&{};\left|\frac{1}{n^i}\sum_{r=1}^{n^i}L^i(y_r^i)-1\right|\leq\lambda,\left|\frac{1}{n^i}\sum_{r=1}^{n^i}L_b^i(y_r^i)-1\right|\leq\lambda\Bigg){}\notag\\
&&{}+P\Bigg(\left|\frac{1}{n^i}\sum_{j=1}^{n^i}\phi\left(\tilde L^i(y_j^i)\frac{(1/n^i)\sum_{r=1}^{n^i}L_b^i(y_r^i)}{(1/n^i)\sum_{r=1}^{n^i}L^i(y_r^i)}\right)\frac{L_b^i(y_j^i)}{(1/n^i)\sum_{r=1}^{n^i}L_b^i(y_r^i)}-\frac{1}{n^i}\sum_{j=1}^{n^i}\phi(\tilde L^i(y_j^i))\frac{L_b^i(y_j^i)}{(1/n^i)\sum_{r=1}^{n^i}L_b^i(y_r^i)}\right|>\frac{\epsilon}{2}{}\notag\\
&&{};\left|\frac{1}{n^i}\sum_{r=1}^{n^i}L^i(y_r^i)-1\right|>\lambda\text{\ or\ }\left|\frac{1}{n^i}\sum_{r=1}^{n^i}L_b^i(y_r^i)-1\right|>\lambda\Bigg)\notag\\
&\leq&P\left(\frac{1}{n^i}\sum_{j=1}^{n^i}(|\phi(\tilde L^i(y_j^i))|+1)O(\lambda)\frac{L_b^i(y_j^i)}{(1/n^i)\sum_{r=1}^{n^i}L_b^i(y_r^i)}>\frac{\epsilon}{2};\left|\frac{1}{n^i}\sum_{r=1}^{n^i}L^i(y_r^i)-1\right|\leq\lambda,\left|\frac{1}{n^i}\sum_{r=1}^{n^i}L_b^i(y_r^i)-1\right|\leq\lambda\right){}\notag\\
&&{}+P\left(\left|\frac{1}{n^i}\sum_{r=1}^{n^i}L^i(y_r^i)-1\right|>\lambda\text{\ or\ }\left|\frac{1}{n^i}\sum_{r=1}^{n^i}L_b^i(y_r^i)-1\right|>\lambda\right)\label{interim revised1}
\end{eqnarray}
where the first term in the last inequality follows from the continuity condition on $\phi$, with $O(\lambda)$ being a deterministic positive function of $\lambda$ that converges to 0 as $\lambda\to0$. This first term is further bounded from above by
\begin{equation}
P\left(\frac{1}{n^i}\sum_{j=1}^{n^i}(|\phi(\tilde L^i(y_j^i))|+1)\frac{L_b^i(y_j^i)}{(1/n^i)\sum_{r=1}^{n^i}L_b^i(y_r^i)}O(\lambda)>\frac{\epsilon}{2}\right)\label{interim revised2}
\end{equation}
By the law of large numbers, we have
$$\frac{1}{n^i}\sum_{j=1}^{n^i}(|\phi(\tilde L^i(y_j^i))|+1)L_b^i(y_j^i)\to E_{Q^i}[(|\phi(\tilde L^i(X^i))|+1)L_b^i(X^i)]=E_{P_b^i}|\phi(\tilde L^i(X^i))|+1\text{\ \ a.s.}$$
by using our assumption $E_{P_b^i}|\phi(\tilde L^i(X^i))|<\infty$. Moreover, by the law of large numbers again, we have $(1/n^i)\sum_{r=1}^{n^i}L_b^i(y_r^i)\to1$ a.s.. Thus,
$$\frac{1}{n^i}\sum_{j=1}^{n^i}(|\phi(\tilde L^i(y_j^i))|+1)\frac{L_b^i(y_j^i)}{(1/n^i)\sum_{r=1}^{n^i}L_b^i(y_r^i)}\to E_{P_b^i}|\phi(\tilde L^i(X^i))|+1\text{\ \ a.s.}$$
When $\lambda$ is chosen small enough relative to $\epsilon/2$, we have \eqref{interim revised2} go to 0 as $n^i\to\infty$.

Since both $\frac{1}{n^i}\sum_{r=1}^{n^i}L^i(y_r^i)\to1$ and $\frac{1}{n^i}\sum_{r=1}^{n^i}L_b^i(y_r^i)\to1$ a.s., the second term in \eqref{interim revised1} also goes to 0 as $n^i\to\infty$. This concludes that the first term in \eqref{interim revised} goes to 0.

For the second term in \eqref{interim revised}, note that
$$\frac{1}{n^i}\sum_{j=1}^{n^i}\phi(\tilde L^i(y_j^i))L_b^i(y_j^i)\to E_{Q^i}[\phi(\tilde L^i(X^i))L_b^i(X^i)]=E_{P_b^i}[\phi(\tilde L^i(X^i))]=d_\phi(P_0^i,P_b^i)\text{\ \ a.s.}$$
by the law of large numbers and the assumption that $E_{P_b^i}|\phi(\tilde L^i(X^i))|<\infty$. Moreover, since $(1/n^i)\sum_{r=1}^{n^i}L_b^i(y_r^i)\to1$, we get
$$\frac{1}{n^i}\sum_{j=1}^{n^i}\phi(\tilde L^i(y_j^i))\frac{L_b^i(y_j^i)}{(1/n^i)\sum_{r=1}^{n^i}L_b^i(y_r^i)}\to d_\phi(P_0^i,P_b^i)\text{\ \ a.s.}$$
Thus, the second term in \eqref{interim revised} goes to 0 as $n^i\to\infty$. Therefore, we conclude that $d_\phi(\tilde P^i,\hat P_b^i)\stackrel{p}{\to}d_\phi(P_0^i,P_b^i)$. Since $d_\phi(P_0^i,P_b^i)<\eta^i$ by our assumption, we have $P(d_\phi(\tilde P^i,\hat P_b^i)\leq\eta^i)\to1$ as $n^i\to\infty$.

Next we consider the objective in \eqref{sample counterpart}. We show that $Z(\tilde P^1,\ldots,\tilde P^m)-Z(P_0^1,\ldots,P_0^m)=O_p(1/\sqrt n)$, following the argument in the theory of differentiable statistical functionals (e.g., \cite{serfling2009approximation}, Chapter 6). For any $\lambda$ between 0 and 1, we write
\begin{eqnarray*}
&&Z(P_0^1+\lambda(\tilde P^1-P_0^1),\ldots,P_0^m+\lambda(\tilde P^m-P_0^m))\\
&=&\int\cdots\int h(\mathbf x^1,\ldots,\mathbf x^m)\prod_{i=1}^m\prod_{t=1}^{T^i}d[P_0^i+\lambda(\tilde P^i-P_0^i)](x_t^i)\\
&=&\sum_{k=0}^T\lambda^k\sum_{u\in\mathcal I^k}\int\cdots\int h(\mathbf x^1,\ldots,\mathbf x^m)\prod_{(i,t)\in(\mathcal S_u^k)^c}dP_0^i(x_t^i)\prod_{(i,t)\in\mathcal S_u^k}d(\tilde P^i-P_0^i)(x_t^i)
\end{eqnarray*}
where $\{\mathcal S_u^k\}_{u\in\mathcal I^k}$ is the collection of all subsets of $\{(i,t):i=1,\ldots,m,t=1,\ldots,T^i\}$ with cardinality $k$, and $\mathcal I^k$ indexes all these subsets. Note that
\begin{eqnarray}
&&\frac{d}{d\lambda}Z(P_0^1+\lambda(\tilde P^1-P_0^1),\ldots,P_0^m+\lambda(\tilde P^m-P_0^m))\Bigg|_{\lambda=0^+}\notag\\
&=&\sum_{i=1}^m\sum_{t=1}^{T^i}\int\cdots\int h(\mathbf x^1,\ldots,\mathbf x^m)\prod_{(j,s):(j,s)\neq(i,t)}dP_0^j(x_s^j)d(\tilde P^i-P_0^i)(x_t^i)\notag\\
&=&\sum_{i=1}^m\int\varphi^i(x;P_0^1,\ldots,P_0^m)d(\tilde P^i-P_0^i)(x)
\label{interim new1}
\end{eqnarray}
where
\begin{equation}
\varphi^i(x;P_0^1,\ldots,P_0^m)=\sum_{t=1}^{T^i}E_{P_0^1,\ldots,P_0^m}[h(\mathbf X^1,\ldots,\mathbf X^m)|X_t^i=x]\label{interim new91}
\end{equation}
By the definition of $L^i$, we can write \eqref{interim new1} as
\begin{eqnarray}
&&\sum_{i=1}^m\left(\frac{\int\varphi^i(x;P_0^1,\ldots,P_0^m)L^i(x)d\hat Q^i(x)}{(1/n^i)\sum_{j=1}^{n^i}L^i(y_j^i)}-\int\varphi^i(x;P_0^1,\ldots,P_0^m)L^i(x)dQ^i(x)\right)\notag\\
&=&\sum_{i=1}^m\left(\frac{\int\varphi^i(x;P_0^1,\ldots,P_0^m)L^i(x)d(\hat Q^i-Q^i)(x)}{(1/n^i)\sum_{j=1}^{n^i}L^i(y_j^i)}-\int\varphi^i(x;P_0^1,\ldots,P_0^m)L^i(x)dQ^i(x)\left(1-\frac{1}{(1/n^i)\sum_{j=1}^{n^i}L^i(y_j^i)}\right)\right)
\label{update interim2}
\end{eqnarray}
where $\hat Q^i$ is the empirical distribution $(1/n^i)\sum_{j=1}^{n^i}\delta(y_j^i)$ on the $n^i$ observations generated from $Q^i$.

Suppose $\varphi^i(x;P_0^1,\ldots,P_0^m)L^i(x)=0$ a.s., then $\int\varphi^i(x;P_0^1,\ldots,P_0^m)L^i(x)d(\hat Q^i-Q^i)(x)=0$ a.s.. Otherwise, using the assumed boundedness of $h$, hence $\varphi^i(x;P_0^1,\ldots,P_0^m)$, and $L^i$, we have, by the central limit theorem,
$$\sqrt{n^i}\left(\int\varphi^i(x;P_0^1,\ldots,P_0^m)L^i(x)d(\hat Q^i-Q^i)(x)\right)\Rightarrow N(0,(\sigma^i)^2)$$
where $(\sigma^i)^2=Var_{Q^i}(\varphi^i(X^i;P_0^1,\ldots,P_0^m)L^i(X^i))>0$ is finite. Since $(1/n^i)\sum_{j=1}^{n^i}L^i(y_j^i)\to1$ a.s. by the law of large numbers, and that $\int\varphi^i(x;P_0^1,\ldots,P_0^m)L^i(x)d\hat Q^i(x)$ is bounded, the second term in \eqref{update interim2} converges to 0 a.s.. Thus, by Slutsky's theorem, each summand in \eqref{update interim2} converges in distribution to $N(0,(\sigma^i)^2)$. Since for each $i$ we have $n^i=nw^i$ for some fixed $w^i>0$, we conclude that \eqref{update interim2} equal $O_p(1/\sqrt n)$.

Now consider
\begin{eqnarray}
&&\frac{d^2}{d\lambda^2}Z(P_0^1+\lambda(\tilde P^1-P_0^1),\ldots,P_0^m+\lambda(\tilde P^m-P_0^m))\notag\\
&=&\sum_{k=2}^Tk(k-1)\lambda^{k-2}\sum_{u\in\mathcal I^k}\int\cdots\int h(\mathbf x^1,\ldots,\mathbf x^m)\prod_{(i,t)\in(\mathcal S_u^k)^c}dP_0^i(x_t^i)\prod_{(i,t)\in\mathcal S_u^k}d(\tilde P^i-P_0^i)(x_t^i)\label{update interim3}
\end{eqnarray}
Fixing each $\mathcal S_u^k$, we define
$$h_{\mathcal S_u^k}(\mathbf x_{\mathcal S_u^k})=\int\cdots\int h(\mathbf x^1,\ldots,\mathbf x^m)\prod_{(i,t)\in(\mathcal S_u^k)^c}dP_0^i(x_t^i)$$
where $\mathbf x_{\mathcal S_u^k}=(x_t^i)_{(i,t)\in\mathcal S_u^k}$. Next define
\begin{eqnarray*}
&&\tilde h_{\mathcal S_u^k}(\mathbf x_{\mathcal S_u^k})\\
&=&h_{\mathcal S_u^k}(\mathbf x_{\mathcal S_u^k})-\sum_{(j,t)\in\mathcal S_u^k}\int h_{\mathcal S_u^k}(\mathbf x_{\mathcal S_u^k})dP_0^j(x_t^j)+\sum_{(j_1,t_1),(j_2,t_2)\in\mathcal S_u^k}\int\int h_{\mathcal S_u^k}(\mathbf x_{\mathcal S_u^k})dP_0^{j_1}(x_{t_1}^{j_1})dP_0^{j_2}(x_{t_2}^{j_2})-\cdots{}\\
&&{}+(-1)^k\int\cdots\int h_{\mathcal S_u^k}(\mathbf x_{\mathcal S_u^k})dP_0^{j_1}(x_{t_1}^{j_1})\cdots dP_0^{j_k}(x_{t_k}^{j_k})
\end{eqnarray*}
where each summation above is over the set of all possible combinations of $(j,t)\in\mathcal S_u^k$ with increasing size. Direct verification shows that $\tilde h_{\mathcal S_u^k}$ has the property that
$$\int\cdots\int\tilde h_{\mathcal S_u^k}(\mathbf x_{\mathcal S_u^k})\prod_{(i,t)\in\mathcal S_u^k}dR^j(x_t^j)=\int\cdots\int h_{\mathcal S_u^k}(\mathbf x_{\mathcal S_u^k})\prod_{(i,t)\in\mathcal S_u^k}d(R^j(x_t^j)-P_0^j(x_t^j))$$
for any probability measures $R^j$'s, and
\begin{equation}
\int\tilde h_{\mathcal S_u^k}(\mathbf x_{\mathcal S_u^k})dP_0^j(x_t^j)=0\label{update interim4}
\end{equation}
for any $(j,t)\in\mathcal S_u^k$. Thus, \eqref{update interim3} is equal to
$$\sum_{k=2}^Tk(k-1)\lambda^{k-2}\sum_{u\in\mathcal I^k}\int\tilde h_{\mathcal S_u^k}(\mathbf x_{\mathcal S_u^k})\prod_{(i,t)\in\mathcal S_u^k}d\tilde P^i(x_t^i)$$
Now, viewing $\tilde P^i$ as randomly generated from $Q^i$, consider
\begin{eqnarray}
&&E_{Q^1,\ldots,Q^m}\left(\sum_{k=2}^Tk(k-1)\lambda^{k-2}\sum_{u\in\mathcal I^k}\int\tilde h_{\mathcal S_u^k}(\mathbf x_{\mathcal S_u^k})\prod_{(i,t)\in\mathcal S_u^k}d\tilde P^i(x_t^i)\right)^2\notag\\
&=&E_{Q^1,\ldots,Q^m}\left[\left(\sum_{k=2}^Tk(k-1)\lambda^{k-2}\sum_{u\in\mathcal I^k}\int\tilde h_{\mathcal S_u^k}(\mathbf x_{\mathcal S_u^k})\prod_{(i,t)\in\mathcal S_u^k}d\tilde P^i(x_t^i)\right)^2;\frac{1}{n^i}\sum_{r=1}^{n^i}L^{i}(Y_r^i)\geq1-\epsilon\text{\ for all\ }i=1,\ldots,m\right]{}\notag\\
&&{}+E_{Q^1,\ldots,Q^m}\left[\left(\sum_{k=2}^Tk(k-1)\lambda^{k-2}\sum_{u\in\mathcal I^k}\int\tilde h_{\mathcal S_u^k}(\mathbf x_{\mathcal S_u^k})\prod_{(i,t)\in\mathcal S_u^k}d\tilde P^i(x_t^i)\right)^2;\frac{1}{n^i}\sum_{r=1}^{n^i}L^{i}(Y_r^i)<1-\epsilon\text{\ for some\ }i=1,\ldots,m\right]\label{interim revised3}
\end{eqnarray}
We analyze the two terms in \eqref{interim revised3}. Note that the first term can be written as
\begin{eqnarray}
&&E_{Q^1,\ldots,Q^m}\Bigg[\left(\sum_{k=2}^Tk(k-1)\lambda^{k-2}\sum_{u\in\mathcal I^k}\frac{1}{n^{i_1}n^{i_2}\cdots n^{i_k}}\sum_{j_1=1}^{n^{i_1}}\cdots\sum_{j_k=1}^{n^{i_k}}\tilde h_{\mathcal S_u^k}(Y_{j_1}^{i_1},\ldots,Y_{j_k}^{i_k})\frac{L^{i_1}(Y_{j_1}^{i_1})\cdots L^{i_k}(Y_{j_k}^{i_k})}{\prod_{s=1}^k((1/n^{i_s})\sum_{r=1}^{n^{i_s}}L^{i_s}(Y_r^{i_s}))}\right)^2;{}\notag\\
&&{}\frac{1}{n^i}\sum_{r=1}^{n^i}L^{i}(Y_r^i)\geq1-\epsilon\text{\ for all\ }i=1,\ldots,m\Bigg]\notag\\
&\leq&\Bigg(\sum_{k=2}^Tk(k-1)\lambda^{k-2}\sum_{u\in\mathcal I^k}\frac{1}{n^{i_1}n^{i_2}\cdots n^{i_k}}\Bigg(E_{Q^1,\ldots,Q^m}\Bigg[\Bigg(\sum_{j_1=1}^{n^{i_1}}\cdots\sum_{j_k=1}^{n^{i_k}}\tilde h_{\mathcal S_u^k}(Y_{j_1}^{i_1},\ldots,Y_{j_k}^{i_k}){}\notag\\
&&{}\frac{L^{i_1}(Y_{j_1}^{i_1})\cdots L^{i_k}(Y_{j_k}^{i_k})}{\prod_{s=1}^k((1/n^{i_s})\sum_{r=1}^{n^{i_s}}L^{i_s}(Y_r^{i_s}))}\Bigg)^2;\frac{1}{n^i}\sum_{r=1}^{n^i}L^{i}(Y_r^i)\geq1-\epsilon\text{\ for all\ }i=1,\ldots,m\Bigg]\Bigg)^{1/2}\Bigg)^2\label{interim revised5}
\end{eqnarray}
by Minkowski's inequality, where we view $Y_j^i$'s as the random variables constituting the observations generated from $Q^i$'s. Since the expression $\prod_{s=1}^k((1/n^{i_s})\sum_{r=1}^{n^{i_s}}L^{i_s}(Y_r^{i_s}))$ inside the expectation in \eqref{interim revised5} does not depend on the $j_s$'s, \eqref{interim revised5} is further bounded from above by
\begin{eqnarray}
&&\Bigg(\sum_{k=2}^Tk(k-1)\lambda^{k-2}\sum_{u\in\mathcal I^k}\frac{1}{n^{i_1}n^{i_2}\cdots n^{i_k}}\left(E_{Q^1,\ldots,Q^m}\left(\sum_{j_1=1}^{n^{i_1}}\cdots\sum_{j_k=1}^{n^{i_k}}\tilde h_{\mathcal S_u^k}(Y_{j_1}^{i_1},\ldots,Y_{j_k}^{i_k})\frac{L^{i_1}(Y_{j_1}^{i_1})\cdots L^{i_k}(Y_{j_k}^{i_k})}{(1-\epsilon)^k}\right)^2\right)^{1/2}\Bigg)^2\notag\\
&=&\Bigg(\sum_{k=2}^T\frac{k(k-1)\lambda^{k-2}}{(1-\epsilon)^k}\sum_{u\in\mathcal I^k}\frac{1}{n^{i_1}n^{i_2}\cdots n^{i_k}}\Bigg(\sum_{j_1=1}^{n^{i_1}}\cdots\sum_{j_k=1}^{n^{i_k}}\sum_{j_1'=1}^{n^{i_1}}\cdots\sum_{j_k'=1}^{n^{i_k}}
E_{Q^1,\ldots,Q^m}\Bigg[\tilde h_{\mathcal S_u^k}(Y_{j_1}^{i_1},\ldots,Y_{j_k}^{i_k})L^{i_1}(Y_{j_1}^{i_1})\cdots L^{i_k}(Y_{j_k}^{i_k}){}\notag\\
&&{}\tilde h_{\mathcal S_u^k}(Y_{j_1'}^{i_1},\ldots,Y_{j_k'}^{i_k})L^{i_1}(Y_{j_1'}^{i_1})\cdots L^{i_k}(Y_{j_k'}^{i_k})\Bigg]\Bigg)^{1/2}\Bigg)^2\label{update interim6}
\end{eqnarray}
Note that
\begin{equation}
E_{Q^1,\ldots,Q^m}\left[\tilde h_{\mathcal S_u^k}(Y_{j_1}^{i_1},\ldots,Y_{j_k}^{i_k})L^{i_1}(Y_{j_1}^{i_1})\cdots L^{i_k}(Y_{j_k}^{i_k})\tilde h_{\mathcal S_u^k}(Y_{j_1'}^{i_1},\ldots,Y_{j_k'}^{i_k})L^{i_1}(Y_{j_1'}^{i_1})\cdots L^{i_k}(Y_{j_k'}^{i_k})\right]=0\label{update interim5}
\end{equation}
if any $Y_j^i$ shows up only once among all those in both $\tilde h_{\mathcal S_u^k}(Y_{j_1}^{i_1},\ldots,Y_{j_k}^{i_k})$ and $\tilde h_{\mathcal S_u^k}(Y_{j_1'}^{i_1},\ldots,Y_{j_k'}^{i_k})$ in the expectation. To see this, suppose without loss of generality that $Y_{j_1}^{i_1}$ appears only once. Then we have
\begin{eqnarray*}
&&E_{Q^1,\ldots,Q^m}\left[\tilde h_{\mathcal S_u^k}(Y_{j_1}^{i_1},\ldots,Y_{j_k}^{i_k})L^{i_1}(Y_{j_1}^{i_1})\cdots L^{i_k}(Y_{j_k}^{i_k})\tilde h_{\mathcal S_u^k}(Y_{j_1'}^{i_1},\ldots,Y_{j_k'}^{i_k})L^{i_1}(Y_{j_1'}^{i_1})\cdots L^{i_k}(Y_{j_k'}^{i_k})\right]\\
&=&E_{Q^1,\ldots,Q^m}\Big[E_{Q^1,\ldots,Q^m}\left[\tilde h_{\mathcal S_u^k}(Y_{j_1}^{i_1},\ldots,Y_{j_k}^{i_k})L^{i_1}(Y_{j_1}^{i_1})\Big|Y_j^{i_2},\ldots,Y_{j_k}^{i_k},Y_{j_1'}^{i_1},\ldots,Y_{j_k'}^{i_k}\right]L^{i_2}(Y_j^{i_2})\cdots L^{i_k}(Y_{j_k}^{i_k}){}\\
&&{}\tilde h_{\mathcal S_u^k}(Y_{j_1'}^{i_1},\ldots,Y_{j_k'}^{i_k})L^{i_1}(Y_{j_1'}^{i_1})\cdots L^{i_k}(Y_{j_k'}^{i_k})\big]\\
&=&E_{Q^1,\ldots,Q^m}\Big[E_{P_0^{i_1}}\left[\tilde h_{\mathcal S_u^k}(Y_{j_1}^{i_1},\ldots,Y_{j_k}^{i_k})\Big|Y_j^{i_2},\ldots,Y_{j_k}^{i_k},Y_{j_1'}^{i_1},\ldots,Y_{j_k'}^{i_k}\right]L^{i_2}(Y_j^{i_2})\cdots L^{i_k}(Y_{j_k}^{i_k}){}\\
&&{}\tilde h_{\mathcal S_u^k}(Y_{j_1'}^{i_1},\ldots,Y_{j_k'}^{i_k})L^{i_1}(Y_{j_1'}^{i_1})\cdots L^{i_k}(Y_{j_k'}^{i_k})\Big]\\
&=&0
\end{eqnarray*}
since $E_{P_0^{i_1}}\left[\tilde h_{\mathcal S_u^k}(Y_{j_1}^{i_1},\ldots,Y_{j_k}^{i_k})\Big|Y_j^{i_2},\ldots,Y_{j_k}^{i_k},Y_{j_1'}^{i_1},\ldots,Y_{j_k'}^{i_k}\right]=0$ by \eqref{update interim4}.

The observation in \eqref{update interim5} implies that the summation in \eqref{update interim6}
$$\sum_{j_1=1}^{n^{i_1}}\cdots\sum_{j_k=1}^{n^{i_k}}\sum_{j_1'=1}^{n^{i_1}}\cdots\sum_{j_k'=1}^{n^{i_k}}E_{Q^1,\ldots,Q^m}\left[\tilde h_{\mathcal S_u^k}(Y_{j_1}^{i_1},\ldots,Y_{j_k}^{i_k})L^{i_1}(Y_{j_1}^{i_1})\cdots L^{i_k}(Y_{j_k}^{i_k})\tilde h_{\mathcal S_u^k}(Y_{j_1'}^{i_1},\ldots,Y_{j_k'}^{i_k})L^{i_1}(Y_{j_1'}^{i_1})\cdots L^{i_k}(Y_{j_k'}^{i_k})\right]$$
contains only $O(n^k)$ non-zero summands. This is because in each non-zero summand only at most $k$ distinct $Y_{j}^{i}$'s can be present inside the expectation, and the cardinality of such combinations is $O(n^k)$. Note that each summand is bounded since $h$, hence $\tilde h_{\mathcal S_u^k}$, and $L^i$ are all bounded by our assumptions.
  Hence \eqref{update interim6} is
\begin{equation}
\left(\sum_{k=2}^T\frac{k(k-1)\lambda^{k-2}}{(1-\epsilon)^k}\binom{T}{k}O\left(\frac{1}{n^{k/2}}\right)\right)^2=O\left(\frac{1}{n^2}\right)\label{update interim7}
\end{equation}
This shows that \eqref{update interim3} is $O_p(1/n)$ for any $\lambda$ between 0 and 1. Therefore, by using Taylor's expansion, and the conclusion that \eqref{update interim2} is $O_p(1/\sqrt n)$, we have
\begin{equation}
Z(\tilde P^1,\ldots,\tilde P^m)=Z(P_0^1,\ldots,P_0^m)+O_p\left(\frac{1}{\sqrt n}\right)=Z_0+O_p\left(\frac{1}{\sqrt n}\right)\label{interim revised4}
\end{equation}

Note that we have shown previously that $P(\tilde P^i\in\hat{\mathcal U}^i)\to1$ for any $i=1,\ldots,m$ in both Cases 1 and 2. Using this and \eqref{interim revised4}, for any given $\epsilon>0$, we can choose $M,N>0$ big enough such that
$$P(\sqrt n(\hat Z_*-Z_0)>M)\leq P(|\sqrt n(Z(\tilde P^1,\ldots,\tilde P^i)-Z_0)|>M)+\sum_{i=1}^mP(\tilde P^i\notin\hat{\mathcal U}^i)<\epsilon$$
and similarly
$$P(\sqrt n(Z_0-\hat Z^*)>M)\leq P(|\sqrt n(Z(\tilde P^1,\ldots,\tilde P^i)-Z_0)|>M)+\sum_{i=1}^mP(\tilde P^i\notin\hat{\mathcal U}^i)<\epsilon$$
for any $n>N$. This concludes that
$$\hat Z_*\leq Z_0+O_p\left(\frac{1}{\sqrt n}\right)\leq\hat Z^*$$
\endproof

\proof{Proof of Theorem \ref{prop:gradient}.}
To prove 1., consider first a mixture of $\mathbf p^i=(p_j^i)_{j=1,\ldots,n^i}$ with an arbitrary $\mathbf q^i\in\mathcal P_{n^i}$, in the form $(1-\epsilon)\mathbf p^i+\epsilon\mathbf q^i$. It satisfies
$$\frac{d}{d\epsilon}Z(\mathbf p^1,\ldots,\mathbf p^{i-1},(1-\epsilon)\mathbf p^i+\epsilon\mathbf q^i,\mathbf p^{i+1},\ldots,\mathbf p^m)\big|_{\epsilon=0}=\nabla^iZ(\mathbf p)'(\mathbf q^i-\mathbf p^i)$$
by the chain rule. In particular, we must have
\begin{equation}
\psi_j^i(\mathbf p)=\nabla^iZ(\mathbf p)'(\mathbf 1_j^i-\mathbf p^i)=\partial_j^iZ(\mathbf p)-\nabla^iZ(\mathbf p)'\mathbf p^i\label{interim gradient}
\end{equation}
where $\partial_j^iZ(\mathbf p)$ denotes partial derivative of $Z$ with respect to $p_j^i$. Writing \eqref{interim gradient} for all $j$ together gives
$$\bm\psi^i(\mathbf p)=\nabla^iZ(\mathbf p)-(\nabla^iZ(\mathbf p)'\mathbf p^i)\mathbf 1^i$$
where $\mathbf 1^i\in\mathbb R^{n^i}$ is a vector of 1.
Therefore
$$\bm\psi^i(\mathbf p)'(\mathbf q^i-\mathbf p^i)=(\nabla^iZ(\mathbf p)-(\nabla^iZ(\mathbf p)'\mathbf p^i)\mathbf 1^i)'(\mathbf q^i-\mathbf p^i)=\nabla^iZ(\mathbf p)'(\mathbf q^i-\mathbf p^i)$$
since $\mathbf q^i,\mathbf p^i\in\mathcal P_{n^i}$. Summing up over $i$, \eqref{gradient equivalence} follows.

To prove 2., note that we have
\begin{align}
\psi_j^i(\mathbf p)&=\frac{d}{d\epsilon}Z(\mathbf p^1,\ldots,\mathbf p^{i-1},(1-\epsilon)\mathbf p^i+\epsilon\mathbf 1_j^i,\mathbf p^{i+1},\ldots,\mathbf p^m)\big|_{\epsilon=0}\notag\\
&=\frac{d}{d\epsilon}E_{\mathbf p^1,\ldots,\mathbf p^{i-1},(1-\epsilon)\mathbf p^i+\epsilon\mathbf 1_j^i,\mathbf p^{i+1},\ldots,\mathbf p^m}[h(\mathbf X)]\Bigg|_{\epsilon=0}\notag\\
&=E_{\mathbf p}[h(\mathbf X)s_j^i(\mathbf X^i)]\label{Gateaux}
\end{align}
where $s_j^i(\cdot)$ is the score function defined as
\begin{equation}
s_j^i(\mathbf x^i)=\sum_{t=1}^{T^i}\frac{d}{d\epsilon}\log((1-\epsilon)p^i(x_t^i)+\epsilon I(x_t^i=y_j^i))\Bigg|_{\epsilon=0}.\label{score function1}
\end{equation}
Here $p^i(x_t^i)=p_j^i$ where $j$ is chosen such that $x_t^i=y_j^i$. The last equality in \eqref{Gateaux} follows from the fact that
$$\frac{d}{d\epsilon}\prod_{t=1}^{T^i}((1-\epsilon)p^i(x_t^i)+\epsilon I(x_t^i=y_j^i))\Bigg|_{\epsilon=0}=\frac{d}{d\epsilon}\sum_{t=1}^{T^i}\log((1-\epsilon)p^i(x_t^i)+\epsilon I(x_t^i=y_j^i))\Bigg|_{\epsilon=0}\cdot\prod_{t=1}^{T^i}p^i(x_t^i)$$
Note that \eqref{score function1} can be further written as
$$\sum_{t=1}^{T^i}\frac{-p^i(x_t^i)+I(x_t^i=y_j^i)}{p^i(x_t^i)}=-T^i+\sum_{t=1}^{T^i}\frac{I(x_t^i=y_j^i)}{p^i(x_t^i)}=-T^i+\sum_{t=1}^{T^i}\frac{I(x_t^i=y_j^i)}{p_j^i}$$
which leads to \eqref{score function}.
\endproof

\proof{Proof of Lemma \ref{prop:var}}
We have
\begin{equation}
Var_{\mathbf p}(h(\mathbf X)s_j^i(\mathbf X^i))\leq E_{\mathbf p}(h(\mathbf X)s_j^i(\mathbf X^i))^2\leq M^2E_{\mathbf p}(s_j^i(\mathbf X^i))^2=M^2(Var_{\mathbf p}(s_j^i(\mathbf X^i))+(E_{\mathbf p}[s_j^i(\mathbf X^i)])^2)\label{score function var}
\end{equation}
Now note that by the definition of $s_j^i(\mathbf X)$ in \eqref{score function2} we have $E_{\mathbf p}[s_j^i(\mathbf X^i)]=0$ and
$$Var_{\mathbf p}(s_j^i(\mathbf X^i))=\frac{T^iVar_{\mathbf p}(I(X_t^i=y_j^i))}{(p_j^i)^2}=\frac{T^i(1-p_j^i)}{p_j^i}$$
Hence, from \eqref{score function var}, we conclude that $Var_{\mathbf p}(h(\mathbf X)s_j^i(\mathbf X^i))\leq M^2T^i(1-p_j^i)/p_j^i$.
\endproof

\proof{Proof of Proposition \ref{phisolution}}
Consider the Lagrangian relaxation
\begin{eqnarray}
&&\max_{\alpha\geq0,\lambda\in\mathbb R}\min_{\mathbf p^i\geq\mathbf 0}\sum_{j=1}^{n^i}p_j^i\xi_j+\alpha\left(\sum_{j=1}^{n^i}p_{b,j}^i\phi\left(\frac{p_j^i}{p_{b,j}^i}\right)-\eta^i\right)+\lambda\left(\sum_{j=1}^{n^i}p_j^i-1\right)\label{Lagrangian phi}\\
&=&\max_{\alpha\geq0,\lambda\in\mathbb R}-\alpha\sum_{j=1}^{n^i}p_{b,j}^i\max_{p_j^i\geq0}\left\{-\frac{\xi_j+\lambda}{\alpha}\frac{p_j^i}{p_{b,j}^i}-\phi\left(\frac{p_j^i}{p_{b,j}^i}\right)\right\}-\alpha\eta^i-\lambda\notag\\
&=&\max_{\alpha\geq0,\lambda\in\mathbb R}-\alpha\sum_{j=1}^{n^i}p_{b,j}^i\phi^*\left(-\frac{\xi_j+\lambda}{\alpha}\right)-\alpha\eta^i-\lambda\notag
\end{eqnarray}
In the particular case that $\alpha^*=0$, the optimal value of \eqref{Lagrangian phi} is the same as
$$\max_{\lambda\in\mathbb R}\min_{\mathbf p^i\geq\mathbf 0}\sum_{j=1}^{n^i}p_j^i\xi_j+\lambda\left(\sum_{j=1}^{n^i}p_j^i-1\right)$$
whose inner minimization is equivalent to $\min_{\mathbf p^i\in\mathcal P^i}\sum_{j=1}^{n^i}p_j^i\xi_j=\min_{j\in\{1,\ldots,n^i\}}\xi_j$. Among all solutions that lead to this objective value, we find the one that solves
\begin{equation}
\min_{p_j^i,j\in\mathcal M^i:\sum_{j\in\mathcal M^i}p_j^i=1}\sum_{j\in\mathcal M^i}p_{b,j}^i\phi\left(\frac{p_j^i}{p_{b,j}^i}\right)\label{opt phi3}
\end{equation}
Now note that by the convexity of $\phi$ and Jensen's inequality, for any $\sum_{j\in\mathcal M^i}p_j^i=1$, we have
\begin{equation}
\sum_{j\in\mathcal M^i}p_{b,j}^i\phi\left(\frac{p_j^i}{p_{b,j}^i}\right)=\sum_{r\in\mathcal M^i}p_{b,r}^i\sum_{j\in\mathcal M^i}\frac{p_{b,j}^i}{\sum_{r\in\mathcal M^i}p_{b,r}^i}\phi\left(\frac{p_j^i}{p_{b,j}^i}\right)\geq\sum_{j\in\mathcal M^i}p_{b,j}^i\phi\left(\frac{1}{\sum_{j\in\mathcal M^i}p_{b,j}^i}\right)=\phi\left(\frac{1}{\sum_{j\in\mathcal M^i}p_{b,j}^i}\right)\label{opt phi4}
\end{equation}
It is easy to see that choosing $p_j^i$ in \eqref{opt phi3} as $q_j^i$ depicted in \eqref{opt phi2} achieves the lower bound in \eqref{opt phi4}, hence concluding the proposition.

\endproof

\proof{Proof of Proposition \ref{KLsolution}}
Consider the Lagrangian for the optimization \eqref{step optimization2}
\begin{equation}
\min_{\mathbf p^i\in\mathcal P^i}\sum_{j=1}^{n^i}\xi_jp_j^i+\alpha\left(\sum_{j=1}^{n^i}p_j^i\log\frac{p_j^i}{p_{b,j}^i}-\eta^i\right)\label{Lagrangian}
\end{equation}
By Theorem 1, P.220 in \cite{luenberger1969optimization}, suppose that one can find $\alpha^*\geq0$ such that $\mathbf q^i=(q_j^i)_{j=1,\ldots,n^i}\in\mathcal P_{n^i}$ minimizes \eqref{Lagrangian} for $\alpha=\alpha^*$ and moreover that $\alpha^*\left(\sum_{j=1}^{n^i}q_j^i\log\frac{q_j^i}{p_{b,j}^i}-\eta^i\right)=0$, then $\mathbf q^i$ is optimal for \eqref{step optimization2}.

Suppose $\alpha^*=0$, then the minimizer of \eqref{Lagrangian} can be any probability distributions that have masses concentrated on the set of indices in $\mathcal M^i$. Any one of these distributions that lies in $\hat{\mathcal U}^i$ will be an optimal solution to \eqref{step optimization2}. To check whether any of them lies in $\hat{\mathcal U}^i$, consider the one that has the minimum $d_\phi(\mathbf q^i,\mathbf p_b^i)$ and see whether it is less than or equal to $\eta^i$. In other words, we want to find $\min_{p_j^i,j\in\mathcal M^i:\sum_{j\in\mathcal M^i}p_j^i=1}\sum_{j\in\mathcal M^i}p_j^i\log(p_j^i/p_{b,j}^i)$. The optimal solution to this minimization is $p_{b,j}^i/\sum_{j\in\mathcal M^i}p_{b,j}^i$ for $j\in\mathcal M^i$, which gives an optimal value $-\log\sum_{j\in\mathcal M^i}p_{b,j}^i$. Thus, if $-\log\sum_{j\in\mathcal M^i}p_{b,j}^i\leq\eta^i$, we find an optimal solution $\mathbf q^i$ to \eqref{step optimization2} given by \eqref{opt2}.

In the case that $\alpha^*=0$ does not lead to an optimal solution, or equivalently $-\log\sum_{j\in\mathcal M^i}p_{b,j}^i>\eta^i$, we consider $\alpha^*>0$. We write the objective value of \eqref{Lagrangian} with $\alpha=\alpha^*$ as
\begin{equation}
\sum_{j=1}^{n^i}\xi_jp_j^i+\alpha^*\sum_{j=1}^{n^i}p_j^i\log\frac{p_j^i}{p_{b,j}^i}-\alpha^*\eta^i\label{interim new10}
\end{equation}
By Jensen's inequality,
$$\sum_{j=1}^{n^i}p_j^ie^{-\xi_j/\alpha^*-\log(p_j^i/p_{b,j}^i)}\geq e^{-\sum_{j=1}^{n^i}\xi_jp_j^i/\alpha^*-\sum_{j=1}^{n^i}p_j^i\log(p_j^i/p_{b,j}^i)}$$
giving
\begin{equation}
\sum_{j=1}^{n^i}\xi_jp_j^i+\alpha^*\sum_{j=1}^{n^i}p_j^i\log\frac{p_j^i}{p_{b,j}^i}\geq-\alpha^*\log\sum_{j=1}^{n^i}p_{b,j}^ie^{-\xi_j/\alpha^*}\label{interim new11}
\end{equation}
It is easy to verify that putting $p_j^i$ as
$$q_j^i=\frac{p_{b,j}^ie^{-\xi_j/\alpha^*}}{\sum_{r=1}^{n^i}p_{b,r}^ie^{-\xi_r/\alpha^*}}$$
gives the lower bound in \eqref{interim new11}. Thus $q_j^i$ minimizes \eqref{interim new10}. Moreover, $\alpha^*>0$ can be chosen such that
$$\sum_{j=1}^{n^i}q_j^i\log\frac{q_j^i}{p_{b,j}^i}=-\frac{\sum_{j=1}^{n^i}\xi_jp_{b,j}^ie^{-\xi_j/\alpha^*}}{\alpha^*\sum_{j=1}^{n^i}p_{b,j}^ie^{-\xi_j/\alpha^*}}-\log\sum_{j=1}^{n^i}p_{b,j}^ie^{-\xi_j/\alpha^*}=\eta^i$$
Letting $\beta=-1/\alpha^*$, we obtain \eqref{opt1} and \eqref{root}. Note that \eqref{root} must bear a negative root because of the following. Note that the left hand side of \eqref{root} is continuous, and goes to 0 when $\beta\to0$. Defining $\xi_*=\min\{\xi_j:j=1,\ldots,n^i\}$, we have, as $\beta\to-\infty$, $\varphi_{\bm\xi}^i(\beta)=\log\sum_{j=1}^{n^i}p_{b,j}^ie^{\beta\xi_j}=\log\left(\sum_{j\in\mathcal M^i}p_{b,j}^ie^{\beta\xi_*}(1+\sum_{j\notin\mathcal M^i}p_{b,j}^ie^{\beta(\xi_j-\xi_*)}/\sum_{j\in\mathcal M^i}p_{b,j}^i)\right)=\beta\xi_*+\log\sum_{j\in\mathcal M^i}p_{b,j}^i+O(e^{c_1\beta})$ for some positive constant $c_1$, and ${\varphi_{\bm\xi}^i}'(\beta)=\sum_{j=1}^{n^i}\xi_jp_{b,j}^ie^{\beta\xi_j}/\sum_{j=1}^{n^i}p_{b,j}^ie^{\beta\xi_j}=\xi_*(1+\sum_{j\notin\mathcal M^i}\xi_jp_{b,j}^ie^{\beta(\xi_j-\xi_*)}/\sum_{j\in\mathcal M^i}p_{b,j}^i)/(1+\sum_{j\notin\mathcal M^i}p_{b,j}^ie^{\beta(\xi_j-\xi_*)}/\sum_{j\in\mathcal M^i}p_{b,j}^i)=\xi_*+O(e^{c\beta})$ for some positive constant $c_2$. So $\beta{\varphi_{\bm\xi}^i}'(\beta)-\varphi_{\bm\xi}^i(\beta)=-\log\sum_{j\in\mathcal M^i}p_{b,j}^i+O(e^{(c_1\wedge c_2)\beta})>\eta^i$ when $\beta$ is negative enough.
\endproof

\proof{Proof of Theorem \ref{as}}
The proof is an adaptation of \cite{blum1954multidimensional}. Recall that $\mathbf p_k=\text{vec}(\mathbf p_k^i:i=1,\ldots,m)$ where we write each component of $\mathbf p_k$ as $p_{k,j}^i$. Let $N=\sum_{i=1}^mn^i$ be the total counts of support points. Since $h(\mathbf X)$ is bounded a.s., we have $|h(\mathbf X)|\leq M$ a.s. for some $M$. Without loss of generality, we assume that $Z(\mathbf p)\geq0$ for all $\mathbf p$. Also note that $Z(\mathbf p)$, as a high-dimensional polynomial, is continuous everywhere in $\hat{\mathcal U}$.

For notational convenience, we write $\mathbf d_k=\mathbf q(\mathbf p_k)-\mathbf p_k$ and $\hat{\mathbf d}_k=\hat{\mathbf q}(\mathbf p_k)-\mathbf p_k$, i.e. $\mathbf d_k$ is the $k$-th step best feasible direction given the exact gradient estimate, and $\hat{\mathbf d}_k$ is the one with estimated gradient.

Now, given $\mathbf p_k$, consider the iterative update $\mathbf p_{k+1}=(1-\epsilon_k)\mathbf p_k+\epsilon_k\hat{\mathbf q}(\mathbf p_k)=\mathbf p_k+\epsilon_k\hat{\mathbf d}_k$. We have, by Taylor series expansion,
$$Z(\mathbf p_{k+1})=Z(\mathbf p_k)+\epsilon_k\nabla Z(\mathbf p_k)'\hat{\mathbf d}_k+\frac{\epsilon_k^2}{2}\hat{\mathbf d}_k'\nabla^2Z(\mathbf p_k+\theta_k\epsilon_k\hat{\mathbf d}_k)\hat{\mathbf d}_k$$
for some $\theta_k$ between 0 and 1. By Theorem \ref{prop:gradient}, we can rewrite the above as
\begin{equation}
Z(\mathbf p_{k+1})=Z(\mathbf p_k)+\epsilon_k\bm\psi(\mathbf p_k)'\hat{\mathbf d}_k+\frac{\epsilon_k^2}{2}\hat{\mathbf d}_k'\nabla^2Z(\mathbf p_k+\theta_k\epsilon_k\hat{\mathbf d}_k)\hat{\mathbf d}_k\label{interim as1}
\end{equation}
Consider the second term in the right hand side of \eqref{interim as1}. We can write
\begin{eqnarray}
\bm\psi(\mathbf p_k)'\hat{\mathbf d}_k&=&\hat{\bm\psi}(\mathbf p_k)'\hat{\mathbf d}_k+(\bm\psi(\mathbf p_k)-\hat{\bm\psi}(\mathbf p_k))'\hat{\mathbf d}_k\notag\\
&\leq&\hat{\bm\psi}(\mathbf p_k)'\mathbf d_k+(\bm\psi(\mathbf p_k)-\hat{\bm\psi}(\mathbf p_k))'\hat{\mathbf d}_k\text{\ \ \ \ by the definition of $\hat{\mathbf d}_k$}\notag\\
&=&\bm\psi(\mathbf p_k)'\mathbf d_k+(\hat{\bm\psi}(\mathbf p_k)-\bm\psi(\mathbf p_k))'\mathbf d_k+(\bm\psi(\mathbf p_k)-\hat{\bm\psi}(\mathbf p_k))'\hat{\mathbf d}_k\notag\\
&=&\bm\psi(\mathbf p_k)'\mathbf d_k+(\hat{\bm\psi}(\mathbf p_k)-\bm\psi(\mathbf p_k))'(\mathbf d_k-\hat{\mathbf d}_k)\label{interim as11}
\end{eqnarray}
Hence \eqref{interim as1} and \eqref{interim as11} together imply
$$Z(\mathbf p_{k+1})\leq Z(\mathbf p_k)+\epsilon_k\bm\psi(\mathbf p_k)'\mathbf d_k+\epsilon_k(\hat{\bm\psi}(\mathbf p_k)-\bm\psi(\mathbf p_k))'(\mathbf d_k-\hat{\mathbf d}_k)+\frac{\epsilon_k^2}{2}\hat{\mathbf d}_k'\nabla^2Z(\mathbf p_k+\theta_k\epsilon_k\hat{\mathbf d}_k)\hat{\mathbf d}_k$$
%
Let $\mathcal F_k$ be the filtration generated by $\mathbf p_1,\ldots,\mathbf p_k$. We then have
\begin{eqnarray}
E[Z(\mathbf p_{k+1})|\mathcal F_k]&\leq& Z(\mathbf p_k)+\epsilon_k\bm\psi(\mathbf p_k)'\mathbf d_k+\epsilon_kE[(\hat{\bm\psi}(\mathbf p_k)-\bm\psi(\mathbf p_k))'(\mathbf d_k-\hat{\mathbf d}_k)|\mathcal F_k]{}\notag\\
&&{}+\frac{\epsilon_k^2}{2}E[\hat{\mathbf d}_k'\nabla^2Z(\mathbf p_k+\theta_k\epsilon_k\hat{\mathbf d}_k)\hat{\mathbf d}_k|\mathcal F_k]
%
\label{interim as}
\end{eqnarray}
We analyze \eqref{interim as} term by term. First, since $Z(\mathbf p)$ is a high-dimensional polynomial and $\hat{\mathcal U}$ is a bounded set, the largest eigenvalue of the Hessian matrix $\nabla^2Z(\mathbf p)$, for any $\mathbf p\in\hat{\mathcal U}$, is uniformly bounded by a constant $H>0$. Hence
\begin{equation}
E[\hat{\mathbf d}_k'\nabla^2Z(\mathbf p_k+\theta_k\epsilon_k\hat{\mathbf d}_k)\hat{\mathbf d}_k|\mathcal F_k]\leq HE[\|\hat{\mathbf d}_k\|^2|\mathcal F_k]\leq V<\infty\label{interim as3}
\end{equation}
for some $V>0$. Now
\begin{eqnarray}
&&E[(\hat{\bm\psi}(\mathbf p_k)-\bm\psi(\mathbf p_k))'(\mathbf d_k-\hat{\mathbf d}_k)|\mathcal F_k]\\
&\leq&\sqrt{E[\|\hat{\bm\psi}(\mathbf p_k)-\bm\psi(\mathbf p_k)\|^2|\mathcal F_k]E[\|\mathbf d_k-\hat{\mathbf d}_k\|^2|\mathcal F_k]}\text{\ \ \ \ by Cauchy-Schwarz inequality}\notag\\
&\leq&\sqrt{E[\|\hat{\bm\psi}(\mathbf p_k)-\bm\psi(\mathbf p_k)\|^2|\mathcal F_k]E[2(\|\mathbf d_k\|^2+\|\hat{\mathbf d}_k\|^2)|\mathcal F_k]}\text{\ \ \ \ by parallelogram law}\notag\\
&\leq&\sqrt{8mE[\|\hat{\bm\psi}(\mathbf p_k)-\bm\psi(\mathbf p_k)\|^2|\mathcal F_k]}\text{\ \ \ \ since $\|\mathbf d_k\|^2,\|\hat{\mathbf d}_k\|^2\leq2m$ by using the fact that $\mathbf p_k,\mathbf q(\mathbf p_k),\hat{\mathbf q}(\mathbf p_k)\in\mathcal P$}\notag\\
&\leq&\sqrt{\frac{8mM^2T}{R_k}\sum_{i,j}\frac{1-p_{k,j}^i}{p_{k,j}^i}}\text{\ \ \ \ by Lemma \ref{prop:var}}\notag\\
&\leq&M\sqrt{\frac{8mTN}{R_k\min_{i,j}p_{k,j}^i}}\label{interim as2}
\end{eqnarray}
Note that by iterating the update rule $(1-\epsilon_k)\mathbf p_k+\epsilon_k\mathbf q_k$, we have
$$\min_{i,j}p_{k,j}^i\geq\prod_{j=1}^{k-1}(1-\epsilon_j)\delta$$
where $\delta=\min_{i,j}p_{1,j}^i>0$. We thus have \eqref{interim as2} less than or equal to
\begin{equation}
M\sqrt{\frac{8mTN}{\delta R_k}}\prod_{j=1}^{k-1}(1-\epsilon_j)^{-1/2}\label{interim as4}
\end{equation}
%

Therefore, noting that $\bm\psi(\mathbf p_k)'\mathbf d_k\leq0$ by the definition of $\mathbf d_k$, from \eqref{interim as} we have
\begin{align}
E[Z(\mathbf p_{k+1})-Z(\mathbf p_k)|\mathcal F_k]\leq\epsilon_kM\sqrt{\frac{8mTN}{\delta R_k}}\prod_{j=1}^{k-1}(1-\epsilon_j)^{-1/2}+\frac{\epsilon_k^2V}{2}\label{interim as5}
\end{align}
and hence
$$\sum_{k=1}^\infty E[E[Z(\mathbf p_{k+1})-Z(\mathbf p_k)|\mathcal F_k]^+]\leq M\sqrt{\frac{8mTN}{\delta}}\sum_{k=1}^\infty\frac{\epsilon_k}{\sqrt{R_k}}\prod_{j=1}^{k-1}(1-\epsilon_j)^{-1/2}+\sum_{k=1}^\infty\frac{\epsilon_k^2V}{2}$$
By Assumptions \ref{tuning} and \ref{sample size tuning}, and Lemma \ref{prelim} (depicted after this proof), we have $Z(\mathbf p_k)$ converge to an integrable random variable.
%

Now take expectation on \eqref{interim as} further to get
\begin{eqnarray*}
E[Z(\mathbf p_{k+1})]&\leq&E[Z(\mathbf p_k)]+\epsilon_kE[\bm\psi(\mathbf p_k)'\mathbf d_k]+\epsilon_kE[(\hat{\bm\psi}(\mathbf p_k)-\bm\psi(\mathbf p_k))'(\mathbf d_k-\hat{\mathbf d}_k)]{}\\
&&{}+\frac{\epsilon_k^2}{2}E[\hat{\mathbf d}_k'\nabla^2Z(\mathbf p_k+\theta_k\epsilon_k\hat{\mathbf d}_k)\hat{\mathbf d}_k]
\end{eqnarray*}
and telescope to get
\begin{eqnarray}
E[Z(\mathbf p_{k+1})]&\leq&E[Z(\mathbf p_1)]+\sum_{j=1}^k\epsilon_jE[\bm\psi(\mathbf p_j)'\mathbf d_j]+\sum_{j=1}^k\epsilon_jE[(\hat{\bm\psi}(\mathbf p_j)-\bm\psi(\mathbf p_j))'(\mathbf d_j-\hat{\mathbf d}_j)]{}\notag\\
&&{}+\sum_{j=1}^k\frac{\epsilon_j^2}{2}E[\hat{\mathbf d}_j'\nabla^2Z(\mathbf p_j+\theta_j\epsilon_j\hat{\mathbf d}_j)\hat{\mathbf d}_j]\label{interim as6}
\end{eqnarray}
Now take the limit on both sides of \eqref{interim as6}. Note that $E[Z(\mathbf p_{k+1})]\to E[Z_\infty]$ for some integrable $Z_\infty$ by dominated convergence theorem. Also $Z(\mathbf p_1)<\infty$, and by \eqref{interim as3} and \eqref{interim as4} respectively, we have
$$\lim_{k\to\infty}\sum_{j=1}^k\frac{\epsilon_j^2}{2}E[\hat{\mathbf d}_j'\nabla Z(\mathbf p_j+\theta_j\epsilon_j\hat{\mathbf d}_j)\hat{\mathbf d}_j]\leq\sum_{j=1}^\infty\frac{\epsilon_j^2V}{2}<\infty$$
and
$$\lim_{k\to\infty}\sum_{j=1}^k\epsilon_jE[(\hat{\bm\psi}(\mathbf p_j)-\bm\psi(\mathbf p_j))'(\mathbf d_j-\hat{\mathbf d}_j)]\leq M\sqrt{\frac{8mTN}{\delta}}\sum_{j=1}^\infty\frac{\epsilon_j}{\sqrt{R_j}}\prod_{i=1}^{j-1}(1-\epsilon_i)^{-1/2}<\infty$$
Therefore, from \eqref{interim as6}, and since $E[\bm\psi(\mathbf p_j)'\mathbf d_j]\leq0$, we must have $\sum_{j=1}^k\epsilon_jE[\bm\psi(\mathbf p_j)'\mathbf d_j]$ converges a.s., which implies that $\limsup_{k\to\infty}E[\bm\psi(\mathbf p_k)'\mathbf d_k]=0$. So there exists a subsequence $k_i$ such that $\lim_{i\to\infty}E[\bm\psi(\mathbf p_{k_i})'\mathbf d_{k_i}]=0$. 
This in turn implies that $\bm\psi(\mathbf p_{k_i})'\mathbf d_{k_i}\stackrel{p}{\to}0$. Then, there exists a further subsequence $l_i$ such that $\bm\psi(\mathbf p_{l_i})'\mathbf d_{l_i}\to0$ a.s..

Consider part \ref{as part1} of the theorem. Let $S^*=\{\mathbf p\in\mathcal P:g(\mathbf p)=0\}$. Since $g(\cdot)$ is continuous, we have $D(\mathbf p_{l_i},S^*)\to0$ a.s.. Since $Z(\cdot)$ is continuous, we have $D(Z(\mathbf p_{l_i}),\mathcal Z^*)\to0$ a.s.. But since we have proven that $Z(\mathbf p_k)$ converges a.s., we have $D(Z(\mathbf p_k),\mathcal Z^*)\to0$ a.s.. This gives part \ref{as part1} of the theorem.

Now consider part \ref{as part2}. By Assumption \ref{main assumption}, since $\mathbf p^*$ is the only $\mathbf p$ such that $g(\mathbf p)=0$ and $g(\cdot)$ is continuous, we must have $\mathbf p_{l_i}\to\mathbf p^*$ a.s.. Since $Z(\cdot)$ is continuous, we have $Z(\mathbf p_{l_i})\to Z(\mathbf p^*)$. But since $Z(\mathbf p_k)$ converges a.s. as shown above, we must have $Z(\mathbf p_k)\to Z(\mathbf p^*)$. Then by Assumption \ref{main assumption} again, since $\mathbf p^*$ is the unique optimizer, we have $\mathbf p_k\to\mathbf p^*$ a.s.. This concludes part \ref{as part2} of the theorem.
\endproof

\begin{lemma}[Adapted from \cite{blum1954multidimensional}]
Consider a sequence of integrable random variable $Y_k,k=1,2,\ldots$. Let $\mathcal F_k$ be the filtration generated by $Y_1,\ldots,Y_k$. Assume
$$\sum_{k=1}^\infty E[E[Y_{k+1}-Y_k|\mathcal F_k]^+]<\infty$$
where $x^+$ denotes the positive part of $x$, i.e. $x^+=x$ if $x\geq0$ and $0$ if $x<0$. Moreover, assume that $Y_k$ is bounded uniformly from above. Then $Y_k\to Y_\infty$ a.s., where $Y_\infty$ is an integrable random variable.\label{prelim}
\end{lemma}
The lemma follows from \cite{blum1954multidimensional}, with the additional conclusion that $Y_\infty$ is integrable, which is a direct consequence of the martingale convergence theorem.

\begin{theorem}[Conditions in Theorem \ref{rate thm}]
Conditions \ref{c1}-\ref{c8} needed in Theorem \ref{rate thm} are:
\begin{enumerate}
\item $$k_0\geq2a\left(\frac{4KMTm}{c^2\tau^2}+\frac{KL\vartheta}{c\tau}\right)$$\label{c1}
\item $$-\left(1-\frac{2KL\vartheta}{c\tau}-\frac{2a\varrho K}{c^2\tau^2k_0}\right)\nu+\frac{2aKL\vartheta\varrho}{c\tau k_0^{1+\gamma}}+\frac{\varrho}{k_0^\gamma}+\frac{2K\nu^2}{c^2\tau^2}\leq0$$\label{c2}
\item $$\frac{2KL\vartheta}{c\tau}+\frac{2K\nu}{c^2\tau^2}<1$$\label{c3}
\item $$\frac{a}{k_0}\left(1-\frac{2KL\vartheta}{c\tau}-\frac{2K\nu}{c^2\tau^2}\right)<1$$\label{c extra}
\item $$k_0\geq\frac{a\rho}{\rho-1}$$\label{c4}
\item $$\beta>\rho a+2\gamma+2$$\label{c5}
\item
\begin{eqnarray*}
&&\prod_{j=1}^{k_0-1}(1-\epsilon_j)^{-1}\frac{M^2TN}{\vartheta^2\delta b}\frac{1}{(\beta-\rho a-1)(k_0-1)^{\beta-1}}{}\\
&&+\prod_{j=1}^{k_0-1}(1-\epsilon_j)^{-1/2}\frac{M}{\varrho}\sqrt{\frac{8mTN}{\delta b}}\frac{1}{((\beta-\rho a)/2-\gamma-1)(k_0-1)^{\beta/2-\gamma-1}}<\varepsilon
\end{eqnarray*}
where $N=\sum_{i=1}^mn^i$ is the total count of all support points.\label{c6}
\item $K>0$ is a constant such that $|\mathbf x'\nabla^2Z(\mathbf p)\mathbf y|\leq K\|\mathbf x\|\|\mathbf y\|$ for any $x,y\in\mathbb R^n$ and $\mathbf p\in\mathcal A$ (which must exist because $Z(\cdot)$ is a polynomial defined over a bounded set).\label{c7}
\item $\delta=\min_{\substack{i=1,\ldots,m\\j=1,\ldots,n^i}}p_{1,j}^i>0$\label{c8}
\end{enumerate}
\end{theorem}

\proof{Proof of Theorem \ref{rate thm}}
We adopt the notation as in the proof of Theorem \ref{as}. In addition, for convenience, we write $\bm\psi_k=\bm\psi(\mathbf p_k)$, $\hat{\bm\psi}_k=\hat{\bm\psi}(\mathbf p_k)$, $\mathbf q_k=\mathbf q(\mathbf p_k)$, $\hat{\mathbf q}_k=\hat{\mathbf q}(\mathbf p_k)$, $g_k=g(\mathbf p_k)=-\bm\psi(\mathbf p_k)'\mathbf d_k$, $\nabla Z_k=\nabla Z(\mathbf p_k)$, and $\nabla^2Z_k=\nabla^2Z(\mathbf p_k)$. Note that $\mathbf p_{k+1}=\mathbf p_k+\epsilon_k\hat{\mathbf d}_k$.

First, by the proof of Theorem \ref{as}, given any $\nu$ and $\tilde k_0$, almost surely there must exists a $k_0\geq\tilde k_0$ such that $g_{k_0}\leq\nu$. If the optimal solution is reached and is kept there, then $g_k=0$ from thereon and the algorithm reaches and remains at optimum at finite time, hence there is nothing to prove. So let us assume that $0<g_{k_0}\leq\nu$. Moreover, let us assume that $\nu$ is chosen small enough so that for any $\mathbf p$ with $g(\mathbf p)\leq\nu$ and $\mathbf p>\mathbf 0$, we have $\bm\psi(\mathbf p)\in\mathcal N_{\Delta-\vartheta}(\bm\psi(\mathbf p^*))$ (which can be done since $g(\cdot)$ is assumed continuous by Assumption \ref{main assumption} and $\bm\psi(\mathbf p)$ is continuous for any $\mathbf p>\mathbf 0$ by the construction in Theorem \ref{prop:gradient}).

We consider the event
$$\mathcal E=\bigcup_{k=k_0}^\infty\mathcal E_k\cup\bigcup_{k=k_0}^\infty\mathcal E_k'$$
where
$$\mathcal E_k=\{\|\hat{\bm\psi}_k-\bm\psi_k\|>\vartheta\}$$
and
$$\mathcal E_k'=\left\{|(\hat{\bm\psi}_k-\bm\psi_k)'(\hat{\mathbf d}_k-\mathbf d_k)|>\frac{\varrho}{k^\gamma}\right\}$$
Note that by the Markov inequality,
$$P(\mathcal E_k)\leq\frac{E\|\hat{\bm\psi}_k-\bm\psi_k\|^2}{\vartheta^2}\leq\frac{M^2T}{\vartheta^2R_k}\sum_{i,j}\frac{1-p_{k,j}^i}{p_{k,j}^i}\leq\frac{M^2TN}{\vartheta^2R_k\delta}\prod_{j=1}^{k-1}(1-\epsilon_j)^{-1}$$
where the second inequality follows from Lemma \ref{prop:var} and the last inequality follows as in the derivation in \eqref{interim as2} and \eqref{interim as4}. On the other hand, we have
\begin{equation}
P(\mathcal E_k')\leq\frac{k^\gamma E|(\hat{\bm\psi}_k-\bm\psi_k)'(\hat{\mathbf d}_k-\mathbf d_k)|}{\varrho}\leq\frac{k^\gamma M}{\varrho}\sqrt{\frac{8mTN}{\delta R_k}}\prod_{j=1}^{k-1}(1-\epsilon_j)^{-1/2}\label{interim remark}
\end{equation}
by following the derivation in \eqref{interim as2} and \eqref{interim as4}.
%
Therefore,
\begin{align}
P(\mathcal E)&\leq\sum_{k=k_0}^\infty P(\mathcal E_k)+\sum_{k=k_0}^\infty P(\mathcal E_k')\notag\\
&\leq\frac{M^2TN}{\vartheta^2\delta}\sum_{k=k_0}^\infty\frac{1}{R_k}\prod_{j=1}^{k-1}(1-\epsilon_j)^{-1}+\frac{M}{\varrho}\sqrt{\frac{8mTN}{\delta}}\sum_{k=k_0}^\infty\frac{k^\gamma}{\sqrt{R_k}}\prod_{j=1}^{k-1}(1-\epsilon_j)^{-1/2}\notag\\
&=\prod_{j=1}^{k_0-1}(1-\epsilon_j)^{-1}\frac{M^2TN}{\vartheta^2\delta}\sum_{k=k_0}^\infty\frac{1}{R_k}\prod_{j=k_0}^{k-1}(1-\epsilon_j)^{-1}+\prod_{j=1}^{k_0-1}(1-\epsilon_j)^{-1/2}\frac{M}{\varrho}\sqrt{\frac{8mTN}{\delta}}\sum_{k=k_0}^\infty\frac{k^\gamma}{\sqrt{R_k}}\prod_{j=k_0}^{k-1}(1-\epsilon_j)^{-1/2}\label{interim rate1}
\end{align}
Now recall that $\epsilon_k=a/k$. Using the fact that $1-x\geq e^{-\rho x}$ for any $0\leq x\leq(\rho-1)/\rho$ and $\rho>1$, we have, for any
$$\frac{a}{k}\leq\frac{\rho-1}{\rho}$$
or equivalently
$$k\geq\frac{a\rho}{\rho-1}$$
we have
$$1-\epsilon_k=1-\frac{a}{k}\geq e^{-\rho a/k}$$
Hence choosing $k_0$ satisfying Condition \ref{c4}, we get
\begin{equation}
\prod_{j=k_0}^{k-1}(1-\epsilon_j)^{-1}\leq e^{\rho a\sum_{k_0}^{k-1}1/j}\leq\left(\frac{k-1}{k_0-1}\right)^{\rho a}\label{interim rate2}
\end{equation}
Therefore, picking $R_k=bk^\beta$ and using \eqref{interim rate2}, we have \eqref{interim rate1} bounded from above by
\begin{eqnarray}
&&\prod_{j=1}^{k_0-1}(1-\epsilon_j)^{-1}\frac{M^2TN}{\vartheta^2\delta b}\sum_{k=k_0}^\infty\frac{1}{(k_0-1)^{\rho a}k^{\beta-\rho a}}+\prod_{j=1}^{k_0-1}(1-\epsilon_j)^{-1/2}\frac{M}{\varrho}\sqrt{\frac{8mTN}{\delta b}}\sum_{k=k_0}^\infty\frac{1}{(k_0-1)^{\rho a/2}k^{(\beta-\rho a)/2-\gamma}}\notag\\
&\leq&\prod_{j=1}^{k_0-1}(1-\epsilon_j)^{-1}\frac{M^2TN}{\vartheta^2\delta b}\frac{1}{(\beta-\rho a-1)(k_0-1)^{\beta-1}}{}\notag\\
&&{}+\prod_{j=1}^{k_0-1}(1-\epsilon_j)^{-1/2}\frac{M}{\varrho}\sqrt{\frac{8mTN}{\delta b}}\frac{1}{((\beta-\rho a)/2-\gamma-1)(k_0-1)^{\beta/2-\gamma-1}}\label{interim remark1}
\end{eqnarray}
if Condition \ref{c5} holds. Then Condition \ref{c6} guarantees that $P(\mathcal E)<\varepsilon$.

The rest of the proof will show that under the event $\mathcal E^c$, we must have the bound \eqref{main}, hence concluding the theorem. To this end, we first set up a recursive representation of $g_k$. Consider
\begin{align}
g_{k+1}&=-\bm\psi_{k+1}'\mathbf d_{k+1}=-\bm\psi_{k+1}'(\mathbf q_{k+1}-\mathbf p_{k+1})\notag\\
&=-\bm\psi_k'(\mathbf q_{k+1}-\mathbf p_{k+1})+(\bm\psi_k-\bm\psi_{k+1})'(\mathbf q_{k+1}-\mathbf p_{k+1})\notag\\
&=-\bm\psi_k'(\mathbf q_{k+1}-\mathbf p_k)+\bm\psi_k'(\mathbf p_{k+1}-\mathbf p_k)+(\bm\psi_k-\bm\psi_{k+1})'(\mathbf q_{k+1}-\mathbf p_{k+1})\notag\\
&\leq g_k+\epsilon_k\bm\psi_k'\hat{\mathbf d}_k+(\bm\psi_k-\bm\psi_{k+1})'\mathbf d_{k+1}\text{\ \ \ \ by the definition of $g_k$, $\hat{\mathbf d}_k$ and $\mathbf d_{k+1}$}\notag\\
&\leq g_k-\epsilon_kg_k+\epsilon_k(\hat{\bm\psi}_k-\bm\psi_k)'(\mathbf d_k-\hat{\mathbf d}_k)+(\bm\psi_k-\bm\psi_{k+1})'\mathbf d_{k+1}\text{\ \ \ \ by \eqref{interim as11}}\notag\\
&=(1-\epsilon_k)g_k+(\nabla Z_k-\nabla Z_{k+1})'\mathbf d_{k+1}+\epsilon_k(\hat{\bm\psi}_k-\bm\psi_k)'(\mathbf d_k-\hat{\mathbf d}_k)\label{interim rate}
\end{align}
Now since $\nabla Z(\cdot)$ is continuously differentiable, we have $\nabla Z_{k+1}=\nabla Z_k+\epsilon_k\nabla^2Z(\mathbf p_k+\tilde{\theta}_k\hat{\mathbf d}_k)\hat{\mathbf d}_k$ for some $\tilde{\theta}_k$ between 0 and 1. Therefore \eqref{interim rate} is equal to
\begin{eqnarray}
&&(1-\epsilon_k)g_k-\epsilon_k\hat{\mathbf d}_k'\nabla^2Z(\mathbf p_k+\tilde{\theta}_k\hat{\mathbf d}_k)\mathbf d_{k+1}+\epsilon_k(\hat{\bm\psi}_k-\bm\psi_k)'(\mathbf d_k-\hat{\mathbf d}_k)\notag\\
&\leq&(1-\epsilon_k)g_k+\epsilon_kK\|\hat{\mathbf d}_k\|\|\mathbf d_{k+1}\|+\epsilon_k(\hat{\bm\psi}_k-\bm\psi_k)'(\mathbf d_k-\hat{\mathbf d}_k)\text{\ \ \ \ by Condition \ref{c7}}\notag\\
&\leq&(1-\epsilon_k)g_k+\epsilon_kK\|\mathbf d_k\|\|\mathbf d_{k+1}\|+\epsilon_kK\|\hat{\mathbf d}_k-\mathbf d_k\|\|\mathbf d_{k+1}\|+\epsilon_k(\hat{\bm\psi}_k-\bm\psi_k)'(\mathbf d_k-\hat{\mathbf d}_k){}\notag\\
&&{}\text{\ \ \ \ by the triangle inequality}\notag\\
&\leq&(1-\epsilon_k)g_k+\epsilon_kK\frac{g_kg_{k+1}}{c^2\|\bm\psi_k\|\|\bm\psi_{k+1}\|}+\epsilon_kKL\|\hat{\bm\psi}_k-\bm\psi_k\|\frac{g_{k+1}}{c\|\bm\psi_{k+1}\|}+\epsilon_k(\hat{\bm\psi}_k-\bm\psi_k)'(\mathbf d_k-\hat{\mathbf d}_k){}\notag\\
&&{}\text{\ \ \ \ by using Assumption \ref{bias} with the fact that $g_k\leq\nu$  and hence $\bm\psi_k,\hat{\bm\psi}_k\in\mathcal N_\Delta(\bm\psi(\mathbf p^*))$, and also}{}\notag\\
&&{}\text{\ \ \ \ Assumption \ref{angle}. The fact $g_k\leq\nu$ will be proved later by induction.}\notag\\
&\leq&(1-\epsilon_k)g_k+\epsilon_k\frac{K}{c^2\tau^2}g_kg_{k+1}+\epsilon_k\frac{KL}{c\tau}\|\hat{\bm\psi}_k-\bm\psi_k\|g_{k+1}+\epsilon_k(\hat{\bm\psi}_k-\bm\psi_k)'(\mathbf d_k-\hat{\mathbf d}_k){}\label{interim rate3}\\
&&{}\text{\ \ \ \ by Assumption \ref{nonzero gradient}}\notag
\end{eqnarray}
Now under the event $\mathcal E^c$, and noting that $\epsilon=a/k$, \eqref{interim rate3} implies that
$$g_{k+1}\leq\left(1-\frac{a}{k}\right)g_k+\frac{aK}{c^2\tau^2k}g_kg_{k+1}+\frac{aKL\vartheta}{c\tau k}g_{k+1}+\frac{a\varrho}{k^{1+\gamma}}$$
or
$$\left(1-\frac{aK}{c^2\tau^2k}g_k-\frac{aKL\vartheta}{c\tau k}\right)g_{k+1}\leq\left(1-\frac{a}{k}\right)g_k+\frac{a\varrho}{k^{1+\gamma}}$$
We claim that $|g_k|=|\bm\psi_k'\mathbf d_k|\leq4MTm$, which can be seen by writing
\begin{align}
\psi_j^i(\mathbf p)&=E_{\mathbf p}[h(\mathbf X)s_j^i(\mathbf X^i)]=\sum_{t=1}^{T^i}E_{\mathbf p}\left[h(\mathbf X)\frac{I(X_t^i=y_j^i)}{p_j^i}\right]-T^iE_{\mathbf p}[h(\mathbf X)]\notag\\
&=\sum_{t=1}^{T^i}E_{\mathbf p}[h(\mathbf X)|X_t=y_j^i]-T^iE_{\mathbf p}[h(\mathbf X)]
\end{align}
so that $|\psi_j^i(\mathbf p)|\leq2MT^i$ for any $\mathbf p$ and $i$. Using this and the fact that $1/(1-x)\leq1+2x$ for any $0\leq x\leq1/2$, we have, for
\begin{equation}
\frac{4aKMTm}{c^2\tau^2k}+\frac{aKL\vartheta}{c\tau k}\leq\frac{1}{2}\label{interim rate4}
\end{equation}
we must have
\begin{equation}
g_{k+1}\leq\left(1+\frac{2aK}{c^2\tau^2k}g_k+\frac{2aKL\vartheta}{c\tau k}\right)\left(\left(1-\frac{a}{k}\right)g_k+\frac{a\varrho}{k^{1+\gamma}}\right)\label{interim rate5}
\end{equation}
Note that \eqref{interim rate4} holds if
$$k\geq2a\left(\frac{4KMTm}{c^2\tau^2}+\frac{KL\vartheta}{c\tau}\right)$$
which is Condition \ref{c1} in the theorem. Now \eqref{interim rate5} can be written as
\begin{align}
g_{k+1}&\leq\left(1-\frac{a}{k}+\frac{2aKL\vartheta}{c\tau k}+\frac{2a^2K\varrho}{c^2\tau^2k^{2+\gamma}}\right)g_k+\frac{a\varrho}{k^{1+\gamma}}+\frac{2a^2KL\vartheta\varrho}{c\tau k^{2+\gamma}}-\frac{2a^2KL\vartheta}{c\tau k^2}g_k+\frac{2aK}{c^2\tau^2k}\left(1-\frac{a}{k}\right)g_k^2\notag\\
&\leq\left(1-\frac{a}{k}+\frac{2aKL\vartheta}{c\tau k}+\frac{2a^2K\varrho}{c^2\tau^2k^{2+\gamma}}\right)g_k+\frac{a\varrho}{k^{1+\gamma}}+\frac{2a^2KL\vartheta\varrho}{c\tau k^{2+\gamma}}+\frac{2aK}{c^2\tau^2k}\left(1-\frac{a}{k}\right)g_k^2\label{interim rate6}
\end{align}
We argue that under Condition \ref{c2}, we must have $g_k\leq\nu$ for all $k\geq k_0$. This can be seen by induction using \eqref{interim rate6}. By our setting at the beginning of this proof we have $g_{k_0}\leq\nu$. Suppose $g_k\leq\nu$ for some $k$. We then have
\begin{align}
g_{k+1}&\leq\left(1-\frac{a}{k}+\frac{2aKL\vartheta}{c\tau k}+\frac{2a^2K\varrho}{c^2\tau^2k^{2+\gamma}}\right)\nu+\frac{a\varrho}{k^{1+\gamma}}+\frac{2a^2KL\vartheta\varrho}{c\tau k^{2+\gamma}}+\frac{2aK}{c^2\tau^2k}\left(1-\frac{a}{k}\right)\nu^2\notag\\
&\leq\nu+\frac{a}{k}\left(\left(-1+\frac{2KL\vartheta}{c\tau}+\frac{2aK\varrho}{c^2\tau^2k^{1+\gamma}}\right)\nu+\frac{\varrho}{k_0^\gamma}+\frac{2aKL\vartheta\varrho}{c\tau k_0^{1+\gamma}}+\frac{2K\nu^2}{c^2\tau^2}\right)\notag\\
&\leq\nu\label{interim rate8}
\end{align}
by Condition \ref{c2}. This concludes our claim.

Given that $g_k\leq\nu$ for all $k\geq k_0$, \eqref{interim rate5} implies that
\begin{align}
g_{k+1}&\leq\left(1-\frac{a}{k}\left(1-\frac{2KL\vartheta}{c\tau}\right)-\frac{2a^2KL\vartheta}{c\tau k^2}+\frac{2aK\nu}{c^2\tau^2k}\left(1-\frac{a}{k}\right)\right)g_k+\frac{a\varrho}{k^{1+\gamma}}+\frac{a^2\varrho}{k^{2+\gamma}}\left(\frac{2K\nu}{c^2\tau^2}+\frac{2KL\vartheta}{c\tau}\right)\notag\\
&\leq\left(1-\frac{a}{k}\left(1-\frac{2KL\vartheta}{c\tau}-\frac{2K\nu}{c^2\tau^2}\right)\right)g_k+\frac{a\varrho}{k^{1+\gamma}}+\frac{a^2\varrho}{k^{2+\gamma}}\left(\frac{2K\nu}{c^2\tau^2}+\frac{2KL\vartheta}{c\tau}\right)\notag\\
&\leq\left(1-\frac{C}{k}\right)g_k+\frac{G}{k^{1+\gamma}}\label{interim rate7}
\end{align}
where
$$C=a\left(1-\frac{2KL\vartheta}{c\tau}-\frac{2K\nu}{c^2\tau^2}\right)$$
and
$$G=a\varrho+\frac{a^2\varrho}{k_0}\left(\frac{2K\nu}{c^2\tau^2}+\frac{2KL\vartheta}{c\tau}\right)$$
Now note that Conditions \ref{c3} and \ref{c extra} imply $C>0$ and $1-C/k>0$ respectively. By recursing the relation \eqref{interim rate7}, we get
\begin{align*}
g_{k+1}&\leq\prod_{j=k_0}^k\left(1-\frac{C}{j}\right)g_{k_0}+\sum_{j=k_0}^k\prod_{i=j+1}^k\left(1-\frac{C}{i}\right)\frac{G}{j^{1+\gamma}}\\
&\leq e^{-C\sum_{j=k_0}^k1/j}g_{k_0}+\sum_{j=k_0}^ke^{-C\sum_{i=j+1}^k1/i}\frac{G}{j^{1+\gamma}}\\
&\leq\left(\frac{k_0}{k+1}\right)^Cg_{k_0}+\sum_{j=k_0}^k\left(\frac{j+1}{k+1}\right)^C\frac{G}{j^{1+\gamma}}\\
&\leq\left(\frac{k_0}{k+1}\right)^Cg_{k_0}+\left(1+\frac{1}{k_0}\right)^CG\times\left\{\begin{array}{ll}
\frac{1}{(C-\gamma)(k+1)^\gamma}&\text{\ if\ }0<\gamma<C\\
\frac{1}{(\gamma-C)(k_0-1)^{\gamma-C}(k+1)^C}&\text{\ if\ }\gamma>C\\
\frac{\log(k/(k_0-1))}{(k+1)^C}&\text{\ if\ }\gamma=C
\end{array}\right.
\end{align*}
which gives \eqref{main}. This concludes the proof.


\endproof

\proof{Proof of Corollary \ref{rate cor}}
We use the notations in the proof of Theorem \ref{rate thm}. Our analysis starts from \eqref{interim as1}, namely
$$Z_{k+1}=Z_k+\epsilon_k\bm\psi_k'\hat{\mathbf d}_k+\frac{\epsilon_k^2}{2}\hat{\mathbf d}_k'\nabla^2Z(\mathbf p_k+\theta_k\epsilon_k\hat{\mathbf d}_k)\hat{\mathbf d}_k$$
for some $\tilde\theta_k$ between 0 and 1. Using the fact that $\bm\psi_k'\hat{\mathbf d}_k\geq\bm\psi_k'\mathbf d_k$ by the definition of $\mathbf d_k$, we have
\begin{align*}
Z_{k+1}&\geq Z_k+\epsilon_k\bm\psi_k'\mathbf d_k+\frac{\epsilon_k^2}{2}\hat{\mathbf d}_k'\nabla^2Z(\mathbf p_k+\theta_k\epsilon_k\hat{\mathbf d}_k)\hat{\mathbf d}_k\\
&=Z_k-\epsilon_kg_k+\frac{\epsilon_k^2}{2}\hat{\mathbf d}_k'\nabla^2Z(\mathbf p_k+\theta_k\epsilon_k\hat{\mathbf d}_k)\hat{\mathbf d}_k
\end{align*}
Now, using \eqref{main}, Condition \ref{c7} in Theorem \ref{rate thm} and $\|\hat{\mathbf d}_k\|^2\leq2$, we have
\begin{align}
Z_{k+1}&\geq Z_k-\epsilon_k\left(\frac{A}{k^C}+B\times\left\{\begin{array}{ll}\frac{1}{(C-\gamma)k^\gamma}&\text{if $0<\gamma<C$}\\
\frac{1}{(\gamma-C)(k_0-1)^{\gamma-C}k^C}&\text{if $\gamma>C$}\\
\frac{\log((k-1)/(k_0-1))}{k^C}&\text{if $\gamma=C$}
\end{array}\right\}\right)-\epsilon_k^2K\notag\\
&=Z_k-\frac{aA}{k^{1+C}}-aB\times\left\{\begin{array}{ll}\frac{1}{(C-\gamma)k^{1+\gamma}}&\text{if $0<\gamma<C$}\\
\frac{1}{(\gamma-C)(k_0-1)^{\gamma-C}k^{1+C}}&\text{if $\gamma>C$}\\
\frac{\log((k-1)/(k_0-1))}{k^{1+C}}&\text{if $\gamma=C$}
\end{array}\right\}-\frac{a^2K}{k^2}\label{interim rate12}
\end{align}
Now iterating \eqref{interim rate12} from $k$ to $l$, we have
$$Z_l\geq Z_k-\sum_{j=k}^{l-1}\frac{aA}{j^{1+C}}-aB\times\left\{\begin{array}{ll}\frac{1}{(C-\gamma)}\sum_{j=k}^{l-1}\frac{1}{j^{1+\gamma}}&\text{if $0<\gamma<C$}\\
\frac{1}{(\gamma-C)(k_0-1)^{\gamma-C}}\sum_{j=k}^{l-1}\frac{1}{j^{1+C}}&\text{if $\gamma>C$}\\
\sum_{j=k}^{l-1}\frac{\log((j-1)/(k_0-1))}{j^{1+C}}&\text{if $\gamma=C$}
\end{array}\right\}-a^2K\sum_{j=k}^{l-1}\frac{1}{j^2}$$
and letting $l\to\infty$, we get
\begin{equation}
Z^*\geq Z_k-\sum_{j=k}^\infty\frac{aA}{j^{1+C}}-aB\times\left\{\begin{array}{ll}\frac{1}{(C-\gamma)}\sum_{j=k}^\infty\frac{1}{j^{1+\gamma}}&\text{if $0<\gamma<C$}\\
\frac{1}{(\gamma-C)(k_0-1)^{\gamma-C}}\sum_{j=k}^\infty\frac{1}{j^{1+C}}&\text{if $\gamma>C$}\\
\sum_{j=k}^\infty\frac{\log((j-1)/(k_0-1))}{j^{1+C}}&\text{if $\gamma=C$}
\end{array}\right\}-a^2K\sum_{j=k}^\infty\frac{1}{j^2}\label{interim rate13}
\end{equation}
where the convergence to $Z^*$ is guaranteed by Theorem \ref{as}. Note that \eqref{interim rate13} implies that
\begin{align*}
Z^*&\geq Z_k-\frac{aA}{C(k-1)^C}-aB\times\left\{\begin{array}{ll}\frac{1}{(C-\gamma)\gamma(k-1)^\gamma}&\text{if $0<\gamma<C$}\\
\frac{1}{(\gamma-C)(k_0-1)^{\gamma-C}C(k-1)^C}&\text{if $\gamma>C$}\\
\frac{\log((k-1)/(k_0-1))}{C(k-1)^C}&\text{if $\gamma=C$}
\end{array}\right\}-\frac{a^2K}{k-1}\\
&\geq Z_k-\frac{D}{k-1}-\frac{E}{(k-1)^C}-F\times\left\{\begin{array}{ll}\frac{1}{(C-\gamma)\gamma(k-1)^\gamma}&\text{if $0<\gamma<C$}\\
\frac{1}{(\gamma-C)(k_0-1)^{\gamma-C}C(k-1)^C}&\text{if $\gamma>C$}\\
\frac{\log((k-1)/(k_0-1))}{C(k-1)^C}&\text{if $\gamma=C$}
\end{array}\right.
\end{align*}
where $D=a^2K$, $E=aA/C$ and $F=aB$. This gives \eqref{main1}.
\endproof

\proof{Proof of Lemma \ref{index max}}
Consider first a fixed $a$. When $a(1-\omega)>1$, \eqref{index} reduces to $\frac{\beta-\rho a-\zeta-2}{2(\beta+1)}\wedge\frac{1}{\beta+1}$. Since $\frac{\beta-\rho a-\zeta-2}{2(\beta+1)}$ is increasing in $\beta$ and $\frac{1}{\beta+1}$ is decreasing in $\beta$, the maximizer of $\frac{\beta-\rho a-\zeta-2}{2(\beta+1)}\wedge\frac{1}{\beta+1}$ occurs at the intersection of $\frac{\beta-\rho a-\zeta-2}{2(\beta+1)}$ and $\frac{1}{\beta+1}$, which is $\beta=\rho a+\zeta+4$. The associated value of \eqref{index} is $\frac{1}{\rho a+\zeta+5}$.

When $a(1-\omega)\leq1$, \eqref{index} reduces to $\frac{a(1-\omega)}{\beta+1}\wedge\frac{\beta-\rho a-\zeta-2}{2(\beta+1)}$. By a similar argument, the maximizer is $\beta=a(2-2\omega+\rho)+\zeta+2$, with the value of \eqref{index} equal to $\frac{a(1-\omega)}{a(2-2\omega+\rho)+\zeta+3}$.

Thus, overall, given $a$, the optimal choice of $\beta$ is $\beta=\rho a+\zeta+2+2((a(1-\omega))\wedge1)$, with the value of \eqref{index} given by $\frac{(a(1-\omega))\wedge1}{\rho a+\zeta+3+2((a(1-\omega))\wedge1)}$. When $a(1-\omega)>1$, the value of \eqref{index} is $\frac{1}{\rho a+\zeta+5}$ which is decreasing in $a$, whereas when $a(1-\omega)\leq1$, the value of \eqref{index} is $\frac{a(1-\omega)}{a(2-2\omega+\rho)+\zeta+3}$ which is increasing in $a$. Thus the maximum occurs when $a(1-\omega)=1$, or $a=\frac{1}{1-\omega}$. The associated value of \eqref{index} is $\frac{1}{\rho/(1-\omega)+\zeta+5}$.
\endproof

\remark
Suppose that Assumption \ref{bias} is replaced by letting
$$\|\mathbf v(\bm\xi_1)-\mathbf v(\bm\xi_2)\|\leq L\|\bm\xi_1-\bm\xi_2\|$$
hold for any $\bm\xi_1,\bm\xi_2\in\mathbb R^N$. Then, in the proof of Theorem \ref{rate thm}, the inequality \eqref{interim remark} can be replaced by
\begin{align*}
P(\mathcal E_k')&\leq\frac{k^\gamma E|(\hat{\bm\psi}_k-\bm\psi_k)'(\hat{\mathbf d}_k-\mathbf d_k)|}{\varrho}\\
&\leq\frac{k^\gamma}{\varrho}\sqrt{E[\|\hat{\bm\psi}_k-\bm\psi_k\|^2]E[\|\mathbf d_k-\hat{\mathbf d}_k\|^2]}\text{\ \ \ \ by the Cauchy-Schwarz inequality}\\
&\leq\frac{k^\gamma L}{\varrho}E[\|\hat{\bm\psi}_k-\bm\psi_k\|^2]\text{\ \ \ \ by the relaxed Assumption \ref{bias}}\\
&\leq\frac{LM^2TNk^\gamma}{R_k\varrho\delta}\prod_{j=1}^{k-1}(1-\epsilon_j)^{-1}\text{\ \ \ \ by following the derivation in \eqref{interim as2} and \eqref{interim as4}}
\end{align*}
Consequently, equation \eqref{interim remark1} becomes
$$\prod_{j=1}^{k_0-1}(1-\epsilon_j)^{-1}\frac{M^2TN}{\delta b}\left(\frac{1}{\vartheta^2(\beta-\rho a-1)(k_0-1)^{\beta-1}}+\frac{L}{\varrho(\beta-\gamma-\rho a-1)(k_0-1)^{\beta-\gamma-1}}\right)$$
if Condition \ref{c5} is replaced by
$$\beta>\gamma+\rho a+1$$
Correspondingly, Condition \ref{c6} needs to be replaced by
$$\prod_{j=1}^{k_0-1}(1-\epsilon_j)^{-1}\frac{M^2TN}{\delta b}\left(\frac{1}{\vartheta^2(\beta-\rho a-1)(k_0-1)^{\beta-1}}+\frac{L}{\varrho(\beta-\gamma-\rho a-1)(k_0-1)^{\beta-\gamma-1}}\right)<\varepsilon$$
The results in Theorem \ref{rate thm} and Corollary \ref{rate cor} then retain. Under these modified Conditions \ref{c5} and \ref{c6}, discussion point 3(b) in Section \ref{sec:local rate} then gives $\beta=\gamma+\rho a+1+\zeta$ for some $\zeta>0$ and $\gamma=\beta-\rho a-\zeta-1$. In discussion point 4, the convergence rate in terms of replications becomes $1/W^{((a(1-\omega))\wedge(\beta-\rho a-\zeta-1)\wedge1)/(\beta+1)}$. By maximizing
\begin{equation}
\frac{(a(1-\omega))\wedge(\beta-\rho a-\zeta-1)\wedge1}{\beta+1}\label{index1}
\end{equation}
like in \eqref{index} by Lemma \ref{index max} (see Lemma \ref{index max1} right after this remark), we get
$$a=\frac{1}{1-\omega},\ \ \beta=\frac{\rho}{1-\omega}+\zeta+2$$
and the optimal value is
$$\frac{1}{\rho/(1-\omega)+\zeta+3}$$
So, following the argument there, we choose $\vartheta$ and $\nu$, and hence $\omega$, to be small, and we choose $\rho$ to be close to 1. This gives rise to the approximate choice that $a\approx1+\omega$ and $\beta\approx3+\zeta+\omega$. The convergence rate is then $O(W^{-1/(4+\zeta+\omega)})$, leading to our claim in Section \ref{sec:local rate} that the complexity can improve to $O(1/\epsilon^{4+\zeta+\omega})$ if Assumption \ref{bias} is relaxed.
\label{remark:rate}
\endremark

\begin{lemma}
The maximizer of \eqref{index1} is given by
$$a=\frac{1}{1-\omega},\ \ \beta=\frac{\rho}{1-\omega}+\zeta+2$$
and the optimal value is
$$\frac{1}{\rho/(1-\omega)+\zeta+3}$$\label{index max1}
\end{lemma}

\proof{Proof of Lemma \ref{index max1}}
Consider first a fixed $a$. When $a(1-\omega)>1$, \eqref{index1} reduces to $\frac{\beta-\rho a-\zeta-1}{\beta+1}\wedge\frac{1}{\beta+1}$. Since $\frac{\beta-\rho a-\zeta-1}{\beta+1}$ is increasing in $\beta$ and $\frac{1}{\beta+1}$ is decreasing in $\beta$, the maximizer of $\frac{\beta-\rho a-\zeta-1}{\beta+1}\wedge\frac{1}{\beta+1}$ occurs at the intersection of $\frac{\beta-\rho a-\zeta-1}{\beta+1}$ and $\frac{1}{\beta+1}$, which is $\beta=\rho a+\zeta+2$. The associated value of \eqref{index1} is $\frac{1}{\rho a+\zeta+3}$.

When $a(1-\omega)\leq1$, \eqref{index1} reduces to $\frac{a(1-\omega)}{\beta+1}\wedge\frac{\beta-\rho a-\zeta-1}{\beta+1}$. By a similar argument, the maximizer is $\beta=a(1-\omega+\rho)+\zeta+1$, with the value of \eqref{index1} equal to $\frac{a(1-\omega)}{a(1-\omega+\rho)+\zeta+2}$.

Thus, overall, given $a$, the optimal choice of $\beta$ is $\beta=\rho a+\zeta+1+(a(1-\omega))\wedge1$, with the value of \eqref{index1} given by $\frac{(a(1-\omega))\wedge1}{\rho a+\zeta+2+(a(1-\omega))\wedge1}$. When $a(1-\omega)>1$, the value of \eqref{index1} is $\frac{1}{\rho a+\zeta+3}$ which is decreasing in $a$, whereas when $a(1-\omega)\leq1$, the value of \eqref{index} is $\frac{a(1-\omega)}{a(1-\omega+\rho)+\zeta+2}$ which is increasing in $a$. Thus the maximum occurs when $a(1-\omega)=1$, or $a=\frac{1}{1-\omega}$. The associated value of \eqref{index1} is $\frac{1}{\rho/(1-\omega)+\zeta+3}$.
\endproof

\section{Additional Details of the Numerical Results}
\subsection{Multi-start Initialization}\label{sec:multi-start}
\begin{figure}[!htbp]
  \centering
    \includegraphics[width=0.6\textwidth]{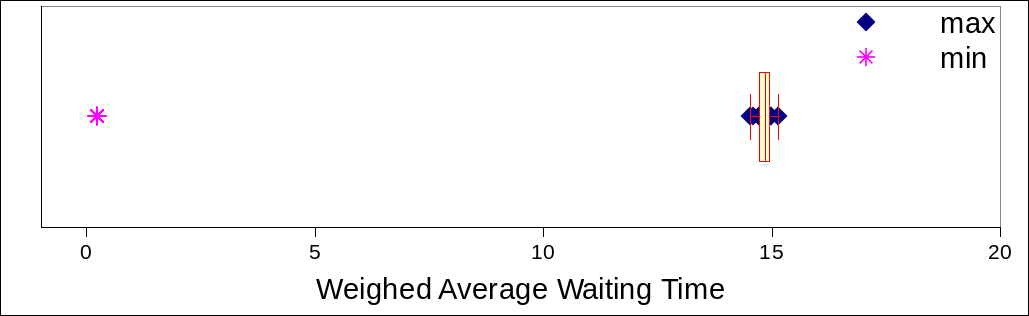}
    \caption{Returned optimal solutions from $10$ runs on $n=30$,
      $M=50$, exponential for discretization \label{fig:multirun}} 
\end{figure}

The results in Section \ref{expt:mltcls} are implemented with an initialization that assigns
equal probabilities to the support points. To test the procedure under different initializations, we repeat ten runs of the FWSA algorithm where the initial
probability masses for the support points (held constant for all runs)
are sampled uniformly independently with appropriate
normalization. Figure~\ref{fig:multirun} provides a box-plot of the identified optima. The sample size for moment constraint generation is
$N_s=50$ and the discretization support size is $n=30$. The returned
optimal solutions for each of the minimization and maximization
formulations all agree up to the first two digits (the box plot shows the small spread of the max values, while the min values are very clustered and they appear to all overlap at the same point). This indicates that the
formulations have a unique global optimal solution or similar local optimal solutions. Note that the bounds generated from this setting are quite loose with a small $N_s$.

\subsection{Details of the Benchmark Steady-State Formulation in Section \ref{sec:gg1example}}\label{sec:numerics appendix}
We consider the depicted $Z(\mathbf p)$ in Section \ref{sec:gg1example}. As $T$ grows, the
average waiting time converges to the corresponding steady-state
value, which, when the
traffic intensity $\rho_{\mbp} = E_{\mbp}[X_t]$ is less than $1$, is given in closed form by
the Pollaczek-Khinchine formula (\cite{klin32}) as:
\[
{Z}_{\infty}(\mbp)=\frac{\rho_{\mbp} E_{\mbp}[X_1] +
  Var_{\mbp}(X_1)}{2(1-\rho_{\mbp} )}.
\]
So, when $T$ is large, an approximation $Z^*_{\infty}$ to the
worst-case performance estimate can be obtained by replacing $Z(\mbp)$ with ${Z}_{\infty}(\mbp)$. (In
experiments, a choice of $T=500$ seems to show close agreement.) With
$E_{\mbp}[X_1]=\sum p_jy_j$ and $E_{\mbp}[X_1^2]=\sum p_jy_j^2$, the steady-state approximation to~(\ref{mg1prob}) is given by {\bf(SS)} below, which is equivalent to {\bf(SS$'$)} via variable substitutions (see
p.191 in~\cite{boyd2009convex}):

\noindent\begin{minipage}{0.45\textwidth}
  \begin{align}
  \min_{\mbp}  &\quad \frac { \sum_j p_jy_j^2} {2(1-
  \sum_j p_jy_j)}\quad\quad\quad \mbox{\bf (SS)} \nonumber\\
s.t.  &\quad \sum_j p_j \log \left(\frac {p_j}{p_{b,j}} \right) \le \eta\nonumber\\
&\quad \sum_j p_j = 1 \nonumber \\
&\quad 0 \le p_j\le 1,\quad \forall j=1,\ldots,n  \nonumber
  \end{align}
\end{minipage}
$\quad \Longrightarrow \quad$
\begin{minipage}{0.35\textwidth}
  \begin{align}
\min_{\mbp} & \quad \sum_j w_j {y_j}^{2}
\quad\quad\quad\quad\quad\quad \mbox{\bf (SS$'$)}
\nonumber\\ s.t.&\quad \sum_j w_j \log \left(\frac {w_j} {tp_{b,j}}\right)
\le \eta t \nonumber\\ & \quad 2t - 2 \sum_j w_jy_j =
1 \nonumber
\\ & \quad \sum_j w_j = t \nonumber
\\ & \quad
0\le w_j \le t \quad \forall j=1,\ldots,n \nonumber 
  \end{align}
\end{minipage}

\subsection{Shape of the Obtained Optimal Distributions in Section \ref{sec:gg1example}}\label{sec:numerics appendix1}
\begin{figure}[!htbp]
  \centering
  \begin{subfigure}[b]{0.45\textwidth}
    \includegraphics[width=\textwidth]{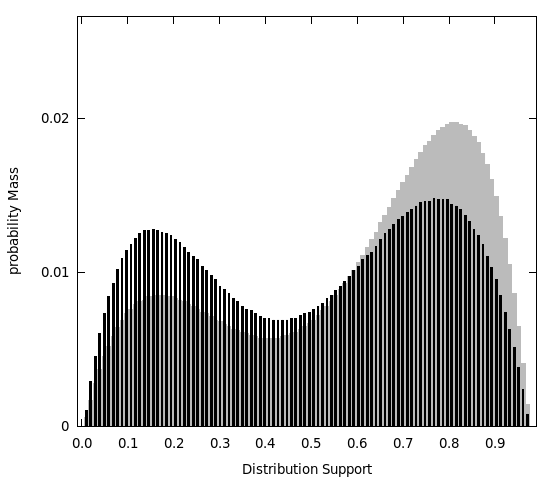}
    \caption{(min) $\mathbf p_b$ from beta-mixture}
    \label{fig:mg1optdistnminbeta}
  \end{subfigure}
\quad\quad
  \begin{subfigure}[b]{0.45\textwidth}
    \includegraphics[width=\textwidth]{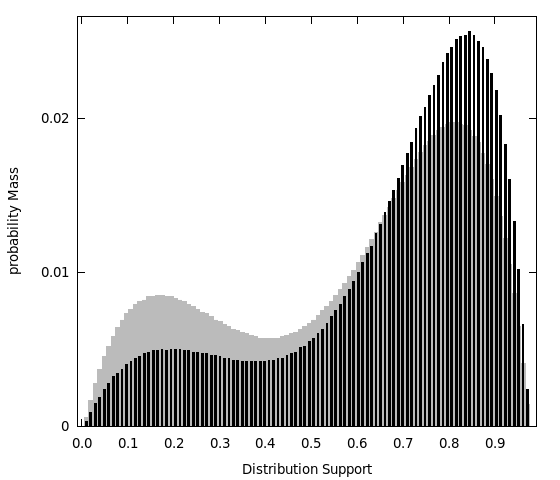}
    \caption{(max) $\mathbf p_b$ from beta-mixture}
    \label{fig:mg1optdistnmaxbeta}
  \end{subfigure}  
  \caption{Optimal solutions $\mbp^*$ identified by
    the FWSA algorithm with $n=100$ and $\eta=0.05$, setting
    $a=1.5,\beta=2.75$. The gray bars represent the baseline p.m.f. $\mbp_b$.
  }
  \label{fig:mg1worstdistn}
\end{figure}

Continuing with the example in Section \ref{sec:gg1example}, Figure~\ref{fig:mg1worstdistn} shows the form of the optimal
distributions $\mbp^*$ identified by the FWSA algorithm for the
minimization (Figure~\ref{fig:mg1optdistnminbeta}) and maximization
(Figure~\ref{fig:mg1optdistnmaxbeta}) problems under~\eqref{mg1prob}. The optimal distributions follow a similar bimodal
structure as the baseline distribution $\mathbf p_b$. The maximization version assigns probability masses in an
unequal manner to the two modes in order to drive up both the mean and the variance of $\mbp$, as {\bf (SS)} (in Appendix \ref{sec:numerics appendix}) leads us to expect, whereas the minimization
version on the other hand makes the mass allocation more equal in order
to minimize the mean and the variance of $\mbp$ while maintaining the
maximum allowed KL divergence.

\end{document}